%% file: main.tex
\documentclass[11pt,CJK]{article}
\usepackage{amsfonts,amsmath,amssymb,amsthm}
\usepackage{mathrsfs}
\usepackage{latexsym}
\usepackage{graphicx}
\usepackage{fancyhdr,cite}
\usepackage{indentfirst}
\usepackage{color}

\theoremstyle{plain}                       
\newtheorem{lemma}{Lemma}[section]
\newtheorem{thm}[lemma]{Theorem}
\newtheorem{cor}[lemma]{Corollary}
\newtheorem{remark}[lemma]{Remark}

\theoremstyle{remark}
\newtheorem*{pf}{Proof}

\numberwithin{equation}{section}

\def\Xint#1{\mathchoice
  {\XXint\displaystyle\textstyle{#1}}%
  {\XXint\textstyle\scriptstyle{#1}}%
  {\XXint\scriptstyle\scriptscriptstyle{#1}}%
  {\XXint\scriptscriptstyle\scriptscriptstyle{#1}}%
  \!\int}
\def\XXint#1#2#3{{\setbox0=\hbox{$#1{#2#3}{\int}$}
  \vcenter{\hbox{$#2#3$}}\kern-.5\wd0}}

\def\dashint{\Xint-}

\begin{document}

\include{article_1}

\clearpage
\end{document}

%% file: article_1.tex
\allowdisplaybreaks
\pagestyle{myheadings}\markboth{$~$ \hfill {\rm Q. Xu,} \hfill $~$} {$~$ \hfill {\rm  } \hfill$~$}
\author{Qiang Xu
\thanks{Email: xuqiang09@lzu.edu.cn.}
\thanks{This work was supported by the Chinese Scholar Council (File No. 201306180043).}
\\
School of Mathematics and Statistics, Lanzhou University, \\
Lanzhou, Gansu 730000, PR China.}

\title{\textbf{Uniform Regularity Estimates in Homogenization Theory of Elliptic System with Lower Order Terms} }
\maketitle
\begin{abstract}
  In this paper,
  we extend the uniform regularity estimates obtained by M. Avellaneda and F. Lin in \cite{MAFHL,MAFHL3}
  to the more general second order elliptic systems
  in divergence form $\{\mathcal{L}_\varepsilon, \varepsilon>0\}$, with rapidly oscillating periodic coefficients.
  We establish not only sharp $W^{1,p}$ estimates, H\"{o}lder estimates, Lipschitz estimates
  and non-tangential maximal function estimates for the Dirichlet problem on a bounded $C^{1,\eta}$ domain,
  but also a sharp $O(\varepsilon)$ convergence rate in $H_0^1(\Omega)$ by virtue of the Dirichlet correctors.
  Moreover, we define the Green's matrix associated with $\mathcal{L}_\varepsilon$ and obtain its decay estimates.
  We remark that the well known compactness methods are not employed here, instead we construct the transformations $\eqref{T:1}$
  to make full use of the results in \cite{MAFHL,MAFHL3}.
\end{abstract}

\section{Introduction and main results}
The main purpose of this paper is to study the uniform regularity estimates for second order elliptic systems with lower order terms,
arising in homogenization theory. More precisely, we consider
\begin{eqnarray*}
\mathcal{L}_{\varepsilon} = -\text{div}\left[A\left(x/\varepsilon\right)\nabla + V(x/\varepsilon)\right]
+ B(x/\varepsilon)\nabla + c(x/\varepsilon) + \lambda I,
\end{eqnarray*}
where $\lambda$ is a constant, and $I=(\delta^{\alpha\beta})$ denotes the identity matrix.
In a special case, let $A = I = 1$, $V=B$, $c=0$, and $\mathcal{W} = \text{div}(V)$, the operator $\mathcal{L}_\varepsilon$ becomes
\begin{equation*}
 \mathfrak{L}_\varepsilon = -\Delta +\frac{1}{\varepsilon}\mathcal{W}(x/\varepsilon) + \lambda,
\end{equation*}
where $\mathcal{W}$ is the rapidly oscillating potential term (see \cite[pp.91]{ABJLGP}).
It is not hard to see that the uniform
regularity estimates obtained in this paper are not trivial generalizations of \cite{MAFHL,MAFHL3},
and they are new even for $\mathfrak{L}_\varepsilon$.

Let $1 \leq i,j \leq d$, $1\leq \alpha,\beta\leq m$,
where $d \geq 3$ denotes the dimension, and $m\geq 1$ is the number of equations in the system.
Suppose that the measurable functions
$A = (a_{ij}^{\alpha\beta}):\mathbb{R}^d \rightarrow \mathbb{R}^{m^2\times d^2}$,
$V=(V_i^{\alpha\beta}):\mathbb{R}^d \rightarrow \mathbb{R}^{m^2\times d}$,
$B=(B_i^{\alpha\beta}):\mathbb{R}^d \rightarrow \mathbb{R}^{m^2\times d}$,
$c=(c^{\alpha\beta}):\mathbb{R}^d \rightarrow \mathbb{R}^{m^2}$
satisfy the following conditions:
\begin{itemize}
 \item the uniform ellipticity condition
\begin{equation}\label{a:1}
 \mu |\xi|^2 \leq a_{ij}^{\alpha\beta}(y)\xi_i^\alpha\xi_j^\beta\leq \mu^{-1} |\xi|^2
 \quad \text{for}~ y\in \mathbb{R}^d~\text{and}~ \xi=(\xi_i^\alpha)\in \mathbb{R}^{md}, ~\text{where}~ \mu>0;
\end{equation}
(The summation convention for repeated indices is used throughout.)
\item the periodicity condition
\begin{equation}\label{a:2}
A(y+z) = A(y), ~ V(y+z) = V(y),~B(y+z) = B(y),~ c(y+z) = c(y) \quad \text{for}~ y\in\mathbb{R}^d~\text{and}~ z\in\mathbb{Z}^d;
\end{equation}
\item the boundedness condition
\begin{equation}\label{a:3}
 \max\{\|V\|_{L^{\infty}(\mathbb{R}^d)},\|B\|_{L^{\infty}(\mathbb{R}^d)},\|c\|_{L^{\infty}(\mathbb{R}^d)}\}
 \leq \kappa_1, \quad \text{where}~\kappa_1>0;
\end{equation}
\item the regularity condition
\begin{equation}\label{a:4}
 \max\{ \|A\|_{C^{0,\tau}(\mathbb{R}^d)},~ \|V\|_{C^{0,\tau}(\mathbb{R}^d)},~\|B\|_{C^{0,\tau}(\mathbb{R}^d)}\} \leq \kappa_2,
 \quad \text{where}~\tau\in(0,1)~\text{and} ~\kappa_2 > 0.
\end{equation}
\end{itemize}
Set $\kappa = \max\{\kappa_1,\kappa_2\}$, and we say $A\in \Lambda(\mu,\tau,\kappa)$
if $A = A(y)$ satisfies conditions $\eqref{a:1}$, $\eqref{a:2}$ and $\eqref{a:4}$.
Throughout this paper, we always assume that $\Omega$ is a bounded $C^{1,\eta}$ domain with $\eta\in [\tau,1)$,
and $L_\varepsilon = - \text{div}[A(x/\varepsilon)\nabla]$ is the elliptic operator from \cite{MAFHL,MAFHL3},
unless otherwise stated.

The main idea of this paper is to find
the transformations $\eqref{T:1}$ between two solutions corresponding to $\mathcal{L}_\varepsilon$ and $L_\varepsilon$ such that the regularity results of
$L_\varepsilon$ can be applied to $\mathcal{L}_\varepsilon$ directly. Particularly,
to handle the boundary Lipschitz estimates, we define the Dirichlet correctors 
$\Phi_{\varepsilon,k} = (\Phi_{\varepsilon,k}^{\alpha\beta})$, $0\leq k\leq d$, associated with $\mathcal{L}_\varepsilon$ as follows:
\begin{equation}\label{pde:1.2}
 L_\varepsilon(\Phi_{\varepsilon,k}) = \text{div}(V_\varepsilon) ~~\text{in}~\Omega, \qquad \Phi_{\varepsilon,k} = I ~~\text{on}~\partial\Omega
\end{equation}
for $k=0$, and
\begin{equation}\label{pde:1.3}
 L_\varepsilon(\Phi_{\varepsilon,k}^\beta) = 0 ~~\text{in}~\Omega, \qquad \Phi_{\varepsilon,k}^\beta = P_k^\beta ~~\text{on}~\partial\Omega
\end{equation}
for $1\leq k\leq d$, where $V_\varepsilon(x) = V(x/\varepsilon)$,
$\Phi_{\varepsilon,k}^\beta = (\Phi_{\varepsilon,k}^{1\beta},\cdots,\Phi_{\varepsilon,k}^{m\beta})\in H^1(\Omega;\mathbb{R}^m)$,
and $P_k^\beta = x_k(0,\cdots,1,\cdots,0)$ with $1$ in the $\beta^{th}$ position.
We remark that $\eqref{pde:1.3}$
was studied in \cite{MAFHL,SZW0}, but $\eqref{pde:1.2}$ has not yet been developed.
Here we show that $\Phi_{\varepsilon,0}$
ought to be of the form in $\eqref{pde:1.2}$, and its properties are shown in section 4.

For Neumann boundary conditions, a significant development
was made by C.E. Kenig, F. Lin and Z. Shen \cite{SZW4}, where they constructed Neumann correctors to verify the Lipschitz estimates
for $L_\varepsilon$. Recently, S.N. Armstrong and Z. Shen \cite{SZ} found a new way to obtain the same results
even without Dirichlet correctors or Neumann correctors in the almost periodic setting.
We plan to study uniform regularity estimates for $\mathcal{L}_\varepsilon$ with Neumann boundary conditions in a forthcoming paper.

The main results are as follows.


\begin{thm}[$W^{1,p}$ estimates]\label{thm:1.1}
Suppose that $A\in \emph{VMO}(\mathbb{R}^d)$ satisfies $\eqref{a:1}$, $\eqref{a:2}$, and other coefficients of $\mathcal{L}_\varepsilon$ satisfy $\eqref{a:3}$.
Let $1<p<\infty$, $f = (f_i^\alpha)\in L^p(\Omega;\mathbb{R}^{md})$, $F\in L^q(\Omega;\mathbb{R}^m)$ and $g\in B^{1-\frac{1}{p},p}(\partial\Omega;\mathbb{R}^m)$,
where $q = \frac{pd}{p+d}$ if $p > \frac{d}{d-1}$, and $q > 1$ if $1<p\leq\frac{d}{d-1}$.
Then the Dirichlet problem
\begin{equation}\label{pde:1.1}
 \left\{\begin{aligned}
 \mathcal{L}_\varepsilon(u_\varepsilon) &=& \emph{div}(f) + F  \quad &\emph{in} ~~\Omega, \\
  u_\varepsilon &=&g \qquad\qquad\quad &\emph{on}~~\partial\Omega
 \end{aligned}\right.
\end{equation}
has a unique weak solution $u_\varepsilon \in W^{1,p}(\Omega;\mathbb{R}^m)$,
whenever $\lambda \geq \lambda_0$ and $\lambda_0=\lambda_0(\mu,\kappa,m,d)$ is sufficiently large.
Furthermore, the solution satisfies the uniform estimate
\begin{equation}\label{pri:1.1}
\|\nabla u_\varepsilon\|_{L^p(\Omega)} \leq C\big\{\|f\|_{L^p(\Omega)} + \|F\|_{L^q(\Omega)} + \|g\|_{B^{1-1/p,p}(\partial\Omega)}\big\},
\end{equation}
where $C$ depends on $\mu,\omega(t),\kappa,\lambda, p, q,d,m$ and $\Omega$.
\end{thm}

Note that $A\in \text{VMO}(\mathbb{R}^d)$ if $A$ satisfies
\begin{equation*}
 \sup_{x\in \mathbb{R}^d\atop 0<\rho< t}\dashint_{B(x,\rho)}\Big|A(y) - \dashint_{B(x,\rho)}A\Big| dy
 \leq \omega(t), \quad \text{and} \quad \lim_{t\to0^+} \omega(t) = 0,
\end{equation*}
and $B^{\alpha,p}(\partial\Omega;\mathbb{R}^m)$ denotes the $L^p$ Besov space of order $\alpha$ (see \cite{RAA}).
We mention that for ease of notations we say the constant $C$ depends on $\omega$ instead of $\omega(t)$ in the rest of the paper.
We will prove Theorem $\ref{thm:1.1}$ in Section 3 by using bootstrap and duality arguments.
We mention that there are no periodicity or regularity assumptions on the coefficients of the lower order terms in Theorems $\ref{thm:1.1}$ and
$\ref{thm:1.4}$. The estimate $\eqref{pri:1.1}$ still holds when $\Omega$ is a bounded $C^1$ domain (see \cite{SZW10}).

Results of the $W^{1,p}$ estimates for elliptic or parabolic equations with VMO coefficients
can be found in \cite{SBLW,SBLW1,LI,GJ,DKNVK,NVK}.
In the periodic setting,
similar estimates for parabolic systems, elasticity systems, and Stokes systems were obtained by \cite{GZ,JGZSLS,SGZWS}, respectively.
Also, the uniform $W^{1,p}$ estimates for $L_\varepsilon$ with almost periodic coefficients were shown in \cite{SZ} recently.

\begin{thm}[H\"{o}lder estimates]\label{thm:1.4}
Suppose that the coefficients of $\mathcal{L}_\varepsilon$ satisfy the same conditions as in Theorem $\ref{thm:1.1}$.
Let $p>d$, $f=(f^\alpha_i)\in L^p(\Omega;\mathbb{R}^{md})$, $F\in L^q(\Omega;\mathbb{R}^m)$, and $g\in C^{0,\sigma}(\partial\Omega;\mathbb{R}^m)$,
where $q=\frac{pd}{p+d}$, and $\sigma = 1-d/p$. Then the weak solution $u_\varepsilon$ to $\eqref{pde:1.1}$
satisfies the uniform estimate
\begin{equation}\label{pri:1.5}
 \|u_\varepsilon\|_{C^{0,\sigma}(\Omega)} \leq C\big\{\|f\|_{L^p(\Omega)} + \|F\|_{L^q(\Omega)} + \|g\|_{C^{0,\sigma}(\partial\Omega)} \big\},
\end{equation}
where $C$ depends on $\mu,\omega,\kappa,\lambda, p, \sigma, d, m$ and $\Omega$.
\end{thm}

The estimate $\eqref{pri:1.5}$ is sharp in terms of the H\"older exponent of $g$.
If $g\in C^{0,1}(\partial\Omega;\mathbb{R}^m)$, $\eqref{pri:1.5}$ is just Corollary $\ref{cor:4.1}$.
The uniform H\"{o}lder estimates for $L_\varepsilon$ were given in \cite{MAFHL} by the compactness method which also works
for non-divergence form elliptic equations (see \cite{MAFHL2}).
However, we can not derive the sharp estimate by simply applying this method.
So we turn to study the Green's matrix $\mathcal{G}_\varepsilon(x,y)$ associated with $\mathcal{L}_\varepsilon$ and obtain the decay estimates
\begin{eqnarray*}
 |\mathcal{G}_\varepsilon(x,y)| \leq \frac{C}{|x-y|^{d-2}}\min\Big\{1,\frac{d^\sigma_x}{|x-y|^\sigma},
 \frac{d^{\sigma^\prime}_y}{|x-y|^{\sigma^\prime}},\frac{d^\sigma_xd^{\sigma^\prime}_y}{|x-y|^{\sigma+\sigma^\prime}}\Big\}, \qquad \forall~ x,y\in \Omega,~~x\not=y,
\end{eqnarray*}
where $\sigma,\sigma^\prime\in(0,1)$, $d_x=\text{dist}(x,\partial\Omega)$ denotes the distance between $x$ and $\partial\Omega$, and $C$ is independent of $\varepsilon$ (see Theorem $\ref{thm:4.2}$).
Then we prove Theorem $\ref{thm:1.4}$ through a subtle argument developed by Z. Shen \cite{SZW8},
where he proved a similar result for $L_\varepsilon$ in the almost periodic setting.

The existence and some related properties of the Green's matrix with respect to $L_1$ were studied by S. Hofmann
and S. Kim \cite{SHSK}.
We also refer the reader to \cite{SHSK2} for parabolic systems, and \cite{MGKOW,WLGSHFW} for the scalar case.

\begin{thm}[Lipschitz estimates]\label{thm:1.2} 
 Suppose that $A\in\Lambda(\mu,\tau,\kappa)$, $V$ satisfies $\eqref{a:2}$, $\eqref{a:4}$, $B$ and $c$ satisfy $\eqref{a:3}$, and $\lambda\geq\lambda_0$.
Let $p>d$ and $0<\sigma\leq\eta$. Then for any $f\in C^{0,\sigma}(\Omega;\mathbb{R}^{md})$, $F\in L^p(\Omega;\mathbb{R}^m)$, and
$g\in C^{1,\sigma}(\partial\Omega;\mathbb{R}^m)$, the weak solution to $\eqref{pde:1.1}$
satisfies the uniform estimate
\begin{equation}\label{pri:1.2}
 \|\nabla u_\varepsilon\|_{L^\infty(\Omega)}
 \leq C \left\{\|f\|_{C^{0,\sigma}(\Omega)}+\|F\|_{L^p(\Omega)}+\|g\|_{C^{1,\sigma}(\partial\Omega)}\right\},
\end{equation}
where $C$ depends on $\mu,\tau,\kappa,\lambda,p,d,m,\sigma,\eta$ and $\Omega$.
\end{thm}
The estimate $\eqref{pri:1.2}$ can not be improved
even if the coefficients of $\mathcal{L}_\varepsilon$ and $\Omega$ are smooth, since
the corrector $\chi_0$ defined in $\eqref{pde:2.1}$ is a counter example.
Here we use two important transformations
\begin{equation}\label{T:1}
  u_\varepsilon = \big[I+\varepsilon\chi_0(x/\varepsilon)\big]v_\varepsilon,
 \qquad\quad \text{and} \qquad  u_\varepsilon = \Phi_{\varepsilon,0} v_\varepsilon  
\end{equation}
to deal with the interior and global Lipschitz estimates, respectively.
We explain the main idea as follows:
\begin{equation*}
 (\text{D}_1)\left\{
 \begin{aligned}
  \mathcal{L}_\varepsilon(u_\varepsilon) & = \text{div}(f) + F& \quad  & \text{in}~\Omega, \\
                          u_\varepsilon  & = g         & \quad  & \text{on}~\partial\Omega.
 \end{aligned}\right.
 \qquad \underrightarrow{~ u_\varepsilon = \Phi_{\varepsilon,0}v_\varepsilon~} \qquad
  (\text{D}_2)\left\{
 \begin{aligned}
  L_\varepsilon(v_\varepsilon) & = \text{div}(\tilde{f}) + \tilde{F} &\quad  & \text{in}~ \Omega, \\
                          v_\varepsilon  & = g                &\quad  & \text{on}~\partial\Omega.
 \end{aligned}\right.
\end{equation*}
Note that
$\Phi_{\varepsilon,0}$ is not periodic, which is the main difficulty to overcome. So we
rewrite $(\text{D}_1)$ as $(\text{D}_2)$ to keep $L_\varepsilon$ periodic, while the price to pay is that the new source term $\tilde{f}$
involves $\nabla u_\varepsilon$. As we mentioned before, there is no uniformly bounded H\"{o}lder estimate for $\nabla u_\varepsilon$.
Fortunately, it follows from Theorem $\ref{thm:1.4}$ that
\begin{equation}\label{f:1.1}
\|\nabla\Phi_{\varepsilon,0}\|_{C^{0,\sigma_1}(\Omega)} = O(\varepsilon^{-\sigma_1})
\quad \text{and}\quad \|\nabla u_\varepsilon\|_{C^{0,\sigma_1}(\Omega)} = O(\varepsilon^{-\sigma_2})
\quad \text{as}~ \varepsilon \to 0,
\end{equation}
where $0< \sigma_1 < \sigma_2 < 1$ are independent of $\varepsilon$ (see Lemma $\ref{cor:5.1}$ and $\ref{lemma:5.7}$).
Together with an important consequence of Lemma $\ref{lemma:5.3}$
\begin{equation}\label{f:1.2}
\|\Phi_{\varepsilon,0} - I\|_{L^\infty(\Omega)} = O(\varepsilon)
\quad \text{as} ~ \varepsilon\to 0,
\end{equation}
we obtain that $\tilde{f}$ is uniformly H\"{o}lder continuous through the observation that the convergence rate in $\eqref{f:1.2}$ is faster than the divergence rate in
$\eqref{f:1.1}$ as $\varepsilon\to 0$.
Also, Theorem $\ref{thm:1.1}$ implies $\tilde{F}\in L^p(\Omega)$ with $p>d$. Thus we can employ the results in \cite{MAFHL} immediately,
and the proof of Theorem $\ref{thm:1.2}$ is finalized by a suitable extension technique.

We remark that
the compactness argument for our elliptic systems is also valid,
however it would be much more complicated.
For more references,
C. E. Kenig and C. Prange \cite{CC} established uniform Lipschtiz estimates with more general source terms in the oscillating boundaries setting,
and the same type of results for parabolic systems and Stokes systems were shown in \cite{GZ,SGZWS}, respectively.

\begin{thm}[Nontangential maximal function estimates]\label{thm:1.5}
Suppose that
$ A\in \Lambda(\mu,\tau,\kappa)$, $V,B$ satisfy $\eqref{a:2}$ and $\eqref{a:4}$, $c$ satisfies $\eqref{a:3}$, and $\lambda\geq\lambda_0$.
Let $1<p<\infty$, and $u_\varepsilon$
be the solution of the $L^p$ Dirichlet problem $\mathcal{L}_\varepsilon(u_\varepsilon) = 0$ in $\Omega$
and $u_\varepsilon = g$ on $\partial\Omega$
with $(u_\varepsilon)^*\in L^p(\partial\Omega)$, where $g\in L^p(\partial\Omega;\mathbb{R}^m)$
and $(u_\varepsilon)^*$ is the nontangential maximal function. Then
\begin{equation}\label{pri:1.6}
 \|(u_\varepsilon)^*\|_{L^p(\partial\Omega)} \leq C_p\|g\|_{L^p(\partial\Omega)},
\end{equation}
where $C_p$ depends on $\mu,\tau,\kappa,\lambda,d,m,p,\eta$ and $\Omega$. Furthermore, if $g\in L^\infty(\partial\Omega;\mathbb{R}^m)$, we have
\begin{equation}\label{pri:1.7}
\|u_\varepsilon\|_{L^\infty(\Omega)} \leq C \|g\|_{L^\infty(\partial\Omega)},
\end{equation}
where $C$ depends on $\mu,\tau,\kappa,\lambda,d,m,\eta$ and $\Omega$.
\end{thm}

The estimate $\eqref{pri:1.7}$ is known as the Agmon-Miranda maximum principle,
and $(u_\varepsilon)^*$ is defined in $\eqref{def:3}$.
We remark that the proof of Theorem $\ref{thm:1.5}$ is motivated by \cite{MAFHL,SZW4,SZW0}.
Define the Poisson kernel associated with $\mathcal{L}_\varepsilon$ as
\begin{equation*}
  \mathcal{P}_\varepsilon^{\gamma\beta}(x,y)
 = -n_j(x)a_{ij}^{\alpha\beta}(x/\varepsilon)\frac{\partial}{\partial x_i}\big\{\mathcal{G}_\varepsilon^{\alpha\gamma}(x,y)\big\}
  -n_j(x)B_{j}^{\alpha\beta}(x/\varepsilon)\mathcal{G}_\varepsilon^{\alpha\gamma}(x,y),
\end{equation*}
where $n_j$ denotes the $j^{th}$ component of the outward unit normal vector of $\partial\Omega$.
Due to Theorem $\ref{thm:1.2}$, we obtain
$|\nabla_x\nabla_y\mathcal{G}_\varepsilon(x,y)| \leq C|x-y|^{-d}$ for $x,y\in\Omega$, and $x\not=y$
(see Lemma $\ref{lemma:5.5}$), which implies the decay estimate of $\mathcal{P}_\varepsilon$ (see $\eqref{f:5.18}$).
Thus the solution $u_\varepsilon$ can be formulated by $\eqref{f:5.11}$.
Note that $\mathcal{P}_\varepsilon$ is actually
closely related to the adjoint operator $\mathcal{L}_\varepsilon^*$ (see Remark $\ref{rm:2.3}$).
That is the reason why we additionally assume $\eqref{a:2}$ and $\eqref{a:4}$ for $B$ in this theorem.
We refer the reader to Remark $\ref{rm:5.1}$ for more references on Theorem $\ref{thm:1.5}$.


\begin{thm}[Convergence rates]\label{thm:1.3}
 Let $\Omega$ be a bounded $C^{1,1}$ domain. Suppose that the coefficients of $\mathcal{L}_\varepsilon$ satisfy the same conditions as in
 Theorem $\ref{thm:1.2}$, and $B, c$ additionally satisfy the periodicity condition $\eqref{a:2}$.
 Let $u_\varepsilon$ be the weak solution to $\mathcal{L}_\varepsilon(u_\varepsilon) = F$ in $\Omega$ and $u_\varepsilon = 0$ on $\partial\Omega$,
 where $F\in L^{2}(\Omega;\mathbb{R}^m)$.
 Then we have
 \begin{equation}\label{pri:1.3}
\big\|u_\varepsilon - \Phi_{\varepsilon,0}u - (\Phi_{\varepsilon,k}^\beta - P_k^\beta)\frac{\partial u^\beta}{\partial x_k}\big\|_{H^{1}_0(\Omega)}
\leq C\varepsilon\|F\|_{L^2(\Omega)},
\end{equation}
 where $u$ satisfies $\mathcal{L}_0(u) = F$ in $\Omega$ and $u= 0$ on $\partial\Omega$.
 Moreover, if the coefficients of $\mathcal{L}_\varepsilon$ satisfy $\eqref{a:1}-\eqref{a:4}$, then
 \begin{equation}\label{pri:1.4}
  \|u_\varepsilon - u\|_{L^q(\Omega)} \leq C\varepsilon \|F\|_{L^{p}(\Omega)}
  \qquad
 \end{equation}
 holds for any $F\in L^p(\Omega;\mathbb{R}^m)$, where $q = \frac{pd}{d-p}$ if $1<p<d$, $q = \infty$ if $p>d$, and $C$ depends on $\mu,\tau, \kappa, \lambda, m,d, p$ and $\Omega$.
\end{thm}

 We mention that the estimates $\eqref{pri:1.3}$ and $\eqref{pri:1.4}$ are sharp in terms of the order of $\varepsilon$.
 The ideas in the proof are mainly inspired by \cite{SZW2,SZW1}. It is easy to see $\|u_\varepsilon - u\|_{L^2(\Omega)}= O(\varepsilon)$
 is a direct corollary of $\eqref{pri:1.3}$ or $\eqref{pri:1.4}$.
 In the case of $p=d$, $\eqref{pri:1.4}$ is shown in Remark $\ref{rm:6.1}$.
 Moreover, in the sense of ``operator error estimates''
 the convergence rate like $\eqref{pri:1.4}$ can also be expressed by
 $\|\mathcal{L}_\varepsilon^{-1}-\mathcal{L}_0^{-1}\|_{L^p\to L^q} \leq C\varepsilon$,
 where $\|\cdot\|_{L^p\to L^q}$ is referred to as the $(L^p\to L^q)$-operator norm.

 The convergence rates are active topics in homogenization theory. Decades ago,
 the $L^2$ convergence rates were obtained in \cite{ABJLGP,VSO}
 for scalar cases due to the maximum principle. At the beginning of 2000's,
 the operator-theoretic (spectral) approach was successfully introduced by M. Sh. Birman and T. A. Suslina \cite{BMSHSTA,BMSHSTA2}
 to investigate the convergence rates (operator error estimates) for the problems in the whole space $\mathbb{R}^d$.
 They obtained the sharp convergence rates $O(\varepsilon)$ in the
 $(L^2\to L^2)$-operator norm and $(L^2\to H^1)$-operator norm for a wide class of matrix strongly elliptic second order self-adjoint operators, respectively.
 These results were extended to second order strongly elliptic systems including lower order terms in \cite{SAT2}.
 Recently, C. E. Kenig, F. Lin and Z. Shen \cite{SZW2}
 developed the $L^2$ convergence rates for elliptic systems on Lipschitz domains with either Dirichlet or Neumann boundary data
 by additionally assuming
 regularity and symmetry conditions, while T.A. Suslina \cite{SAT,SAT1} also obtained similar results on a bounded $C^{1,1}$ domain without
 any regularity assumption on the coefficients.
 We refer the reader to \cite{CDJ,GG1,GG2,ODVB,ZVVPSE} and references therein for more results.

 In the end, we comment that the above five theorems are still true for $d=1,2$.
 Since we usually have a different method to treat the cases $d\geq 3$ and $d=1, 2$
 (for example, see \cite[pp.2]{SHSK}), we omit the discussion about the cases of $d=1,2$ here.

\section{Preliminaries}
Define the correctors $\chi_k = (\chi_{k}^{\alpha\beta})$, $0\leq k\leq d$, associated with $\mathcal{L}_\varepsilon$ as follows:
\begin{equation}\label{pde:2.1}
\left\{ \begin{aligned}
 &L_1(\chi_k) = \text{div}(V)  \quad \text{in}~ \mathbb{R}^d, \\
 &\chi_k\in H^1_{per}(Y;\mathbb{R}^{m^2})~~\text{and}~\int_Y\chi_k dy = 0
\end{aligned}
\right.
\end{equation}
for $k=0$, and
\begin{equation}\label{pde:2.2}
 \left\{ \begin{aligned}
  &L_1(\chi_k^\beta + P_k^\beta) = 0 \quad \text{in}~ \mathbb{R}^d~, \\
  &\chi_k^\beta \in H^1_{per}(Y;\mathbb{R}^m)~~\text{and}~\int_Y\chi_k^\beta dy = 0
 \end{aligned}
 \right.
\end{equation}
for $1\leq k\leq d$, where $Y = [0,1)^d \cong \mathbb{R}^d/\mathbb{Z}^d$, and $H^1_{per}(Y;\mathbb{R}^m)$ denotes the closure
of $C^\infty_{per}(Y;\mathbb{R}^m)$ in $H^1(Y;\mathbb{R}^m)$.
Note that $C^\infty_{per}(Y;\mathbb{R}^m)$ is the subset of $C^\infty(Y;\mathbb{R}^m)$, which collects all $Y$-periodic vector-valued functions
(see \cite[pp.56]{ACPD}). By asymptotic expansion arguments, we obtain the homogenized operator
\begin{equation*}
 \mathcal{L}_0 = -\text{div}(\widehat{A}\nabla+ \widehat{V}) + \widehat{B}\nabla + \widehat{c} + \lambda I,
\end{equation*}
where $\widehat{A} = (\hat{a}_{ij}^{\alpha\beta})$, $\widehat{V}=(\widehat{V}_i^{\alpha\beta})$,
$\widehat{B} = (\widehat{B}_i^{\alpha\beta})$ and $\widehat{c}= (\hat{c}^{\alpha\beta})$ are given by
\begin{equation}\label{f:2.2}
\begin{aligned}
\hat{a}_{ij}^{\alpha\beta} = \int_Y \big[a_{ij}^{\alpha\beta} + a_{ik}^{\alpha\gamma}\frac{\partial\chi_j^{\gamma\beta}}{\partial y_k}\big] dy, \qquad
\widehat{V}_i^{\alpha\beta} = \int_Y \big[V_i^{\alpha\beta} + a_{ij}^{\alpha\gamma}\frac{\partial\chi_0^{\gamma\beta}}{\partial y_j}\big] dy, \\
\widehat{B}_i^{\alpha\beta} = \int_Y \big[B_i^{\alpha\beta} + B_j^{\alpha\gamma}\frac{\partial\chi_i^{\gamma\beta}}{\partial y_j}\big] dy, \qquad
\hat{c}^{\alpha\beta} = \int_Y \big[c^{\alpha\beta} + B_i^{\alpha\gamma}\frac{\partial\chi_0^{\gamma\beta}}{\partial y_i}\big] dy.
\end{aligned}
\end{equation}
The proof is left to readers (see \cite[pp.103]{ABJLGP} or \cite[pp.31]{VSO}).

\begin{remark}\label{rm:2.1}
\emph{ It follows from the conditions $\eqref{a:1}$ and $\eqref{a:2}$ that $\mu|\xi|^2\leq \hat{a}_{ij}^{\alpha\beta}\xi_i^\alpha\xi_j^\beta\leq \mu_1|\xi|^2$ holds
for any $\xi=(\xi_i^\alpha)\in \mathbb{R}^{md}$,
where $\mu_1$ depends only on $\mu$. Moreover, if $a_{ij}^{\alpha\beta}=a_{ji}^{\beta\alpha}$,
then $\mu|\xi|^2\leq \hat{a}_{ij}^{\alpha\beta}\xi_i^\alpha\xi_j^\beta\leq \mu^{-1}|\xi|^2$ (see \cite[pp.23]{ABJLGP}).
This illustrates that the operator $\mathcal{L}_0$ is still elliptic.}
\end{remark}

\begin{remark}\label{rm:2.6}
\emph{ We introduce the following notations for simplicity.
We write $\chi_{k,\varepsilon}(x) = \chi_k(x/\varepsilon)$, $A_\varepsilon(x) = A(x/\varepsilon)$, $V_\varepsilon(x) = V(x/\varepsilon)$,
$B_\varepsilon(x) = B(x/\varepsilon)$, $c_\varepsilon(x) = c(x/\varepsilon)$, and their components follow the same abbreviated way.
Note that
the abbreviations are not applied to $\Phi_{\varepsilon,k}(x)$ or $\Psi_{\varepsilon,k}(x)$.}

\emph{Let $B = B(x,r) = B_r(x)$, and $KB = B(x,Kr)$ denote the concentric balls as $K>0$ varies, where $r<1$ in general.
We say that $\Omega$ is a bounded $C^{1,\eta}$ domain, if there exist
$r_0>0$, $M_0>0$ and $\{P_i:i=1,2,\cdots,n_0\}\subset\partial\Omega$
such that $\partial\Omega\subset\cup_{i=1}^{n_0}B(P_i,r_0)$ and for each $i$, there
exists a function $\psi_i\in C^{1,\eta}(\mathbb{R}^{d-1})$ and a coordinate system, such that
$B(P_i, C_0r_0)\cap\Omega = B(0, C_0r_0)
\cap\{(x^\prime,x_d)\in\mathbb{R}^d:x^\prime\in\mathbb{R}^{d-1} ~\text{and}~x_d>\psi_i(x^\prime)\}$,
where $C_0 = 10(M_0+1)$ and $\psi_i$ satisfies
\begin{equation}
 \psi_i(0) = 0, \quad \text{and} \quad \|\psi_i\|_{C^{1,\eta}(\mathbb{R}^{d-1})} \leq M_0.
\end{equation}
We set
$D(r)=D(r,\psi) = \{(x^\prime,x_d)\in\mathbb{R}^d:|x^\prime|<r~\text{and}~\psi(x^\prime)< x_d < \psi(x^\prime) + C_0r\}$ and
$\Delta(r) = \Delta(r,\psi) = \{(x^\prime,\psi(x^\prime))\in\mathbb{R}^d:|x^\prime|<r\}$.
In the paper, we say the constant $C$ depends on $\Omega$,
which means $C$ involves both $M_0$ and $|\Omega|$. This is especially important when we do near boundary regularity estimates.
Here $|\Omega|$ denotes the volume of $\Omega$. We also mention that
for any $E\subset \Omega$, we write $\bar{f}_E = \dashint_{E} f(x) dx = \frac{1}{|E|}\int_E f(x) dx$,
and the subscript of $\bar{f}_E$ is usually omitted.
}
\end{remark}

\begin{remark}\label{rm:2.3}
\emph{Let $\mathcal{L}_\varepsilon^*$ be the adjoint of $\mathcal{L}_\varepsilon$, given by
\begin{equation*}
 \big[\mathcal{L}_\varepsilon^*(v_\varepsilon)\big]^{\beta}
 = -\frac{\partial}{\partial x_j}\Big\{a_{ij}^{\alpha\beta}(x/\varepsilon)\frac{\partial v_\varepsilon^\alpha}{\partial x_i}
 + B^{\alpha\beta}_{j}(x/\varepsilon)v_\varepsilon^\alpha \Big\}
 + V_{i}^{\alpha\beta}(x/\varepsilon)\frac{\partial v_\varepsilon^\alpha}{\partial x_i}
 + c^{\alpha\beta}(x/\varepsilon) v_\varepsilon^\alpha + \lambda v_\varepsilon^\beta.
\end{equation*}
Then  we define
the bilinear form associated with $\mathcal{L}_{\varepsilon}$ as
\begin{eqnarray*}
 \mathrm{B}_\varepsilon[u_\varepsilon,\phi]
 = \int_\Omega \Big\{a_{ij,\varepsilon}^{\alpha\beta}\frac{\partial u_\varepsilon^\beta}{\partial x_j} + V_{i,\varepsilon}^{\alpha\beta} u_\varepsilon^\beta\Big\}
 \frac{\partial \phi^\alpha}{\partial x_i} dx
 + \int_\Omega \Big\{B_{i,\varepsilon}^{\alpha\beta}\frac{\partial u_\varepsilon^\beta}{\partial x_i}
  + c_\varepsilon^{\alpha\beta} u_\varepsilon^\beta + \lambda u_\varepsilon^\alpha \Big\} \phi^\alpha dx,
\end{eqnarray*}
and the conjugate bilinear form with respect to $\mathcal{L}_\varepsilon^*$ as
\begin{eqnarray}\label{f:2.1}
 \mathrm{B}_\varepsilon^*[v_\varepsilon,\phi]
 = \int_\Omega \Big\{a_{ij,\varepsilon}^{\alpha\beta}\frac{\partial v_\varepsilon^\alpha}{\partial x_i}
 + B_{j,\varepsilon}^{\alpha\beta} v_\varepsilon^\alpha\Big\} \frac{\partial \phi^\beta}{\partial x_j} dx
 + \int_\Omega \Big\{V_{i,\varepsilon}^{\alpha\beta}\frac{\partial v_\varepsilon^\beta}{\partial x_i}
  + c_\varepsilon^{\alpha\beta} v_\varepsilon^\alpha + \lambda v_\varepsilon^\beta \Big\} \phi^\beta dx
\end{eqnarray}
for any $u_\varepsilon, v_\varepsilon, \phi \in H^1_0(\Omega;\mathbb{R}^m)$. It follows that
$\mathrm{B}_\varepsilon[u_\varepsilon, v_\varepsilon] = \mathrm{B}_\varepsilon^*[v_\varepsilon,u_\varepsilon]$ and
\begin{equation}
\begin{aligned}
& < \mathcal{L}_\varepsilon(u_\varepsilon),v_\varepsilon > =
\int_\Omega \big(A_\varepsilon\nabla u_\varepsilon + V_\varepsilon u_\varepsilon\big)\nabla v_\varepsilon dx
- \int_\Omega u_\varepsilon \text{div}\big(B_\varepsilon v_\varepsilon\big) dx
+ \int_\Omega \big(c_\varepsilon + \lambda\big) u_\varepsilon v_\varepsilon dx
= < u_\varepsilon, \mathcal{L}^*_\varepsilon(v_\varepsilon) >.
\end{aligned}
\end{equation} }
\end{remark}

\begin{lemma}\label{lemma:2.2}
Let $\Omega$ be a Lipschitz domain. Suppose that $A$ satisfies the ellipticity condition $\eqref{a:1}$, and other coefficients of $\mathcal{L}_\varepsilon$ satisfy $\eqref{a:3}$.
Then we have the following properties: for any $u,v\in H^1_0(\Omega;\mathbb{R}^m)$,
\begin{equation}\label{pri:2.1}
 \big|\mathrm{B}_\varepsilon[u,v]\big| \leq C \|u\|_{H^1_0(\Omega)}\|v\|_{H^1_0(\Omega)},
 \quad \text{and} \quad c_0 \| u\|_{H^1_0(\Omega)}^2 \leq \mathrm{B}_\varepsilon[u,u],
 ~~ \text{whenever} ~\lambda\geq\lambda_0,
\end{equation}
where  $\lambda_0 = \lambda_0(\mu,\kappa,m,d)$ is sufficiently large. Note that $C$ depends on $\mu,\kappa,\lambda,m,d,\Omega$, while
$c_0$ depends on
$\mu,\kappa,m,d,\Omega$.
\end{lemma}

\begin{thm}\label{thm:2.1}
The coefficients of $\mathcal{L}_\varepsilon$ and $\lambda_0$ are as in Lemma $\ref{lemma:2.2}$.
Suppose $F\in H^{-1}(\Omega;\mathbb{R}^m)$ and $g\in H^{\frac{1}{2}}(\partial\Omega;\mathbb{R}^m)$. Then
the Dirichlet boundary value problem
$\mathcal{L}_{\varepsilon}(u_\varepsilon) = F$ in $\Omega$ and $u_\varepsilon = g$ on $\partial\Omega$
has a unique weak solution $u_\varepsilon \in H^1(\Omega)$, whenever $\lambda\geq\lambda_0$,
and the solution satisfies the uniform estimate
\begin{equation}\label{pri:2.3}
\|u_\varepsilon\|_{H^1(\Omega)} \leq C \big\{\|F\|_{H^{-1}(\Omega)}+ \|g\|_{H^{1/2}(\partial\Omega)}\big\},
\end{equation}
where $C$ depends only on $\mu,\kappa,m,d$ and $\Omega$. Moreover, with one more the periodicity condition $\eqref{a:2}$ on the coefficients of
$\mathcal{L}_\varepsilon$,
we then have $u_\varepsilon\rightharpoonup u$ weakly in $H^1(\Omega;\mathbb{R}^m)$ and strongly in $L^2(\Omega;\mathbb{R}^m)$ as $\varepsilon\to 0$,
where $u$ is the weak solution to the homogenized problem
$\mathcal{L}_0(u) = F$ in $\Omega$ and $u = g$ on $\partial\Omega$.
\end{thm}

\begin{remark}
 \emph{ The proof of Lemma $\ref{lemma:2.2}$ follows from the same argument in the scalar case (see \cite{LCE,DGNT}).
 Theorem $\ref{thm:2.1}$ involves the uniqueness and existence of the weak solution to $\eqref{pde:1.1}$,
 and the so-called homogenization theorem associated with $\mathcal{L}_\varepsilon$.
 The proof of Theorem $\ref{thm:2.1}$ follows from
 Lemma $\ref{lemma:2.2}$, Lax-Milgram theorem, Tartar's method of oscillating test functions (see \cite[pp.103]{ABJLGP} or \cite[pp.31]{VSO}).
 We refer the reader to \cite{ACPD} for more details on the Tartar's method. We also mention that
 all the results in Lemma $\ref{lemma:2.2}$ and Theorem $\ref{thm:2.1}$ still hold for
$\mathcal{L}_\varepsilon^*$ and $\mathrm{B}^*_\varepsilon$. }
\end{remark}

\begin{lemma}[Cacciopolli's inequality]\label{lemma:4.4}
Suppose that $A$ satisfies $\eqref{a:1}$, and other coefficients satisfy $\eqref{a:3}$.
Assume that $u_\varepsilon\in H^1_{loc}(\Omega;\mathbb{R}^m)$ is
a weak solution to $\mathcal{L}_\varepsilon(u_\varepsilon) = \emph{div}(f) + F$ in $\Omega$,
where $f=(f_i^\alpha)\in L^2(\Omega;\mathbb{R}^{md})$
and $F\in L^q(\Omega;\mathbb{R}^{m})$ with $q = \frac{2d}{d+2}$.
Then for any $B\subset 2B\subset\Omega$, we have the uniform estimate
\begin{equation}\label{pri:4.10}
\Big(\dashint_B |\nabla u_\varepsilon|^2 dx\Big)^{\frac{1}{2}}
\leq \frac{C}{r}\Big(\dashint_{2B} |u_\varepsilon|^2 dx\Big)^{\frac{1}{2}} + C\Big(\dashint_{2B} |f|^2 dx\Big)^{\frac{1}{2}}
+ rC\Big(\dashint_{2B} |F|^q dx\Big)^{\frac{1}{q}},
\end{equation}
where $C$ depends only on $\mu, \kappa, \lambda, m,d$.
\end{lemma}

\begin{pf}
The proof is standard, and we provide a proof for the sake of completeness. Let $\phi\in C^1_0(\Omega)$ be a cut-off function satisfying
$\phi = 1$ in $B$, $\phi = 0$ outside $2B$, and $|\nabla\phi| \leq 2/r$.
Then let $\varphi = \phi^2u_\varepsilon$ be a test function, it follows that 
\begin{eqnarray*}
\int_\Omega \big[A_\varepsilon^{\alpha\beta}\nabla u_\varepsilon^\beta + V_\varepsilon^{\alpha\beta} u_\varepsilon^\beta\big]\nabla u_\varepsilon^\alpha\phi^2 dx
&+& 2\int_\Omega \big[A_\varepsilon^{\alpha\beta}\nabla u_\varepsilon^{\beta}
+V_\varepsilon^{\alpha\beta}u_\varepsilon^\beta\big]\nabla\phi u_\varepsilon^\alpha\phi dx
+ \int_\Omega B_\varepsilon^{\alpha\beta}\nabla u_\varepsilon^\beta u_\varepsilon^\alpha\phi^2 dx \\
+ \int_\Omega c_\varepsilon^{\alpha\beta}u^\beta_\varepsilon u_\varepsilon^\alpha\phi^2 + \lambda |u_\varepsilon|^2\phi^2 dx
&=& \int_\Omega F^\alpha u_\varepsilon^\alpha\phi^2 dx
-\int_\Omega f^\alpha\nabla u_\varepsilon^\alpha\phi^2 dx - 2\int_\Omega f^\alpha\nabla\phi u_\varepsilon^\alpha\phi dx, \qquad \text{in}~~\Omega.
\end{eqnarray*}
By using the ellipticity condition and Young's inequality with $\delta$,  we have
\begin{eqnarray}\label{f:4.9}
\frac{\mu}{4}\int_\Omega \phi^2|\nabla u_\varepsilon|^2 dx & + & (\lambda - C^\prime)\int_\Omega \phi^2 |u_\varepsilon|^2 dx \nonumber\\
&\leq&  C\int_\Omega |\nabla\phi|^2|u_\varepsilon|^2 dx
+C\int_\Omega \phi^2|f|^2 dx + \int_\Omega \phi^2|F||u_\varepsilon| dx,
\end{eqnarray}
where $C^\prime = C^\prime(\mu,\kappa,m,d)$. This together with
\begin{eqnarray*}
\int_\Omega \phi^2 |F||u_\varepsilon| dx
&\leq& \Big(\int_\Omega \big(\phi|u_\varepsilon|\big)^{2^*}dx\Big)^{\frac{1}{2^*}}
\Big(\int_\Omega \big(\phi|F|\big)^q dx\Big)^{\frac{1}{q}}
\leq  C\Big(\int_\Omega \big|\nabla(\phi u_\varepsilon)\big|^2dx\Big)^{\frac{1}{2}}
\Big(\int_\Omega \big(\phi|F|\big)^q dx\Big)^{\frac{1}{q}}\\
&\leq& \frac{\mu}{8}\int_\Omega |\nabla u_\varepsilon|^2\phi^2 dx + \frac{\mu}{8}\int_\Omega |\nabla\phi|^2 |u_\varepsilon|^2 dx
+ C\big(\int_\Omega |\phi F|^q dx\big)^{2/q}
\end{eqnarray*}
gives $\eqref{pri:4.10}$,
where $2^* = 2d/(d-2)$, and we use H\"{o}lder's inequality, Sobolev's inequality, and Young's inequality in order.
\qed
\end{pf}

\begin{remark}\label{rm:4.1}
\emph{ In fact, $\eqref{pri:4.10}$ is the interior $W^{1,2}$ estimate. By the same argument, we can also derive
the near boundary  Cacciopolli's inequality for the weak solution to
$\mathcal{L}_\varepsilon(u_\varepsilon) = \text{div}(f) + F$ in $D(4r)$ and $u_\varepsilon = 0$ on $\Delta(4r)$,
\begin{equation}\label{pri:4.11}
\Big(\dashint_{D(r)} |\nabla u_\varepsilon|^2 dx\Big)^{\frac{1}{2}}
\leq \frac{C}{r}\Big(\dashint_{D(2r)} |u_\varepsilon|^2 dx\Big)^{\frac{1}{2}} + C\Big(\dashint_{D(2r)} |f-\bar{f}|^2 dx\Big)^{\frac{1}{2}}
+ rC\Big(\dashint_{D(2r)} |F|^q dx\Big)^{\frac{1}{q}}.
\end{equation}
We point out that the constant $C$ in $\eqref{pri:4.10}$ or $\eqref{pri:4.11}$ depends only on $\mu,\kappa,m,d$,
whenever $\lambda \geq \lambda_0 \geq C^\prime$. }
\end{remark}

\begin{remark}\label{rm:2.4}
\emph{Suppose that $A\in\Lambda(\mu,\tau,\kappa)$, and $V$ satisfies $\eqref{a:2}$ and $\eqref{a:4}$.
In view of the interior Schauder estimate (see \cite{MGLM}),  we obtain
\begin{equation}\label{f:2.11}
\max_{0\leq k\leq d}\big\{\|\chi_k\|_{L^\infty(Y)},\|\nabla\chi_k\|_{L^\infty(Y)},
[\nabla\chi_k]_{C^{0,\tau}(Y)}\big\} \leq C(\mu,\tau,\kappa,m,d),
\end{equation}
where $[\nabla\chi_k]_{C^{0,\tau}(Y)}$ is the $\tau^{\text{th}}$-H\"older seminorm of $\nabla \chi_k$ (see
\cite[pp.254]{LCE} for the definition).
If $A\in \text{VMO}(\mathbb{R}^d)$, and $V$ satisfies $\eqref{a:2}$ and $\eqref{a:3}$,
then for any $2\leq p<\infty$ and $\varsigma\in(0,1)$,
\begin{equation}\label{f:2.9}
\max_{0\leq k\leq d}\big\{\|\chi_k\|_{L^\infty(Y)}, [\chi_k]_{C^{0,\varsigma}(Y)},
\|\nabla\chi_k\|_{L^p(Y)}\big\}
\leq C(\mu,\omega,\kappa,p,\varsigma,m,d)
\end{equation}
can be derived from the interior $W^{1,p}$ estimate (see \cite{SBLW1,LI,MGLM}). Note that
in the case of $m = 1$,
due to the De Giorgi-Nash-Moser theorem,
the conditions $\eqref{a:1}$ and $\eqref{a:3}$
are sufficient to derive H\"{o}lder estimates for some $\varsigma\in(0,1)$,
but insufficient to get $W^{1,p}$ estimates
(see \cite{SBLW,DGNT,HQLFH}). In another special case of $d=2$, K.O. Widman \cite{KOW} obtained the H\"{o}lder estimate
by the hole filling technique without any regularity assumption on $A$ (or see \cite{MGLM}). }

\emph{ To handle the convergence rates, some auxiliary functions and their estimates are necessary.
Let $1\leq i,j,l,k\leq d$, $1\leq\alpha,\beta,\gamma\leq m$, and $y = x/\varepsilon$. Define
\begin{equation}\label{f:2.10}
b_{ij}^{\alpha\gamma}(y) = \hat{a}_{ij}^{\alpha\gamma} - a_{ij}^{\alpha\gamma}(y)
- a_{ik}^{\alpha\beta}(y)\frac{\partial}{\partial y_k}\big\{\chi_{j}^{\beta\gamma}\big\}, \qquad
U_i^{\alpha\gamma}(y) =
\widehat{V}_i^{\alpha\gamma} - V_{i}^{\alpha\gamma}(y)-a_{ij}^{\alpha\beta}(y)\frac{\partial}{\partial y_j}\big\{\chi_0^{\beta\gamma}\big\}.
\end{equation}
Note that C.E. Kenig, F. Lin, and Z. Shen \cite{SZW2,SZW1} showed that there exists
$E_{lij}^{\alpha\gamma}\in H^1_{per}(Y)$, such that
$b_{ij}^{\alpha\beta} = \frac{\partial}{\partial y_l}\{E_{lij}^{\alpha\beta}\}$, $E_{lij}^{\alpha\beta} = -E_{ilj}^{\alpha\beta}$
and $\|E_{lij}^{\alpha\beta}\|_{L^\infty(Y)}\leq C(\mu,\omega,m,d)$.
Here we obtain similar results for $U_i^{\alpha\gamma}$: there exists $F_{ki}^{\alpha\gamma}\in H_{per}^1(Y)$ such that
\begin{equation}\label{f:2.4}
U_i^{\alpha\gamma} = \frac{\partial}{\partial y_k}\{F_{ki}^{\alpha\gamma}\},
\quad F_{ki}^{\alpha\gamma}= - F_{ik}^{\alpha\gamma}
\quad \text{and}~~\|F_{ki}^{\alpha\gamma}\|_{L^\infty(Y)}\leq C(\mu,\omega,\kappa,m,d).
\end{equation}
We give a proof of $\eqref{f:2.4}$ for the sake of completeness.
In view of $\eqref{pde:2.1}$ and $\eqref{f:2.2}$, we have $\int_Y U_i^{\alpha\gamma}(y)dy = 0$
and $\frac{\partial}{\partial y_i}\{U_i^{\alpha\gamma}\} = 0$.
Then there exists a unique solution $\theta_i^{\alpha\gamma}\in H_{per}^1(Y)$ satisfying
\begin{equation*}
\Delta \theta_i^{\alpha\gamma} = U_i^{\alpha\gamma} \quad \text{in} \quad \mathbb{R}^d, \quad\qquad \int_Y \theta_i^{\alpha\gamma} = 0
\end{equation*}
(see \cite{ABJLGP,ACPD}).
Let $F_{ki}^{\alpha\gamma} = \frac{\partial}{\partial y_k}\{\theta_i^{\alpha\gamma}\} - \frac{\partial}{\partial y_i}\{\theta_k^{\alpha\gamma}\}$,
and obviously
$F_{ki}^{\alpha\gamma} = - F_{ik}^{\alpha\gamma}$. We mention that $\theta_i^{\alpha\gamma}\in H^2_{per}(Y)$,
which implies $F_{ki}^{\alpha\gamma} \in H_{per}^1(Y)$. Next, we verify $F_{ki}^{\alpha\gamma}$ is bounded, which is equivalent to
\begin{equation}\label{f:2.3}
\max_{1\leq i\leq d}\big\{\|\nabla\theta_i^{\alpha\gamma}\|_{L^\infty(Y)}\big\}\leq C(\mu,\omega,\kappa,m,d).
\end{equation}
Observe that $\eqref{f:2.9}$ gives $U_i^{\alpha\gamma}\in L^p_{per}(Y)$, then the estimate
$\eqref{f:2.3}$ follows from the $L^p$ estimates and the Sobolev embedding theorem, provided $p>d$.
Moreover, we have
\begin{equation*}
 \frac{\partial}{\partial y_k}\big\{F_{ki}^{\alpha\gamma}\big\} = \Delta\{\theta_i^{\alpha\gamma}\}
 - \frac{\partial^2}{\partial y_k\partial y_i}\big\{\theta_k^{\alpha\gamma}\big\} = U_i^{\alpha\gamma}.
\end{equation*}
Note that $\frac{\partial}{\partial y_k}\{U_k^{\alpha\gamma}\} = 0$ in $\mathbb{R}^d$,
which implies $\Delta\frac{\partial}{\partial y_k}\{\theta_k^{\alpha\gamma}\} = 0$.
In view of Liouvill's theorem (see \cite{LCE}), we have $\frac{\partial}{\partial y_k}\{\theta_k^{\alpha\gamma}\} = C$, therefore
$\frac{\partial}{\partial y_k}\frac{\partial}{\partial y_i}\{\theta_k^{\alpha\gamma}\}
= \frac{\partial}{\partial y_i}\frac{\partial}{\partial y_k}\{\theta_k^{\alpha\gamma}\} = 0$, and we complete the proof of $\eqref{f:2.4}$.}

\emph{ In addition, we define the auxiliary functions $\vartheta_{i}^{\alpha\gamma}$ and $\zeta^{\alpha\gamma}$ as follows:
\begin{equation}\label{f:2.6}
\left.\begin{aligned}
\Delta\vartheta_{i}^{\alpha\gamma}
 & = W_i^{\alpha\gamma} := \widehat{B}_i^{\alpha\gamma} - B_{i}^{\alpha\gamma}(y)
- B_{j}^{\alpha\beta}(y)\frac{\partial}{\partial y_j}\big\{\chi_i^{\beta\gamma}\big\} & \quad \text{in} ~~\mathbb{R}^d,
&\qquad \int_Y \vartheta_i^{\alpha\gamma} dy = 0 ,\\
\Delta \zeta^{\alpha\gamma}
&= Z^{\alpha\gamma} := \widehat{c}^{\alpha\gamma} - c^{\alpha\gamma}(y)
- B_{i}^{\alpha\beta}(y)\frac{\partial}{\partial y_i}\big\{\chi_0^{\beta\gamma}\big\} & \quad \text{in}~~\mathbb{R}^d ,
&\qquad \int_Y \zeta^{\alpha\gamma} dy = 0 .
\end{aligned}\right.
\end{equation}
It follows from $\eqref{f:2.2}$ that $\int_Y W_i^{\alpha\gamma}(y)dy = 0$ and $\int_Y Z^{\alpha\gamma}(y)dy = 0$, which
implies the existence of $\vartheta_{i}^{\alpha\gamma}$ and $\zeta^{\alpha\gamma}$. By the same argument, it follows from $\eqref{f:2.9}$ that
$\max\big\{\|\nabla\vartheta_i^{\alpha\gamma}\|_{L^\infty(Y)}, \|\nabla\zeta^{\alpha\gamma}\|_{L^\infty(Y)}\big\}\leq C(\mu,\omega,\kappa,d,m)$. }

\emph{ We end this remark by a summary. Suppose that $A\in \text{VMO}(\mathbb{R}^d)$,
and the coefficients of $\mathcal{L}_\varepsilon$ satisfy $\eqref{a:1}-\eqref{a:3}$, then we have
\begin{equation}\label{f:2.7}
\max_{
  1\leq i,j,l,k\leq d \atop
  1\leq \alpha,\gamma\leq m} \big\{\|E_{lij}^{\alpha\gamma}\|_{L^\infty(Y)},
 \|F_{ki}^{\alpha\gamma}\|_{L^\infty(Y)},
 \|\nabla\vartheta_i^{\alpha\gamma}\|_{L^\infty(Y)}, \|\nabla\zeta^{\alpha\gamma}\|_{L^\infty(Y)} \big\} \leq C(\mu,\omega,\kappa,m,d).
\end{equation}
In the special case of $m = 1$ or $d=2$,
the estimate $\eqref{f:2.7}$ still holds without
any regularity assumption on $A$. }
\end{remark}

We now introduce the Lipschitz estimate and Schauder estimate that will be frequently employed later.
Let $L(u) = -\text{div}(A\nabla u)$ and $\mathcal{L}(u) = -\text{div}(A\nabla u + Vu) + B\nabla u + (c+\lambda I)u$. Then we have the following results:

\begin{lemma}\label{lemma:2.3} Let $\Omega$ be a bounded $C^{1,\tau}$ domain.
Suppose $A$ satisfies $\eqref{a:1}$ and $\eqref{a:4}$. Let $u$ be the weak solution to $L(u) = \emph{div}(f) + F$ in $\Omega$ and $u=0$ on $\partial\Omega$,
where $f\in C^{0,\sigma}(\Omega;\mathbb{R}^{md})$ with $\sigma\in (0,\tau]$, and $F\in L^p(\Omega;\mathbb{R}^m)$ with $p>d$. Then we have:

\emph{(i)} the Schauder estimate
\begin{equation}\label{pri:2.5}
\big[\nabla u\big]_{C^{0,\sigma}(\Omega)} \leq C\big\{\|f\|_{C^{0,\sigma}(\Omega)} + \|F\|_{L^p(\Omega)}\big\};
\end{equation}

\emph{(ii)} the Lipschitz estimate
\begin{equation}\label{pri:2.4}
\|\nabla u\|_{L^\infty(\Omega)} \leq C\big\{[f]_{C^{0,\sigma}(\Omega)} + \|F\|_{L^p(\Omega)}\big\},
\end{equation}
where $C$ depends on $\mu,\tau,\kappa,\sigma,m,d,p$ and $\Omega$. Moreover, if $u =g$ on $\partial\Omega$
with $g\in C^{1,\sigma}(\partial\Omega;\mathbb{R}^m)$, then we have
\begin{equation}\label{pri:2.8}
\big[\nabla u\big]_{C^{0,\sigma}(\Omega)} \leq C\big\{\|f\|_{C^{0,\sigma}(\Omega)} + \|F\|_{L^p(\Omega)}
+ \|g\|_{C^{1,\sigma}(\partial\Omega)}\big\}.
\end{equation}
\end{lemma}

\begin{pf}
The results are standard, and we provide a proof for the sake of completeness.
For (i), we refer the reader to \cite[pp.75-95]{MGLM} for the details. For (ii),
due to the properties of Green function (denoted by $G(x,y)$) associated with
$L$: $|G(x,y)|\leq C|x-y|^{2-d}$, $|\nabla_x G(x,y)|\leq C|x-y|^{1-d}$ and $|\nabla_x\nabla_y G(x,y)|\leq C|x-y|^{-d}$
(the existence of $G(x,y)$ is included in \cite[Theorem 4.1]{SHSK}), we have
$u(x) = -\int_\Omega \nabla_y G(x,y) f(y)dy + \int_\Omega G(x,y) F(y) dy$. Differentiating both sides with respect to $x$ gives
\begin{equation*}
\begin{aligned}
\nabla u(x) & = -\int_\Omega \nabla_x\nabla_y G(x,y)\big[f(y)-f(x)\big]dy - f(x)\int_\Omega \nabla_x\nabla_y G(x,y) dy
+\int_\Omega \nabla_xG(x,y) F(y) dy \\
& = -\int_\Omega \nabla_x\nabla_y G(x,y)\big[f(y)-f(x)\big]dy -f(x)\int_{\partial\Omega} n(y)\nabla_x G(x,y) dS(y)
+ \int_\Omega \nabla_xG(x,y) F(y) dy \\
& = -\int_\Omega \nabla_x\nabla_y G(x,y)\big[f(y)-f(x)\big]dy + \int_\Omega \nabla_xG(x,y) F(y) dy
\end{aligned}
\end{equation*}
for any $x\in\Omega$. Note that we use the integration by parts in the second equality,
and the fact of $\nabla_x G(x,\cdot) = 0$ on $\partial\Omega$ in the last equality.
This implies $\eqref{pri:2.4}$. Finally we can use the extension technique to obtain the estimate $\eqref{pri:2.8}$.
(Its proof is similar to that in the proof of Theorem $\ref{thm:1.1}$, and we refer the reader to
\cite[pp.136-138]{DGNT} for the extension lemmas.)
\qed
\end{pf}

\begin{remark}\label{rm:2.8}
\emph{Set $U(x,r) = \Omega\cap B(x,r)$
for any $x\in\Omega$. Let $\varphi\in C^\infty_0(2B)$ be a cut-off function
satisfying $\varphi = 1$ in $B$, $\varphi = 0$ outside $3/2B$, and $|\nabla \varphi|\leq C/r$.
Let $w=u\varphi$, where $u$ is given in Lemma $\ref{lemma:2.3}$. Then we have
\begin{equation*}
 L(w) = \text{div}(f\varphi) -f\cdot\nabla\varphi + F\varphi - \text{div}(A\nabla\varphi u) - A\nabla u\nabla\varphi
 \quad \text{in}~\Omega,
 \qquad w = 0 \quad \text{on} ~\partial\Omega.
\end{equation*}
We apply the estimate $\eqref{pri:2.4}$ to the above equation with $r=1$, and obtain
\begin{equation}\label{f:2.20}
\begin{aligned}
\|\nabla u\|_{L^\infty(U(x,1))}
&\leq C\big\{\|u\|_{C^{0,\sigma}(U(x,2))}+\|u\|_{L^{p}(U(x,2))}+\|f\|_{L^\infty(U(x,2))}
+ [f]_{C^{0,\sigma}(U(x,2))} + \|F\|_{L^p(U(x,2))} \big\}\\
&\leq C\big\{\|u\|_{W^{1,s}(U(x,2))}+\|f\|_{L^\infty(U(x,2))}
+ [f]_{C^{0,\sigma}(U(x,2))} + \|F\|_{L^p(U(x,2))} \big\}
\end{aligned}
\end{equation}
where $s=\max\{p,[\frac{d}{1-\sigma}]+1\}$ ($[\frac{d}{1-\sigma}]$ is the integer part of $\frac{d}{1-\sigma}$),
and we employ the fact of $\|u\|_{C^{0,\sigma}(U(x,2))}\leq C \|u\|_{W^{1,s}(U(x,2))}$.
On account of
\begin{equation*}
  \|\nabla u\|_{L^{p}(U(x,2))}\leq C\{\|\nabla u\|_{L^2(U(x,4))}+ \|f\|_{L^p(U(x,4))} + \|F\|_{L^q(U(x,4))}\},
\end{equation*}
where $q = \frac{pd}{d+p}$ (see \cite[Theorem 7.2]{MGLM}), we have
\begin{equation}\label{f:2.21}
\|u\|_{W^{1,s}(U(x,2))}
\leq C\{\|u\|_{W^{1,2}(U(x,4))}+ \|f\|_{L^p(U(x,4))} + \|F\|_{L^q(U(x,4))}\}
\end{equation}
where we use $\|u\|_{L^s(U(x,2))}\leq C\big\{\|\nabla u\|_{L^s(U(x,2))} + \|u\|_{L^2(U(x,2))}\big\}$
in the above inequality (see $\eqref{f:4.7}$). 
Combining $\eqref{f:2.20}$ and $\eqref{f:2.21}$, we have
\begin{equation*}
\|\nabla u\|_{L^\infty(U(x,1))}
\leq C\big\{\|\nabla u\|_{L^{2}(U(x,4))}+\|u\|_{L^{2}(U(x,4))}+\|f\|_{L^\infty(U(x,4))} + [f]_{C^{0,\sigma}(U(x,4))} + \|F\|_{L^p(U(x,4))} \big\}.
\end{equation*}
Note that if $U(x,4)\subset\Omega$, then $v=u-\bar{u}$ is still a solution to $L(u) = \text{div}(f) + F$ in $\Omega$,
where $\bar{u} = \dashint_{U(x,4)} u dy$.
In this case, the above estimate becomes
\begin{equation*}
\|\nabla u\|_{L^\infty(U(x,1))}
\leq C\big\{\|\nabla u\|_{L^{2}(U(x,4))}+\|u-\bar{u}\|_{L^{2}(U(x,4))}+\|f\|_{L^\infty(U(x,4))}
+ [f]_{C^{0,\sigma}(U(x,4))} + \|F\|_{L^p(U(x,4))} \big\}.
\end{equation*}
If $|U(x,4)\cap\partial\Omega|\geq \kappa_0>0$,
then we have $\|u\|_{L^2(U(x,4))}\leq C\|\nabla u\|_{L^2(U(x,4))}$ due to $u=0$ on $\partial\Omega$.
(It is another type of Poincar\'e's inequality, which can be derived from the trace theorem coupled with Rellich's theorem
by a contradiction argument, and the proof is left to the reader.)
In all, we are able to apply Poincar\'e's inequality to both of the two cases, and then obtain
\begin{equation}\label{f:2.22}
\|\nabla u\|_{L^\infty(U(x,1))}
\leq C\big\{\|\nabla u\|_{L^{2}(U(x,4))}+\|f\|_{L^\infty(U(x,4))} + [f]_{C^{0,\sigma}(U(x,4))} + \|F\|_{L^p(U(x,4))} \big\}.
\end{equation}
Next, we let $v(y) = u(ry)$, where $y\in U(x,4)$. Then we have $\tilde{L}(v) = \text{div}(\tilde{f})+\tilde{F}$ in $ U(x,4)$
and $v = 0$ on $\partial\Omega\cap U(x,4)$, where
$\tilde{L} = \frac{\partial}{\partial y_i}a_{ij}^{\alpha\beta}(ry)\frac{\partial}{\partial y_j}$, $\tilde{f}(y) = rf(ry)$ and
$\tilde{F}(y) =r^2F(ry)$. It follows from $\eqref{f:2.22}$ that
\begin{equation*}
\|\nabla v\|_{L^\infty(U(x,1))}
\leq C\big\{\|\nabla v\|_{L^{2}(U(x,4))}+\|\tilde{f}\|_{L^\infty(U(x,4))} + [\tilde{f}]_{C^{0,\sigma}(U(x,4))}
+ \|\tilde{F}\|_{L^p(U(x,4))} \big\},
\end{equation*}
and by change of variable, we have
\begin{equation*}
 \|\nabla u\|_{L^\infty(U(x,r))} \leq C\bigg\{\Big(\dashint_{U(x,4r)}|\nabla u|^2dy\Big)^{1/2}+\|f\|_{L^\infty(U(x,4r))}
 +r^{\sigma}\big[f\big]_{C^{0,\sigma}(U(x,4r))} + r\Big(\dashint_{U(x,4r)}|F|^pdy\Big)^{1/p}\bigg\}.
\end{equation*}
By a covering technique (shown in the proof of Theorem $\ref{thm:5.1}$), we have
\begin{equation}\label{f:2.15}
 \|\nabla u\|_{L^\infty(U)} \leq C\bigg\{\Big(\dashint_{2U}|\nabla u|^2dy\Big)^{1/2}+\|f\|_{L^\infty(2U)}
 +r^{\sigma}\big[f\big]_{C^{0,\sigma}(2U)} + r\Big(\dashint_{2U}|F|^pdy\Big)^{1/p}\bigg\}.
\end{equation}
Here $U$ is the abbreviation of $U(x,r)$ and $2U = U(x,2r)$.
We mention that all the above proof is so-called localization argument,
which gives a way to obtain ``local estimates'' (such as interior estimates and boundary estimates)
from corresponding ``global estimates''. The main point is based on cut-off function coupled with rescaling technique.
So, on account of the estimate $\eqref{pri:2.5}$, following the same procedure as before, it is not hard to derive
\begin{equation}\label{f:2.16}
\big[\nabla u\big]_{C^{0,\sigma}(U)} \leq Cr^{-\sigma}\bigg\{\Big(\dashint_{2U}|\nabla u|^2dy\Big)^{1/2}
+\|f\|_{L^\infty(2U)}
 +r^{\sigma}\big[f\big]_{C^{0,\sigma}(2U)} + r\Big(\dashint_{2U}|F|^pdy\Big)^{1/p}\bigg\},
\end{equation}
and by $\eqref{pri:2.8}$,
\begin{equation}\label{f:2.19}
\big[\nabla u\big]_{C^{0,\sigma}(D(r))} \leq Cr^{-\sigma}\bigg\{\Big(\dashint_{D(2r)}|\nabla u|^2dy\Big)^{1/2}
+r^{-1}\|g\|_{L^\infty(\Delta(2r))}
 + \|\nabla g\|_{L^\infty(\Delta(2r))}
 + r^{\sigma}\big[\nabla g\big]_{C^{0,\sigma}(\Delta(2r))} \bigg\}
\end{equation}
holds for $u$ satisfying $L(u) = 0$ in $D(2r)$ and $u=g$ on $\Delta(2r)$.
We note that the estimate $\eqref{f:2.15}$ is of help to arrive at $\eqref{f:2.16}$,
and the extension technique (see \cite[pp.136]{DGNT}) is used in $\eqref{f:2.19}$.
The details of the proof are omitted.
Finally we remark that $\eqref{f:2.19}$ is exactly the Schauder estimate at boundary,
which can be directly proved (see \cite[Theorem 5.21]{MGLM}).}
\end{remark}

\begin{lemma}\label{lemma:2.4} Let $\Omega$ be a bounded $C^{1,\tau}$ domain, and $\sigma\in(0,\tau]$.
Suppose that $A,V$ satisfy $\eqref{a:1}$ and $\eqref{a:4}$, and $B,c$ satisfy $\eqref{a:3}$.
Let $u$ be the weak solution to $\mathcal{L}(u) = \emph{div}(f) + F$ in $\Omega$ and $u=0$ on $\partial\Omega$,
where $f,F$ satisfy the same conditions as in Lemma $\ref{lemma:2.3}$. Then we have

\emph{(i)} the Lipschitz estimate
\begin{equation}\label{pri:2.6}
\|\nabla u\|_{L^\infty(\Omega)}\leq C\big\{\|f\|_{C^{0,\sigma}(\Omega)} + \|F\|_{L^p(\Omega)}\big\};
\end{equation}

\emph{(ii)} the Schauder estimate
\begin{equation}\label{pri:2.7}
\big[\nabla u\big]_{C^{0,\sigma}(\Omega)} \leq C\big\{\|f\|_{C^{0,\sigma}(\Omega)} + \|F\|_{L^p(\Omega)}\big\},
\end{equation}
where $C$ depends on $\mu,\tau,\kappa,m,d,\sigma,p$ and $\Omega$.
\end{lemma}

\begin{pf}
The results are classical, and we offer a sketch of the proof. First we rewrite $\mathcal{L}(u) = \text{div}(f) + F$
as $L(u) = \text{div}(f+Vu) - B\nabla u -(c+\lambda I)u+ F$ in $\Omega$. It follows from
$\eqref{pri:2.4}$ that
\begin{equation}\label{f:2.12}
\begin{aligned}
\|\nabla u\|_{L^{\infty}(\Omega)}
& \leq C\big\{[f]_{C^{0,\sigma}(\Omega)}+\|F\|_{L^p(\Omega)} + \|u\|_{C^{0,\sigma}(\Omega)}+\|u\|_{W^{1,p}(\Omega)}\big\} \\
& \leq C\big\{[f]_{C^{0,\sigma}(\Omega)}+\|F\|_{L^p(\Omega)} + 2\|\nabla u\|_{L^\infty(\Omega)}^\sigma\|u\|_{L^\infty(\Omega)}^{1-\sigma}
+\|u\|_{W^{1,p}(\Omega)}\big\} \\
& \leq C\big\{[f]_{C^{0,\sigma}(\Omega)}+\|F\|_{L^p(\Omega)} + \|u\|_{L^{\infty}(\Omega)}+\|u\|_{W^{1,p}(\Omega)}\big\}
+ \frac{1}{2}\|\nabla u\|_{L^\infty(\Omega)},
\end{aligned}
\end{equation}
where we use Young's inequality in the last inequality. By the Sobolev embedding theorem we have
$\|u\|_{L^\infty(\Omega)}\leq C\|u\|_{W^{1,p}(\Omega)}$ for $p>d$. This
together with $\eqref{f:2.12}$ leads to
\begin{equation*}
\begin{aligned}
\|\nabla u\|_{L^{\infty}(\Omega)}
& \leq C\big\{[f]_{C^{0,\sigma}(\Omega)}+\|F\|_{L^p(\Omega)} +\|u\|_{W^{1,p}(\Omega)}\big\} \\
& \leq C\big\{\|f\|_{C^{0,\sigma}(\Omega)}+\|F\|_{L^p(\Omega)}\big\},
\end{aligned}
\end{equation*}
where we use the $W^{1,p}$ estimate with $1<p<\infty$ in the last inequality,
which can be derived by a similar argument as in the proof of Theorem $\ref{thm:1.1}$, or see \cite{SBLW1,MGLM,DGNT}.

It remains to show (ii). In view of $\eqref{pri:2.5}$ and $\eqref{pri:2.6}$, we obtain
\begin{equation*}
\begin{aligned}
\big[\nabla u\big]_{C^{0,\sigma}(\Omega)}
& \leq C\big\{\|f\|_{C^{0,\sigma}(\Omega)}+\|F\|_{L^p(\Omega)}
+\|u\|_{L^\infty(\Omega)}+\|\nabla u\|_{L^\infty(\Omega)}+\|u\|_{W^{1,p}(\Omega)}\big\} \\
& \leq C\big\{\|f\|_{C^{0,\sigma}(\Omega)}+\|F\|_{L^p(\Omega)}\big\},
\end{aligned}
\end{equation*}
where we also use $W^{1,p}$ estimate in the last inequality. The proof is completed.
\qed
\end{pf}

\begin{remark}\label{rm:2.7}
\emph{Let $u$ be given in Lemma $\ref{lemma:2.4}$.
Applying the localization argument (see Remark $\ref{rm:2.8}$) to the estimates $\eqref{pri:2.6}$ and $\eqref{pri:2.7}$,
we can similarly give the corresponding local estimates:
\begin{equation}\label{f:2.17}
 \|\nabla u\|_{L^\infty(U)} \leq C\bigg\{r^{-1}\Big(\dashint_{2U}|u|^2dy\Big)^{1/2}+\|f\|_{L^\infty(2U)}
 +r^{\sigma}\big[f\big]_{C^{0,\sigma}(2U)} + r\Big(\dashint_{2U}|F|^pdy\Big)^{1/p}\bigg\},
\end{equation}
and
\begin{equation}\label{f:2.18}
\big[\nabla u\big]_{C^{0,\sigma}(U)} \leq Cr^{-\sigma}\bigg\{r^{-1}\Big(\dashint_{2U}|u|^2dy\Big)^{1/2}+\|f\|_{L^\infty(2U)}
 +r^{\sigma}\big[f\big]_{C^{0,\sigma}(2U)} + r\Big(\dashint_{2U}|F|^pdy\Big)^{1/p}\bigg\},
\end{equation}
where $C$ depends on $\mu,\tau,\kappa,m,d,p,\sigma$ and $M_0$.
We mention that in the proof of $\eqref{f:2.17}$,
we also need $W^{1,p}$ estimate like \cite[Theorem 7.2]{MGLM} for $\mathcal{L}$.
It can be established by using the bootstrap method which is exactly shown in the proof of Theorem $\ref{thm:4.1}$,
so we do not repeat them. The remainder of the argument is analogous to that in Remark $\ref{rm:2.8}$. }
\end{remark}

\section{$W^{1,p}$ estimates \& H\"older estimates}

\begin{lemma}\label{lemma:4.1}
Let $2\leq p< \infty$. Suppose that $A\in \emph{VMO}(\mathbb{R}^d)$ satisfies $\eqref{a:1}$ and $\eqref{a:2}$.
Assume $f=(f_i^\alpha)\in L^p(\Omega;\mathbb{R}^{md})$,
and $F\in L^q(\Omega;\mathbb{R}^m)$
with $q = \frac{pd}{p+d}$. Then the weak solution to $L_\varepsilon(u_\varepsilon) = \emph{div}(f) + F$ in $\Omega$
and $u_\varepsilon = 0$ on $\partial \Omega$ satisfies the uniform estimate
\begin{equation}\label{pri:4.1}
\|\nabla u_\varepsilon\|_{L^p(\Omega)} \leq C\{\|f\|_{L^p(\Omega)} + \|F\|_{L^{q}(\Omega)} \},
\end{equation}
where $C$ depends only on $\mu,\omega,\kappa, p, q, d, m$ and $\Omega$~.
\end{lemma}

\begin{remark}
\emph{ The estimate $\eqref{pri:4.1}$ actually holds for $1<p<\infty$, where $q = \frac{pd}{p+d}$ if $p>\frac{d}{d-1}$, and $q> 1$ if $1< p \leq \frac{d}{d-1}$.
In the case of $F = 0$, $\eqref{pri:4.1}$ is shown in \cite{MAFHL3}. If $F\not=0$, we can derive the above result
by the duality argument applied in Lemma $\ref{lemma:4.3}$. The same method may be found in \cite{GJ,SZW4}.
Besides, we refer the reader to \cite{SZW10} for the
sharp range of $p$'s on Lipschitz domains. }
\end{remark}

\begin{thm}\label{thm:4.1}
$(\text{Interior}~W^{1,p}~\text{estimates})$~\textbf{.}
Let $2\leq p < \infty$. Suppose that
$A\in \emph{VMO}(\mathbb{R}^d)$ satisfies $\eqref{a:1}$, $\eqref{a:2}$, and other coefficients satisfy $\eqref{a:3}$.
Assume that $u_\varepsilon\in H^1_{loc}(\Omega;\mathbb{R}^m)$ is a weak solution to $\mathcal{L}_\varepsilon(u_\varepsilon) = \emph{div}(f) + F$ in $\Omega$,
where $f\in L^p(\Omega;\mathbb{R}^{md})$ and $F\in L^{q}(\Omega;\mathbb{R}^m)$ with $q = \frac{pd}{p+d}$. Then,
we have $|\nabla u_\varepsilon|\in L^p(B)$ and the uniform estimate
\begin{equation}\label{pri:4.7}
 \left(\dashint_{B} |\nabla u_\varepsilon|^{p}dx\right)^{\frac{1}{p}}
 \leq \frac{C}{r}\left(\dashint_{2B} |u_\varepsilon|^2 dx\right)^{\frac{1}{2}} + C\left\{\Big(\dashint_{2B} |f|^p dx\Big)^{\frac{1}{p}}
 + r\Big(\dashint_{2B} |F|^q dx\Big)^{\frac{1}{q}}\right\}
\end{equation}
for any $B\subset 2B\subset\Omega$ with $0<r\leq 1$, where $C$ depends only on $\mu, \omega, \kappa, \lambda, p, m, d$.
\end{thm}

\begin{pf}
By rescaling we may assume $r = 1$. In the case of $p = 2$, $\eqref{pri:4.7}$ follows from Lemma $\ref{lemma:4.4}$.
Next, we will prove $\eqref{pri:4.7}$ for $p\in [2,p_*]$, where $p_* = \frac{2d}{d-2k_0}$, and $k_0 < d/2 \leq k_0+1$.
To do so, let $w_\varepsilon = \varphi u_\varepsilon$,
where $\varphi \in C_0^\infty(2B)$ is a cut-off function satisfying $\varphi = 1$ in $B$, $\varphi = 0$ outside $3/2B$, and
$|\nabla \varphi|\leq C$.
We rewrite the original systems as
\begin{equation*}
 -\text{div}(A_\varepsilon\nabla w_\varepsilon) = \text{div}(f\varphi) - f\cdot\nabla\varphi + F\varphi + \tilde{F}
 \quad \text{in} ~\Omega,
\end{equation*}
where
\begin{equation*}
\tilde{F}^\alpha = \text{div}(V_\varepsilon^{\alpha\beta} w_\varepsilon^\beta  - A_\varepsilon^{\alpha\beta}\nabla\varphi u_\varepsilon^\beta)
 -A_\varepsilon^{\alpha\beta}\nabla\varphi\nabla u_\varepsilon^\beta - V_\varepsilon^{\alpha\beta}\nabla\varphi u_\varepsilon^\beta
 - B_\varepsilon^{\alpha\beta}\nabla u_\varepsilon^\beta\varphi
 - c_\varepsilon^{\alpha\beta}w_\varepsilon^\beta - \lambda w_\varepsilon^\alpha.
\end{equation*}
Hence it follows from $\eqref{pri:4.10}$, $\eqref{pri:4.1}$, $\eqref{f:4.6}$, and H\"{o}lder's inequality that
\begin{eqnarray*}
\|\nabla u_\varepsilon\|_{L^{p}(B)}
&\leq & C\left\{\|u_\varepsilon\|_{L^{p}(2B)} + \|\nabla u_\varepsilon\|_{L^{q}(2B)} + \|u_\varepsilon\|_{L^{q}(2B)}
                 + \|f\|_{L^{p}(2B)} + \|F\|_{L^{q}(2B)}\right\}  \\
&\leq & C\left\{ \|\nabla u_\varepsilon\|_{L^{q}(2B)} + \|u_\varepsilon\|_{L^2(2B)} + \|f\|_{L^{p}(2B)} + \|F\|_{L^{q}(2B)} \right\} \\
&\leq & C\left\{\|u_\varepsilon\|_{L^2(2^{k_0+1}B)} + \|f\|_{L^{p}(2^{k_0+1}B)} + \|F\|_{L^{q}(2^{k_0+1}B)}\right\}.
\end{eqnarray*}
We first check the special case of $p=p_*$ to obtain the final step $k_0$ of iteration,
and then verify the above inequality for any $p\in [2,\frac{2d}{d-2}]$. Second, it is not hard to extend the range of $p$'s to $[2,p_*]$ by
at most $k_0$ times of iteration. The rest of the proof is to extend the $p$'s range to $p_*<p<\infty$.
Indeed it is true, since
\begin{eqnarray*}
\|\nabla u_\varepsilon\|_{L^{p}(B)}
&\leq & C\left\{ \|\nabla u_\varepsilon\|_{L^{q}(2B)} + \|u_\varepsilon\|_{L^2(2B)} + \|f\|_{L^{p}(2B)} + \|F\|_{L^{q}(2B)} \right\}
\end{eqnarray*}
and $q < d\in [2,p_*]$, which is exactly the start point for iterations due to the previous case.

Hence, let $N = 4^{k_0}$, we have proved
\begin{equation}\label{f:4.18}
 \|\nabla u_\varepsilon\|_{L^p(B)} \leq C\big\{\|u_\varepsilon\|_{L^2(NB)} + \|f\|_{L^p(NB)} + \|F\|_{L^q(NB)}\big\}
\end{equation}
for any $2\leq p<\infty$ in the case of $r=1$. We remark that (i) the estimate $\eqref{f:4.18}$ uniformly holds for $\varepsilon >0$;
(ii) the constant in $\eqref{f:4.18}$ can be given by
$C\leq C(\mu,\omega,m,d,p)\big\{\|A\|_{L^\infty(\mathbb{R}^d)}+\|V\|_{L^\infty(\mathbb{R}^d)}
+ \|B\|_{L^\infty(\mathbb{R}^d)} + \|c\|_{L^\infty(\mathbb{R}^d)}+\lambda\big\}^{k_0+2}$.
The two points make the rescaling argument valid when we study the estimate $\eqref{pri:4.7}$ for $0<r<1$.

We now let $v_\varepsilon(x) = u_\varepsilon(rx)$, where $x\in B_1$ and $r\in (0,1)$. Hence we have
\begin{equation}
\tilde{\mathcal{L}}_{\frac{\varepsilon}{r}}(v_\varepsilon)=-\text{div}\big[A(rx/\varepsilon)\nabla v_\varepsilon + \tilde{V}(rx/\varepsilon)v_\varepsilon \big] + \tilde{B}(rx/\varepsilon)\nabla v_\varepsilon +
\tilde{c}(rx/\varepsilon)v_\varepsilon+ \tilde{\lambda} v_\varepsilon = \text{div}(\tilde{f}) + \tilde{F}   \quad \text{in} \quad NB_1,
\end{equation}
where
\begin{eqnarray*}
\tilde{V}(x) = rV(x),
\quad \tilde{B} = rB(x), \quad \tilde{c} = r^2c(x), \quad \tilde{\lambda} = r^2\lambda,
\quad \tilde{f} = rf(rx), \quad \tilde{F} = r^2F(rx).
\end{eqnarray*}
It is clear to see that the coefficients of $\tilde{\mathcal{L}}$ satisfy the same assumptions as $\mathcal{L}$ in this theorem.
Set $\varepsilon^\prime = \varepsilon/r$, and applying $\eqref{f:4.18}$ directly, we obtain
\begin{eqnarray*}
  \|\nabla v_{r\varepsilon^\prime}\|_{L^p(B_1)} \leq C\big\{\|v_{r\varepsilon^\prime}\|_{L^2(NB_1)}
  + \|\tilde{f}\|_{L^p(NB_1)} + \|\tilde{F}\|_{L^q(NB_1)}\big\},
\end{eqnarray*}
where $C$ is the same constant as in $\eqref{f:4.18}$. This implies
\begin{eqnarray*}
\|\nabla u_\varepsilon\|_{L^p(B_r)} \leq C\big\{r^{-1+\frac{d}{p} -\frac{d}{2}}\|u_\varepsilon\|_{L^2(NB_{r})}
+ \|f\|_{L^p(NB_{r})} + r^{1+\frac{d}{p}-\frac{d}{q}}\|F\|_{L^q(NB_{r})} \big\}.
\end{eqnarray*}

Finally, for any $B$ with $0<r\leq 1$, we choose the small ball with $r/N$ radius to cover $B_r$. Hence we have
\begin{eqnarray*}
\|\nabla u_\varepsilon\|_{L^p(B_r)}\leq C\big\{r^{-1+\frac{d}{p} -\frac{d}{2}}\|u_\varepsilon\|_{L^2(B_{2r})}
+ \|f\|_{L^p(B_{2r})} + r^{1+\frac{d}{p}-\frac{d}{q}}\|F\|_{L^q(B_{2r})} \big\},
\end{eqnarray*}
and this gives $\eqref{pri:4.7}$. We complete the proof.
\qed
\end{pf}

\begin{remark}\label{rm:4.2}
\emph{ Here we introduce two elementary interpolation inequalities used in the above proof. Let $u\in W^{1,p}(\Omega;\mathbb{R}^m)$ with
$2\leq p<\infty$,
then for any $\delta>0$,
there exists a constant $C_\delta$ depending on $\delta,p,d,m$ and $\Omega$, such that
\begin{equation}\label{f:4.7}
 \|u\|_{L^p(\Omega)} \leq \delta\|\nabla u\|_{L^p(\Omega)} + C_\delta\|u\|_{L^2(\Omega)}.
\end{equation}
The estimate $\eqref{f:4.7}$ can be easily derived by contradiction argument (or see \cite{RAA}). As a result, we have
\begin{equation}\label{f:4.6}
\|u\|_{L^{p}(\Omega)} \leq C\|u\|_{W^{1,q}(\Omega)} \leq C\|\nabla u\|_{L^q(\Omega)} + C\|u\|_{L^2(\Omega)}
\end{equation}
for $1 \leq p < \infty$ and  $q = \frac{pd}{p+d}$, where $C$ depends on $p,d,m$ and $\Omega$. }
\end{remark}

\begin{cor}\label{cor:4.2}
Suppose that the coefficients of $\mathcal{L}_\varepsilon$ satisfy the same conditions as in Theorem $\ref{thm:4.1}$.
Let $p>d$ and $\sigma = 1 - d/p$.
Assume that $u_\varepsilon\in H^1_0(\Omega;\mathbb{R}^m)$ is a weak solution of $\mathcal{L}_\varepsilon(u_\varepsilon) = \emph{div}(f) + F$ in $\Omega$,
where $f\in L^p(\Omega;\mathbb{R}^{md})$
and $F \in L^q(\Omega;\mathbb{R}^m)$ with $q = \frac{pd}{p+d}>\frac{d}{2}$.
Then we have
\begin{equation}\label{pri:4.9}
\|u_\varepsilon\|_{C^{0,\sigma}(B)}
 \leq Cr^{-\sigma}\left\{\Big(\dashint_{2B} |u_\varepsilon|^2 dx\Big)^{\frac{1}{2}}
 + r\Big(\dashint_{2B} |f|^p dx\Big)^{\frac{1}{p}} + r^2\Big(\dashint_{2B} |F|^q dx\Big)^{\frac{1}{q}}\right\}
\end{equation}
for any $B\subset 2B\subset\Omega$ with $0<r\leq 1$. In particular, for any $s>0$,
\begin{equation}\label{pri:4.8}
\|u_\varepsilon\|_{L^\infty(B)}
\leq C\left\{\Big(\dashint_{2B} |u_\varepsilon|^s dx\Big)^{\frac{1}{s}}
+ r\Big(\dashint_{2B} |f|^p dx\Big)^{\frac{1}{p}} + r^2\Big(\dashint_{2B} |F|^q dx\Big)^{\frac{1}{q}}\right\},
\end{equation}
where $C$ depends only on $\mu, \omega, \kappa, \lambda, p, m, d$.
\end{cor}

\begin{pf}
 Assume $r=1$. It follows from the Sobolev embedding theorem and Remark $\ref{rm:4.2}$ that
 \begin{eqnarray*}
  \|u_\varepsilon\|_{C^{0,\sigma}(B)} \leq C \|u_\varepsilon\|_{W^{1,p}(B)}
  \leq C\|\nabla u_\varepsilon\|_{L^p(B)} + C\|u_\varepsilon\|_{L^2(B)}.
  \end{eqnarray*}
 Then it follows from $\eqref{pri:4.7}$ and rescaling arguments that
  \begin{eqnarray*}
  \|u_\varepsilon\|_{C^{0,\sigma}(B)}\leq Cr^{-\sigma}\left(\dashint_{2B} |u_\varepsilon|^2 dx\right)^{1/2}
  + Cr^{1-\sigma}\left\{\Big(\dashint_{2B} |f|^p dx\Big)^{\frac{1}{p}} + r\Big(\dashint_{2B} |F|^q dx\Big)^{\frac{1}{q}}\right\},
 \end{eqnarray*}
 where $\sigma = 1-d/p$. Moreover, for any $x\in B$ we have
 \begin{eqnarray*}
  |u_\varepsilon(x)| \leq |u_\varepsilon(x)-\bar{u}_\varepsilon| + |\bar{u}_\varepsilon|
  \leq \big[u_\varepsilon\big]_{C^{0,\sigma}(B)}r^{\sigma} + \left(\dashint_{B}|u_\varepsilon(y)|^2dy\right)^{1/2}
 \end{eqnarray*}
 from H\"{o}lder's inequality. This gives
 \begin{equation*}
 \|u_\varepsilon\|_{L^\infty(B)}
 \leq C \left(\dashint_{2B} |u_\varepsilon|^2 dx\right)^{1/2}
 + Cr\left\{\Big(\dashint_{2B} |f|^p dx\Big)^{\frac{1}{p}} + r\Big(\dashint_{2B} |F|^q dx\Big)^{\frac{1}{q}}\right\}.
 \end{equation*}
 Moreover, by the iteration method (see \cite[pp.184]{MGLM}),
 we have $\eqref{pri:4.8}$.
 \qed
\end{pf}


\begin{lemma}\label{lemma:4.2}
Suppose that the coefficients of $\mathcal{L}_\varepsilon$ satisfy the same assumptions as in Theorem $\ref{thm:4.1}$.
Let $1<p<\infty$, $f = (f_i^\alpha)\in L^p(\Omega;\mathbb{R}^{md})$.
Then there exists a unique $u_\varepsilon\in W^{1,p}_0(\Omega;\mathbb{R}^m)$ such that
$\mathcal{L}_\varepsilon(u_\varepsilon) = \emph{div}(f)$ in $\Omega$ and
$u_\varepsilon = 0$ on $\partial\Omega$, whenever
$\lambda \geq \lambda_0$. 
Moreover, the solution satisfies the uniform estimate
\begin{equation}\label{pri:4.2}
  \|\nabla u_\varepsilon\|_{L^p(\Omega)} \leq C \|f\|_{L^p(\Omega)},
\end{equation}
where $\lambda_0$ is given in Lemma $\ref{lemma:2.2}$, and $C$ depends on $\mu, \omega, \kappa, \lambda, p, d, m$ and $\Omega$.
\end{lemma}

\begin{pf}
In the case of $p=2$, it follows from Theorem $\ref{thm:2.1}$ that there exists a unique solution
$u_\varepsilon\in H^1_0(\Omega;\mathbb{R}^m)$ satisfying the uniform estimate $\|\nabla u_\varepsilon\|_{L^2(\Omega)}\leq C\|f\|_{L^2(\Omega)}$.
For $p>2$, the uniqueness and existence of the weak solution is reduced to the case of $p=2$. We rewrite the original systems as
\begin{equation*}
L_{\varepsilon}(u_\varepsilon) = \text{div}(f+V_\varepsilon u_\varepsilon)
- B_\varepsilon\nabla u_\varepsilon - (c_\varepsilon+\lambda) u_\varepsilon.
\end{equation*}
Applying $\eqref{pri:2.3}$, $\eqref{pri:4.1}$ and Sobolev's inequality, we obtain
\begin{eqnarray}\label{f:4.17}
 \|\nabla u_\varepsilon\|_{L^{p_{k_{0}}}(\Omega)}
 &\leq& C\{ \|f\|_{L^{p_{k_0}}(\Omega)} + \|u_\varepsilon\|_{L^{p_{k_0}}(\Omega)} + \|\nabla u_\varepsilon\|_{L^{p_{{k_{0}}-1}}(\Omega)} \}
 \leq C\{ \|f\|_{L^{p_{k_{0}}}(\Omega)} + \|\nabla u_\varepsilon\|_{L^{p_{{k_{0}}-1}}(\Omega)} \}  \nonumber\\
 &\leq& C\{ \|f\|_{L^{p_{k_{0}}}(\Omega)} + \|f\|_{L^{2}(\Omega)}\} \leq C\|f\|_{L^{p_{k_{0}}}(\Omega)},
\end{eqnarray}
where $p_{k_0} = \frac{2d}{d-2k_0}$, and $k_0$ is a positive integer such that $k_0< d/2\leq k_0+1$.
We claim that $\eqref{f:4.17}$ holds for any $p\in[2,p_{k_0}]$.
Indeed, let $T(f) = \nabla u_\varepsilon$, then together with $\|T\|_{L^2\to L^2}\leq C$ and $\|T\|_{L^{p_{k_0}}\to L^{p_{k_0}}}\leq C$,
the Marcinkiewicz interpolation theorem gives $\|T\|_{L^p\to L^p}\leq C$ (see \cite{S}).
Moreover, for any $p>p_{k_0}$, we still have
$\|\nabla u_\varepsilon\|_{L^p(\Omega)} \leq C\{\|f\|_{L^p(\Omega)} + \|\nabla u_\varepsilon\|_{L^q(\Omega)}\}\leq C\|f\|_{L^p(\Omega)}$,
since $q<d$.


By the duality argument, we can derive $\eqref{pri:4.2}$ for $p\in(1,2)$.
Let $h = (h_i^\beta)\in C^1_0(\Omega;\mathbb{R}^{md})$, and $v_\varepsilon$ be the weak solution to
$\mathcal{L}_\varepsilon^*(v_\varepsilon) = \text{div}(h)$ in $\Omega$ and $v_\varepsilon =0$ on $\partial\Omega$.
Hence, in view of the previous result,
we have $\|\nabla v_\varepsilon\|_{L^{p^\prime}(\Omega)} \leq C\|h\|_{L^{p^\prime}(\Omega)}$ for any $p^\prime > 2$. Moreover,
if $f\in C^1_0(\Omega;\mathbb{R}^{md})$, there exists the weak solution $u_\varepsilon\in H^1_0(\Omega;\mathbb{R}^{m})$
to the original systems.
Then it follows from Remark
$\ref{rm:2.3}$ that
\begin{eqnarray*}
\int_\Omega \nabla u_\varepsilon h dx  =   -\int_\Omega u_\varepsilon \mathcal{L}_\varepsilon^*(v_\varepsilon) dx
 =  -\int_\Omega \mathcal{L}_\varepsilon(u_\varepsilon)v_\varepsilon dx  = \int_\Omega f\nabla v_\varepsilon dx.
\end{eqnarray*}
This gives $\|\nabla u_\varepsilon\|_{L^p(\Omega)}\leq C\|f\|_{L^p(\Omega)}$, where $p = p^\prime/(p^\prime -1)$.
By the density argument, we can verify the existence of solutions in $W^{1,p}_0(\Omega)$ for general $f\in L^p(\Omega;\mathbb{R}^{md})$,
as well as the uniqueness for $1<p<2$. The proof is complete.
\qed

\end{pf}

\begin{lemma}\label{lemma:4.3}
Suppose that the coefficients of $\mathcal{L}_\varepsilon$ satisfy the same conditions as in Lemma $\ref{lemma:4.2}$.
Let $1<p<\infty$.
Then for any $F\in L^q(\Omega;\mathbb{R}^m)$,
where $q = pd/(p+d)$ if $p> d/(d-1)$, and $q>1$ if $1<p\leq d/(d-1)$, there exists a unique solution $u_\varepsilon\in W^{1,p}_0(\Omega;\mathbb{R}^m)$
to $\mathcal{L}_\varepsilon(u_\varepsilon) = F$ in $\Omega$ and $u_\varepsilon =0$ on $\partial\Omega$, satisfying
the uniform estimate
\begin{equation}\label{pri:4.3}
 \|\nabla u_\varepsilon\|_{L^p(\Omega)} \leq C\|F\|_{L^{q}(\Omega)},
\end{equation}
where $C$ depends only on $\mu, \omega, \kappa,\lambda, p, q, d,m$ and $\Omega$.
\end{lemma}

\begin{pf}
We prove this lemma by the duality argument. The uniqueness is clearly contained in Lemma $\ref{lemma:4.2}$, and the existence of
the solution $u_\varepsilon$ follows from the density and Theorem $\ref{thm:2.1}$.
The rest of the proof is to establish $\eqref{pri:4.3}$.

Consider the dual problem for any $f\in C^1_0(\Omega;\mathbb{R}^{md})$, there exists the unique $v_\varepsilon\in H^1_0(\Omega;\mathbb{R}^m)$ to
$\mathcal{L}_\varepsilon^*(v_\varepsilon)=\text{div}(f)$ in $\Omega$ and $v_\varepsilon = 0$ on $\partial\Omega$. Note that
it follows from Lemma $\ref{lemma:4.2}$ that
$\|\nabla v_\varepsilon\|_{L^{p^\prime}(\Omega)} \leq C\|f\|_{L^{p^\prime}(\Omega)}$. Then we have
\begin{eqnarray*}
\int_\Omega \nabla u_\varepsilon f dx = -\int_\Omega u_\varepsilon\mathcal{L}_\varepsilon^*(v_\varepsilon)dx
= - \int_\Omega \mathcal{L}_\varepsilon(u_\varepsilon)v_\varepsilon dx = -\int_\Omega Fv_\varepsilon dx,
\end{eqnarray*}
and
\begin{eqnarray*}
 \left|\int_\Omega \nabla u_\varepsilon f dx\right| &\leq& \|F\|_{L^{q}(\Omega)}\|v_\varepsilon\|_{L^{q^\prime}(\Omega)}
 \leq C \|F\|_{L^{q}(\Omega)}\|\nabla v_\varepsilon\|_{L^{p^\prime}(\Omega)}
 \leq C  \|F\|_{L^{q}(\Omega)}\|f\|_{L^{p^\prime}(\Omega)}.
\end{eqnarray*}
Note that
$\frac{1}{q^\prime} = \frac{1}{p^\prime} - \frac{1}{d}$ if $p^\prime < d$, $1< q^\prime < \infty$ if $p^\prime = d$,
and $q^\prime = \infty$ if $p^\prime > d$.
In other words, $q = \frac{pd}{p+d}$ if $p>\frac{d}{d-1}$, and $q>1$ if $1<p\leq \frac{d}{d-1}$.
Finally we obtain $\|\nabla u_\varepsilon\|_{L^p(\Omega)} \leq C\|F\|_{L^{q}(\Omega)}$.
\qed
\end{pf}

\begin{flushleft}
\textbf{Proof of Theorem \ref{thm:1.1}}\textbf{.}\quad
In the case of $g = 0$, we write $v_\varepsilon = u_{\varepsilon,1} + u_{\varepsilon,2}$,
where $u_{\varepsilon,1}$ and $u_{\varepsilon,2}$ are  the solutions in Lemma $\ref{lemma:4.2}$ and $\ref{lemma:4.3}$, respectively.
Then we have
\begin{eqnarray}\label{f:4.4}
\|\nabla v_\varepsilon\|_{L^p(\Omega)} &\leq& \|\nabla u_{\varepsilon,1}\|_{L^p(\Omega)} + \|\nabla u_{\varepsilon,2}\|_{L^p(\Omega)}
\leq C\{\|f\|_{L^p(\Omega)} + \|F\|_{L^q(\Omega)} \}.
\end{eqnarray}
\end{flushleft}

For $g\not=0$, consider the homogeneous Dirichlet problem
$\mathcal{L}_\varepsilon(w_\varepsilon) = 0$ in $\Omega$ and $w_\varepsilon = g$ on $\partial\Omega$,
where $g\in B^{1-1/p,p}(\partial\Omega;\mathbb{R}^m)$. By the properties of boundary Besov space, there exists $G\in W^{1,p}(\Omega;\mathbb{R}^m)$
such that
$G = g$ on $\partial\Omega$ and $\|G\|_{W^{1,p}(\Omega)}\leq C\|g\|_{B^{1-1/p,p}(\partial\Omega)}$. Let $h_\varepsilon = w_\varepsilon-G$, we have
\begin{equation*}
\mathcal{L}_\varepsilon(h_\varepsilon) = \text{div}\big(A_\varepsilon\nabla G + V_\varepsilon G\big)
- B_\varepsilon\nabla G  - (c_\varepsilon +\lambda)G~~\text{in}~\Omega, \qquad h_\varepsilon = 0 ~~\text{on}~\partial\Omega.
\end{equation*}
Recall the case of $g = 0$, in which there exists the unique weak solution $h_\varepsilon\in W^{1,p}_0(\Omega;\mathbb{R}^m)$,
satisfying the uniform estimate
$\|\nabla h_\varepsilon\|_{L^p(\Omega)}
\leq C \|G\|_{W^{1,p}(\Omega)} \leq C \|g\|_{B^{1-1/p,p}(\partial\Omega)}$
for $1<p<\infty$. This implies
\begin{equation}\label{f:4.3}
\|\nabla w_\varepsilon\|_{L^p(\Omega)} \leq \|\nabla h_\varepsilon\|_{L^p(\Omega)} + \|\nabla G\|_{L^p(\Omega)}
\leq C \|g\|_{B^{1-1/p,p}(\partial\Omega)}.
\end{equation}
Finally, let $u_\varepsilon = v_\varepsilon + w_\varepsilon$.
Combining $\eqref{f:4.4}$ and $\eqref{f:4.3} $, we have
\begin{eqnarray*}
\|\nabla u_\varepsilon\|_{L^p(\Omega)}
&\leq& C\{\|f\|_{L^p(\Omega)} + \|F\|_{L^q(\Omega)} + \|g\|_{B^{1-1/p,p}(\partial\Omega)} \},
\end{eqnarray*}
where $C$ depends only on $\mu,\omega,\kappa,\lambda,p,q,d,m$ and $\Omega$. We complete the proof.
\qed

\begin{cor}\label{cor:4.1}
Suppose that the coefficients of $\mathcal{L}_\varepsilon$ satisfy the same conditions as in Theorem $\ref{thm:1.1}$.
Set $d<p< \infty$ and $\sigma = 1 - d/p$.
Let $f = (f_i^\alpha)\in L^p(\Omega;\mathbb{R}^{md})$, $F\in L^q(\Omega;\mathbb{R}^m)$ with $q = \frac{pd}{p+d}$,
and $g\in C^{0,1}(\partial\Omega;\mathbb{R}^m)$.
Then the unique solution $u_\varepsilon$ to $\eqref{pde:1.1}$
satisfies the uniform estimate
\begin{equation}\label{pri:4.4}
\|u_\varepsilon\|_{C^{0,\sigma}(\Omega)} \leq C\big\{\|f\|_{L^p(\Omega)} + \|F\|_{L^q(\Omega)} + \|g\|_{C^{0,1}(\partial\Omega)}\big\},
\end{equation}
where $C$ depends only on $\mu,\omega,\kappa,\lambda, p, q, d,m,$ and $\Omega$.
\end{cor}

\begin{pf}
Due to the extension theorem (see \cite[pp.136]{DGNT}), there exists an extension function $G\in C^{0,1}{(\Omega;\mathbb{R}^m)}$
such that $G = g$ on $\partial\Omega$ and
$\|G\|_{C^{0,1}(\Omega)} \leq C\|g\|_{C^{0,1}(\partial\Omega)}$.
This also implies $\|G\|_{W^{1,p}(\Omega)}\leq C\|g\|_{C^{0,1}(\partial\Omega)}$ for any $p\geq 1$.
Let $v_\varepsilon, w_\varepsilon$ be the weak solutions to the following Dirichlet problems:
\begin{equation}\label{pde:3.1}
  (\text{i})\left\{ \begin{array}{rcll}
  \mathcal{L}_\varepsilon(v_\varepsilon) &=& \text{div}(f) + F &~\text{in}~\Omega, \\
   v_\varepsilon &=& 0 &~\text{on}~\partial\Omega,
  \end{array} \right.
  \qquad
    (\text{ii})\left\{ \begin{array}{rcll}
  \mathcal{L}_\varepsilon(w_\varepsilon) &=& 0 &~\text{in}~\Omega, \\
   w_\varepsilon &=& g &~\text{on}~\partial\Omega,
  \end{array} \right.
\end{equation}
respectively. For $(\text{i})$, it follows from the Sobolev embedding theorem and Theorem $\ref{thm:1.1}$ that
$
\|v_\varepsilon\|_{C^{0,\sigma}(\Omega)} \leq C\|\nabla v_\varepsilon\|_{L^p(\Omega)}
\leq C\{\|f\|_{L^p(\Omega)} +  \|F\|_{L^q(\Omega)}\}
$.
For $(\text{ii})$, by setting $h_\varepsilon = w_\varepsilon -G$,
we have $\mathcal{L}_\varepsilon(h_\varepsilon) = -\mathcal{L}_\varepsilon(G)$ in
$\Omega$ and $h_\varepsilon = 0$ on $\partial\Omega$. Therefore
$\|h_\varepsilon\|_{C^{0,\sigma}(\Omega)}\leq C\|G\|_{W^{1,p}(\Omega)}\leq C\|g\|_{C^{0,1}(\partial\Omega)}$, which implies
$
 \|w_\varepsilon\|_{C^{0,\sigma}(\Omega)} \leq C \|g\|_{C^{0,1}(\partial\Omega)}.
$
Let $u_\varepsilon = v_\varepsilon + w_\varepsilon$. Combining the estimates related to $v_\varepsilon$ and $w_\varepsilon$, we derive the estimate $\eqref{pri:4.4}$.
\qed
\end{pf}

\begin{remark}
\emph{ Assume the same conditions as in Corollary $\ref{cor:4.1}$,
let $u_\varepsilon$ be a weak solution to $\eqref{pde:1.1}$.
Then by the same argument as in the proof of Theorem $\ref{thm:4.1}$,
we obtain the near boundary H\"{o}lder estimate
\begin{equation}\label{pri:4.12}
\|u_\varepsilon\|_{C^{0,\sigma}(D(r))}
\leq Cr^{-\sigma}\left\{\Big(\dashint_{D(2r)}|u_\varepsilon|^2dx\Big)^{\frac{1}{2}}
+ r \|g\|_{C^{0,1}(\Delta(2r))} + r\Big(\dashint_{D(2r)}|f|^pdx\Big)^{\frac{1}{p}}
+ r^2 \Big(\dashint_{D(2r)}|F|^qdx\Big)^{\frac{1}{q}} \right\},
\end{equation}
where $\sigma = 1 - d/p$, and $C$ depends on $\mu,\omega,\kappa,\lambda,p,d,m$ and $\Omega$. }
\end{remark}

\begin{remark}
\emph{ In the following, we frequently use the abbreviated writing like
$\dashint_{\Omega} F(\cdot,y) = \dashint_{\Omega} F(x,y) dx$.}
\end{remark}

\begin{thm}\label{thm:4.2}
Suppose the same conditions on $\mathcal{L}_\varepsilon$ as in Theorem $\ref{thm:1.1}$, and $\lambda\geq \lambda_0$.
Then there exists a unique Green's matrix $\mathcal{G}_\varepsilon:\Omega \times \Omega \to\mathbb{R}^{m^2}\cup\{\infty\}$,
such that
$\mathcal{G}_\varepsilon(\cdot,y)\in H^1(\Omega\setminus B_r(y);\mathbb{R}^{m^2})\cap W_0^{1,s}(\Omega;\mathbb{R}^{m^2})$
with $s\in[1,\frac{d}{d-1})$ for each $y\in \Omega$ and $r>0$, and
\begin{equation}\label{pri:4.15}
\mathrm{B}^*_\varepsilon\big[\mathcal{G}_\varepsilon^\gamma(\cdot,y),\phi\big]
= \phi^\gamma(y), \qquad \forall~ \phi\in W_0^{1,p}(\Omega;\mathbb{R}^m),~p>d.
\end{equation}
Particularly, for any $F\in L^q(\Omega;\mathbb{R}^m)$ with $q>d/2$,
\begin{equation}\label{pri:4.16}
 u^\gamma_\varepsilon(y) = \int_\Omega \mathcal{G}_\varepsilon^{\gamma\alpha}(x,y)F^\alpha(x)dx
\end{equation}
satisfies
$\mathcal{L}_\varepsilon(u_\varepsilon) = F$ in $\Omega$ and $u_\varepsilon = 0$ on $\partial\Omega$.
Moreover, let $^{*}\mathcal{G}_\varepsilon(\cdot,x)$ be the adjoint Green's matrix of $\mathcal{G}_\varepsilon(\cdot,y)$,
then $\mathcal{G}_\varepsilon(x,y) = [~^{*}\mathcal{G}_\varepsilon(y,x)~]^*$
and for any $\sigma,\sigma^\prime\in(0,1)$, the following estimates
\begin{eqnarray}\label{pri:4.17}
 |\mathcal{G}_\varepsilon(x,y)| \leq \frac{C}{|x-y|^{d-2}}\min\Big\{1,\frac{d^{\sigma}_x}{|x-y|^{\sigma}},
 \frac{d^{\sigma^\prime}_y}{|x-y|^{\sigma^\prime}},\frac{d^{\sigma}_xd^{\sigma^\prime}_y}{|x-y|^{\sigma+\sigma^\prime}}\Big\}
\end{eqnarray}
hold for any $x,y\in \Omega$ and $x\not=y$, where $d_x = \emph{dist}(x,\partial\Omega)$,
and $C$ depends only on $\mu,\omega,\kappa,\lambda,d,m$ and $\Omega$.
\end{thm}

\begin{lemma}[Approximating Green's matrix]\label{lemma:4.5}
Assume the same conditions as in Theorem $\ref{thm:4.2}$. Define the approximating Green's matrix
$G_{\rho,\varepsilon}(\cdot,y)$ as
\begin{equation}\label{def:1}
 \mathrm{B}_\varepsilon^*[G_{\rho,\varepsilon}^{\gamma}(\cdot,y),u] = \dashint_{\Omega_\rho(y)} u^\gamma dx, \qquad \forall ~u\in H_0^1(\Omega;\mathbb{R}^m),
\end{equation}
where $1\leq \gamma\leq m$, and $\Omega_\rho(y) = \Omega\cap B_\rho(y)$.
Then if $|x-y|<d_y/2$, we have the uniform estimate
\begin{eqnarray}\label{pri:4.13}
 |G_{\rho,\varepsilon}^\gamma(x,y)| \leq \frac{C}{|x-y|^{d-2}},\qquad \forall~\rho < |x-y|/4,
\end{eqnarray}
where $C$ depends only on $\mu,\omega,\kappa,\lambda,d,m$ and $\Omega$.
Moreover, for any $s\in [1,\frac{d}{d-1})$, we have
\begin{equation}\label{pri:4.14}
 \sup_{\rho>0}\big\|G_{\rho,\varepsilon}^\gamma(\cdot,y)\big\|_{W_0^{1,s}(\Omega)} \leq C(\mu,\omega,\kappa,\lambda,d,m,s,\Omega,d_y).
\end{equation}
\end{lemma}

\begin{pf}
First of all, we show $G_{\rho,\varepsilon}(x,y)= [G_{\rho,\varepsilon}^{\gamma\theta}(x,y)]$ is well defined.
Let $I(u) = \dashint_{\Omega_\rho(y)} u^\gamma dx$, then $I\in H^{-1}(\Omega;\mathbb{R}^m)$ and
$|I(u)|\leq C|\Omega_\rho(y)|^{-1/2^*}\|u\|_{H^1_0(\Omega)}$ with $2^* = \frac{2d}{d-2}$.
It follows from Theorem $\ref{thm:2.1}$ that there exists a unique
$G_{\rho,\varepsilon}^{\gamma}(\cdot,y)\in H^1_0(\Omega;\mathbb{R}^m)$ satisfying $\eqref{def:1}$ and
\begin{equation}\label{f:4.8}
 \|\nabla G_{\rho,\varepsilon}^{\gamma}(\cdot,y)\|_{L^2(\Omega)} \leq C|\Omega_\rho(y)|^{-\frac{1}{2^*}}.
\end{equation}
For any $F\in C^\infty_0(\Omega;\mathbb{R}^m)$, consider
$\mathcal{L}_\varepsilon(u_\varepsilon) = F$ in $\Omega$ and $u_\varepsilon = 0$ on $\partial\Omega$.
There exists the unique solution $u_\varepsilon\in H^1_0(\Omega;\mathbb{R}^m)$ such that
\begin{equation}\label{f:4.16}
 \int_\Omega F G_{\rho,\varepsilon}^\gamma(\cdot,y) = \mathrm{B}_\varepsilon[u_\varepsilon, G_{\rho,\varepsilon}^\gamma(\cdot,y)]
 = \mathrm{B}_\varepsilon^*[G_{\rho,\varepsilon}^\gamma(\cdot,y),u_\varepsilon] = \dashint_{\Omega_{\rho}(y)} u^\gamma_\varepsilon.
\end{equation}
Suppose $\mathrm{supp}(F) \subsetneqq B\subsetneqq\Omega$, where $B= B_R(y)$. 
Then it follows from $\eqref{pri:2.3}$, $\eqref{pri:4.8}$ and $\eqref{f:4.16}$ that
\begin{eqnarray*}
\Big|\int_{\Omega} F G_{\rho,\varepsilon}^\gamma(\cdot,y)\Big|
\leq  \|u_\varepsilon\|_{L^\infty(1/4B)}
\leq C\Big[\Big(\dashint_{1/2B}|u_\varepsilon|^2\Big)^{\frac{1}{2}} + R^2\Big(\dashint_{1/2B}|F|^p\Big)^{\frac{1}{p}}\Big]
\leq CR^2\Big(\dashint_{B}|F|^p\Big)^{\frac{1}{p}}
\end{eqnarray*}
for any $\rho < R/4$ and $p > d/2$. This implies
\begin{equation*}
 \Big(\dashint_{B}\big|G_{\rho,\varepsilon}^\gamma(\cdot,y)\big|^q\Big)^{\frac{1}{q}} \leq C R^{2-d},
 \qquad \forall ~R\leq d_y,\quad\forall~ q\in\big[1,\frac{d}{d-2}\big).
\end{equation*}

Now we turn to $\eqref{pri:4.13}$. Set $r=|x-y|$, and $r\leq d_y/2$.
In view of $\eqref{def:1}$, $G_{\rho,\varepsilon}(x,y)$ actually satisfies
$\mathcal{L}_\varepsilon\big[G_{\rho,\varepsilon}^\gamma(\cdot,y)\big] =0$ in $B_{\frac{r}{2}}(y)\setminus B_{\frac{r}{4}}(y)$.
By using $\eqref{pri:4.8}$ again, we obtain
\begin{equation*}
 \big|G_{\rho,\varepsilon}^\gamma(x,y)\big| \leq C \dashint_{B_{\frac{r}{2}}(x)} \big|G_{\rho,\varepsilon}^\gamma(\cdot,y)\big|
 \leq C \dashint_{B_{2r}(y)} \big|G_{\rho,\varepsilon}^\gamma(\cdot,y)\big| \leq C |x-y|^{2-d}
\end{equation*}
for any $\rho < |x-y|/4$, where $C$ depends only on $\mu,\omega,\kappa,\lambda,d,m$ and $\Omega$.

Then we will prove $\eqref{pri:4.14}$. Step one, we verify the following estimates,
\begin{eqnarray}\label{f:4.10}
\int_{\Omega\setminus B(y,R)}|\nabla G_{\rho,\varepsilon}^\gamma(\cdot,y)|^2  \leq CR^{2-d}, \quad
\int_{\Omega\setminus B(y,R)}|G_{\rho,\varepsilon}^\gamma(\cdot,y)|^{2^*}  \leq CR^{-d}, \qquad \forall~\rho>0,\quad\forall~R<d_y/4.
\end{eqnarray}
On the one hand, let $\varphi\in C^1_0(\Omega)$ be a cut-off function satisfying $\varphi\equiv 0$ on $B(y,R)$, $\varphi \equiv 1$ outside $B(y,2R)$,
and $|\nabla\varphi| \leq C/R$. Choose $u = \varphi^2 G_{\rho,\varepsilon}^\gamma(\cdot,y)$ in $\eqref{def:1}$ and $\lambda\geq \lambda_0$.
It follows from $\eqref{f:4.9}$ and $\eqref{pri:4.13}$ that
\begin{equation}\label{f:4.30}
 \int_\Omega \varphi^2 |\nabla G_{\rho,\varepsilon}^\gamma(\cdot,y)|^2
 \leq C\int_\Omega |\nabla\varphi|^2|G_{\rho,\varepsilon}^\gamma(\cdot,y)|^2
 \leq \frac{C}{R^2}\int_{2B\setminus B} |x-y|^{2(2-d)}
 \leq CR^{2-d}, \qquad \forall~\rho <R/4.
\end{equation}
On the other hand, it follows from $\eqref{f:4.8}$ that
\begin{eqnarray*}
\int_{\Omega\setminus B(y,R)}|\nabla G_{\rho,\varepsilon}^\gamma(\cdot,y)|^2
\leq  \int_{\Omega}|\nabla G_{\rho,\varepsilon}^\gamma(\cdot,y)|^2
\leq CR^{2-d}, \qquad \forall~\rho \geq R/4.
\end{eqnarray*}
Thus we have the first inequality of $\eqref{f:4.10}$ for all $\rho>0$.

For the second estimate in $\eqref{f:4.10}$,
we observe
\begin{equation}\label{f:4.31}
 \int_{\Omega}\big|\varphi G_{\rho,\varepsilon}^\gamma(\cdot,y)\big|^{2^*}
 \leq C\Big(\int_\Omega \big|\nabla\big(\varphi G_{\rho,\varepsilon}^\gamma(\cdot,y)\big)\big|^2 \Big)^{\frac{d}{d-2}}
 \leq C\Big(\int_\Omega \big|\nabla\varphi G_{\rho,\varepsilon}^\gamma(\cdot,y)\big|^2
 + \big|\varphi\nabla G_{\rho,\varepsilon}^\gamma(\cdot,y)\big|^2\Big)^{\frac{d}{d-2}}
 \leq CR^{-d}
\end{equation}
for any $\rho < R/4$, where we use Sobolev's inequality in the first inequality and $\eqref{f:4.30}$ in the last inequality. We remark that
the constant $C$ does not involve $R$. In the case of $\rho \geq R/4$,
since $G_{\rho,\varepsilon}^\gamma(\cdot,y) = 0$ on $\partial\Omega$, we have
\begin{equation*}
\int_{\Omega\setminus B(y,R)}|G_{\rho,\varepsilon}^\gamma(\cdot,y)|^{2^*}
\leq \int_{\Omega}|G_{\rho,\varepsilon}^\gamma(\cdot,y)|^{2^*}
\leq C\Big(\int_\Omega |\nabla G_{\rho,\varepsilon}^\gamma(\cdot,y)|^{2}\Big)^{\frac{d}{d-2}}
\leq C R^{-d},
\end{equation*}
where we use Sobolev's inequality in the second inequality and $\eqref{f:4.8}$ in the last inequality.
This together with $\eqref{f:4.31}$ leads to
\begin{equation*}
\int_{\Omega\setminus B(y,R)}|G_{\rho,\varepsilon}^\gamma(\cdot,y)|^{2^*}  \leq CR^{-d}, \qquad \forall~\rho>0,\quad\forall~R<d_y/4.
\end{equation*}

We now address ourselves to the uniform estimates of $G_{\rho,\varepsilon}^\gamma(\cdot,y)$ and $\nabla G_{\rho,\varepsilon}^\gamma(\cdot,y)$
with respect to parameter $\rho$. In the case of $t>(d_y/4)^{1-d}$, we obtain
\begin{eqnarray}\label{f:4.14}
\big|\big\{x\in\Omega: |\nabla G_{\rho,\varepsilon}^\gamma(\cdot,y)|>t\big\}\big|
\leq CR^{d} + t^{-2}\int_{\Omega\setminus B(y,R)}|\nabla G_{\rho,\varepsilon}^\gamma(\cdot,y)|^2
\leq C t^{-\frac{d}{d-1}},\qquad \forall ~\rho>0.
\end{eqnarray}
For $t>(d_y/4)^{2-d}$, it follows that
\begin{eqnarray}\label{f:4.15}
\big|\big\{x\in\Omega: |G_{\rho,\varepsilon}^\gamma(\cdot,y)|>t\big\}\big|
\leq C t^{-\frac{d}{d-2}},\qquad \forall ~\rho>0.
\end{eqnarray}
Then in view of $\eqref{f:4.14}$ and $\eqref{f:4.15}$, we have
\begin{eqnarray*}
 \int_{\Omega}|G_{\rho,\varepsilon}^\gamma(\cdot,y)|^s
 \leq C d_y^{s(2-d)} + C\int_{(d_y/4)^{2-d}}^\infty t^{s-1}\cdot t^{-\frac{d}{d-2}} dt
 \leq C\big[d_y^{s(2-d)}+d_y^{s(2-d)+d}\big]
\end{eqnarray*}
for $s\in[1,\frac{d}{d-2})$, and
\begin{eqnarray*}
 \int_{\Omega}|\nabla G_{\rho,\varepsilon}^\gamma(\cdot,y)|^s
 \leq C\big[d_y^{s(1-d)}+d_y^{s(1-d)+d}\big]
\end{eqnarray*}
for $s\in [1,\frac{d}{d-1})$,
where $C$ depends only on $\mu,\omega,\kappa,\lambda,d,m,s$ and $\Omega$. Combining the two inequalities above,
we have $\eqref{pri:4.14}$, and the proof is complete.
\qed
\end{pf}

\begin{flushleft}
\textbf{Proof of Theorem \ref{thm:4.2}}\textbf{.}\quad From the uniform estimate $\eqref{pri:4.14}$, it follows that there exist a subsequence of
$\{G_{\rho_n,\varepsilon}^\gamma(\cdot,y)\}_{n=1}^\infty$ and $\mathcal{G}_{\varepsilon}^\gamma(\cdot,y)$ such that for any $s\in(1,\frac{d}{d-1})$,
\begin{eqnarray}
 G_{\rho_n,\varepsilon}^\gamma(\cdot,y)~\rightharpoonup~ \mathcal{G}_{\varepsilon}^\gamma(\cdot,y)
 \quad \text{weakly in}~ W^{1,s}_0(\Omega;\mathbb{R}^m)
 \quad \text{as}~~n \to \infty.
\end{eqnarray}
\end{flushleft}
Hence, we have
\begin{equation*}
\mathrm{B}_\varepsilon^*[\mathcal{G}_\varepsilon^\gamma(\cdot,y),\phi]
= \lim_{n\to\infty}  \mathrm{B}_\varepsilon^*[G_{\rho_n,\varepsilon}^\gamma(\cdot,y),\phi]
= \lim_{n\to\infty} \dashint_{\Omega_{\rho_n}(y)} \phi^\gamma
= \phi^\gamma(y)
\end{equation*}
for any $\phi\in W^{1,p}_0(\Omega;\mathbb{R}^m)$ with $p>d$,
where we use the definition of the approximating Green's matrix. Due to Theorem $\ref{thm:1.1}$,
there exists the weak solution $u_\varepsilon\in W_0^{1,p}(\Omega;\mathbb{R}^m)$ satisfying
$\mathcal{L}_\varepsilon(u_\varepsilon) = F$ in $\Omega$ and $u_\varepsilon = 0$ on $\partial\Omega$
for any $F\in L^{\frac{p}{2}}(\Omega;\mathbb{R}^m)$ with $p>d$. Thus we obtain
\begin{equation*}
 u^\gamma_\varepsilon(y)
 = \mathrm{B}_\varepsilon^*[\mathcal{G}_\varepsilon^\gamma(\cdot,y),u_\varepsilon]
 = \mathrm{B}_\varepsilon[u_\varepsilon,\mathcal{G}_\varepsilon^\gamma(\cdot,y)]
 = \int_\Omega \mathcal{G}_\varepsilon^\gamma(\cdot,y)F.
\end{equation*}

We now verify the uniqueness. If  $\widetilde{\mathcal{G}}_\varepsilon^\gamma(\cdot,y)$ is another Green's matrix of $\mathcal{L}_\varepsilon$, we then
have $\tilde{u}_\varepsilon^\gamma = \int_\Omega \widetilde{\mathcal{G}}_\varepsilon^\gamma(\cdot,y) F$.
It follows from the uniqueness of the weak solution that
$ \int_\Omega \big[\widetilde{\mathcal{G}}_\varepsilon^\gamma(\cdot,y) - \mathcal{G}_\varepsilon^\gamma(\cdot,y)\big] F = 0$
for any $F\in L^{\frac{p}{2}}(\Omega;\mathbb{R}^m)$, hence
$\widetilde{\mathcal{G}}_\varepsilon^\gamma(\cdot,y) = \mathcal{G}_\varepsilon^\gamma(\cdot,y)$ a.e. in $\Omega$.

Next, let $^*G_{\varrho,\varepsilon}(\cdot,x)$ denote the approximating adjoint of $G_{\rho,\varepsilon}(\cdot,y)$,
which satisfies
\begin{equation}\label{def:2}
 \mathrm{B}_\varepsilon[^*G_{\varrho,\varepsilon}^{\theta}(\cdot,x),u] = \dashint_{\Omega_\varrho(x)} u^\theta,
 \qquad \forall ~u\in H_0^1(\Omega;\mathbb{R}^m).
\end{equation}
By the same argument, we can derive the existence and uniqueness of $^*\mathcal{G}_\varepsilon(\cdot,x)$,
as well as
the estimates similar to $\eqref{pri:4.13}$ and $\eqref{pri:4.14}$. Thus for any $\rho,\varrho>0$, we obtain
\begin{eqnarray*}
 \dashint_{\Omega_\rho(y)}{^*G}^{\gamma\theta}_{\varrho,\varepsilon}(z,x) dz
 = \mathrm{B}_\varepsilon^*[G_{\rho,\varepsilon}^\gamma(\cdot,y),{^*G}^{\theta}_{\varrho,\varepsilon}(\cdot,x)]
 = \mathrm{B}_\varepsilon[{^*G}^{\theta}_{\varrho,\varepsilon}(\cdot,x),G_{\rho,\varepsilon}^\gamma(\cdot,y)]
 =  \dashint_{\Omega_\varrho(x)}{G}^{\theta\gamma}_{\rho,\varepsilon}(z,y) dz .
\end{eqnarray*}
Note that $\mathcal{L}_\varepsilon[{^*G}_{\varrho,\varepsilon}^{\theta}(\cdot,x)] = 0$ in $\Omega\setminus B_\varrho(x)$ and
$\mathcal{L}_\varepsilon^*[{G}_{\rho,\varepsilon}^{\gamma}(\cdot,y)] = 0$ in $\Omega\setminus B_\rho(y)$.
In view of Corollary $\ref{cor:4.2}$, ${^*G}_{\varrho,\varepsilon}^{\theta}(\cdot,x)$
and ${G}_{\rho,\varepsilon}^{\gamma}(\cdot,y)$ are locally H\"{o}lder continuous.
Therefore, we have
${^*\mathcal{G}}_{\varepsilon}^{\gamma\theta}(y,x) = {\mathcal{G}}_{\varepsilon}^{\theta\gamma}(x,y)$ as $\rho\to 0$ and $\varrho\to 0$,
which implies $\mathcal{G}_\varepsilon(x,y) = [{^*\mathcal{G}}_\varepsilon(y,x)]^*$ for any $x,y\in\Omega$ and $x\not=y$.

Let $r=|x-y|$ and $F\in C^\infty_0(\Omega_{\frac{r}{3}}(x))$. Assume $u_\varepsilon$ is the solution of
$\mathcal{L}_\varepsilon(u_\varepsilon) = F$ in $\Omega$ and $u_\varepsilon =0$ on $\partial\Omega$. Then we have
$u_\varepsilon(y) = \int_\Omega \mathcal{G}_\varepsilon(z,y)F(z) dz$. Since $\mathcal{L}_\varepsilon(u_\varepsilon) = 0$ in $\Omega\setminus\Omega_{\frac{r}{3}(x)}$,
it follows from Corollary $\ref{cor:4.2}$ that
\begin{eqnarray*}
 |u_\varepsilon(y)|
 & \leq  &C\Big(\dashint_{\Omega_{\frac{r}{3}}(y)}|u_\varepsilon(z)|^2 dz\Big)^{\frac{1}{2}}
  \leq  Cr^{1-\frac{d}{2}}\Big(\int_{\Omega}|u_\varepsilon(z)|^{2^*} dz\Big)^{\frac{1}{2^*}}\\
 & \leq  &Cr^{1-\frac{d}{2}}\Big(\int_{\Omega}|\nabla u_\varepsilon(z)|^{2} dz\Big)^{\frac{1}{2}}
  \leq  Cr^{1-\frac{d}{2}} \|F\|_{L^{\frac{2d}{d+2}}(\Omega)}
  \leq  C r^{2-\frac{d}{2}}\|F\|_{L^2(\Omega_{\frac{r}{3}}(x))},
\end{eqnarray*}
where we use H\"older's inequality in the second inequality and  Sobolev's inequality in the third inequality,
as well as the estimate $\eqref{pri:4.3}$ with $p=2$ in the fourth inequality.
This implies
\begin{equation}\label{f:4.11}
 \Big(\dashint_{\Omega_{\frac{r}{3}}(x)}|\mathcal{G}_\varepsilon(z,y)|^2 dz\Big)^{\frac{1}{2}} \leq Cr^{2-d}.
\end{equation}
Note that $\mathcal{L}_\varepsilon[\mathcal{G}_\varepsilon(\cdot,y)] = 0$ in $\Omega\setminus B(y,r)$ for any $r>0$.
So in the case of $r\leq 3d_x$,
it follows from $\eqref{pri:4.8}$ and $\eqref{f:4.11}$ that
\begin{eqnarray*}
 |\mathcal{G}_\varepsilon(x,y)| \leq C\Big(\dashint_{\Omega_{\frac{r}{3}}(x)}|\mathcal{G}_\varepsilon(z,y)|^2 dz\Big)^{\frac{1}{2}}
 \leq \frac{C}{|x-y|^{d-2}}.
\end{eqnarray*}
For $r > 3d_x$, in view of $\eqref{pri:4.12}$ and $\eqref{f:4.11}$, for any $\sigma\in(0,1)$, we have
\begin{eqnarray*}
 |\mathcal{G}_\varepsilon(x,y)| = |\mathcal{G}_\varepsilon(x,y) - \mathcal{G}_\varepsilon(\bar{x},y)|
 \leq C(\frac{|x-\bar{x}|}{r})^\sigma \Big(\dashint_{\Omega_{\frac{r}{3}}(x)}|\mathcal{G}_\varepsilon(z,y)|^2 dz\Big)^{\frac{1}{2}}
 \leq C\frac{d_x^\sigma}{|x-y|^{d-2+\sigma}},
\end{eqnarray*}
where $\bar{x}\in \partial\Omega$ such that $d_x = |x-\bar{x}|$. By the same argument, we
can obtain similar results for $^*\mathcal{G}_\varepsilon(\cdot,x)$.
Since $|\mathcal{G}_\varepsilon(x,y)|=|{^*\mathcal{G}}_\varepsilon(y,x)|$ for any $x,y\in \Omega$ and $x\not=y$,
the following results are easily derived by the same arguments:
\begin{eqnarray*}
 |\mathcal{G}_\varepsilon(x,y)| \leq \frac{Cd_y^{\sigma^\prime}}{|x-y|^{d-2+\sigma^\prime}}\min\Big\{1,\frac{d_x^\sigma}{|x-y|^\sigma}\Big\}
\end{eqnarray*}
for any $\sigma,\sigma^\prime\in (0,1)$, where $C$ depends only on $\mu,\omega,\kappa,\lambda,d,m,\sigma$ and $\Omega$. The proof is complete.
\qed
\begin{remark}
\emph{We will see in Section 4 that the estimates $\eqref{pri:4.17}$ actually hold for $\sigma=\sigma^\prime =1$, which are
\begin{equation*}
 |\mathcal{G}_\varepsilon(x,y)|\leq
 \frac{C}{|x-y|^{d-2}}\min\Big\{1,\frac{d_y}{|x-y|},\frac{d_x}{|x-y|},\frac{d_xd_y}{|x-y|^2}\Big\}.
\end{equation*}
Let $r=|x-y|$, and $\bar{y}\in\partial\Omega$ such that $d_y = |y-\bar{y}|$.
In the case of $d_y < r/6$, due to $\mathcal{G}_\varepsilon(x,\cdot) = 0$ on $\partial\Omega$, we have
\begin{equation}\label{f:4.2}
\begin{aligned}
|\mathcal{G}_\varepsilon(x,y)|
&= |\mathcal{G}_\varepsilon(x,y)-\mathcal{G}_\varepsilon(x,\bar{y})|\\
&\leq \|\nabla\mathcal{G}_\varepsilon(x,\cdot)\|_{L^\infty(\Omega_{\frac{r}{6}}(y))}|y-\bar{y}|
\leq \frac{Cd_y}{r}\Big(\dashint_{\Omega_{\frac{r}{3}}(y)} |\mathcal{G}_\varepsilon(x,z)|^2dz\Big)^{\frac{1}{2}}
\leq \frac{Cd_y}{|x-y|^{d-1}},
\end{aligned}
\end{equation}
where we employ the estimate $\eqref{f:5.26}$ in the second inequality,
and $\eqref{f:4.11}$ in the last one. In the case of $d_y\geq r/6$, we can straightforward derive the above estimate from
the estimate $|\mathcal{G}_\varepsilon(x,y)|\leq C|x-y|^{2-d}$. Similarly, we can derive
\begin{equation*}
|\mathcal{G}_\varepsilon(x,y)| = |^*\mathcal{G}_\varepsilon(y,x)|\leq Cd_x|x-y|^{1-d}.
\end{equation*}
Then we plug the above estimate back into the last inequality of $\eqref{f:4.2}$,
and obtain $|\mathcal{G}_\varepsilon(x,y)|\leq Cd_xd_y|x-y|^{-d}$. }
\end{remark}

\begin{remark}
 \emph{ The main idea in the proofs of Theorem $\ref{thm:4.2}$ and Lemma $\ref{lemma:4.5}$ can be found in \cite{MAFHL,SHSK}.
 We comment that the indices $\sigma$ and $\sigma^\prime\in (0,1)$ can be equal,
 which actually come from the H\"older estimate with zero boundary data.
 Equipped with the estimate $\eqref{pri:4.17}$, it is possible to arrive at the sharp H\"older estimate with nonzero boundary data. }
\end{remark}

\begin{flushleft}
\textbf{Proof of Theorem \ref{thm:1.4}}\textbf{.}\quad We first assume that $u_{\varepsilon,1}$ satisfies
$\mathcal{L}_\varepsilon(u_{\varepsilon,1}) = 0$ in $\Omega$ and $u_{\varepsilon,1} = g$ on $\partial\Omega$.
Let $v$ be the extension function of $g$, satisfying
$\Delta v^\alpha = 0$ in $\Omega$ and $v^\alpha = g^\alpha$ on $\partial\Omega$.
For any $x\in\Omega$, set $B = B(x,d_x)$. We have the estimate
\end{flushleft}
\begin{equation}\label{f:4.32}
|\nabla v(x)|
\leq C\Big(\dashint_{\frac{1}{4}B}|\nabla v|^2dy\Big)^{\frac{1}{2}}
\leq \frac{C}{~d_x~}\Big(\dashint_{\frac{1}{2}B}|v(y)-v(x)|^2dy\Big)^{\frac{1}{2}}
\leq Cd_x^{\sigma-1}[v]_{C^{0,\sigma}(\Omega)}
\leq Cd_x^{\sigma-1}\|g\|_{C^{0,\sigma}(\partial\Omega)}
\end{equation}
for any $\sigma\in(0,1)$, where we use the (interior) Lipschitz estimate $\eqref{f:2.15}$ in the first inequality, Cacciopolli's inequality in the second inequality,
and the H\"older estimate: $[v]_{C^{0,\sigma}(\Omega)}\leq C\|g\|_{C^{0,\sigma}(\partial\Omega)}$ in the last inequality.
By normalization we may assume $\|g\|_{C^{0,\sigma}(\partial\Omega)}=1$.
Let $w_\varepsilon = u_{\varepsilon,1} - v$, then $\mathcal{L}_\varepsilon(w_\varepsilon) = - \mathcal{L}_\varepsilon(v)$ in $\Omega$
and $w_\varepsilon = 0$ on
$\partial\Omega$. It follows from $\eqref{pri:4.16}$ that
\begin{eqnarray*}
w_\varepsilon(y) = -\int_\Omega \nabla\mathcal{G}_\varepsilon(x,y)\big[A_\varepsilon\nabla v + V_\varepsilon v\big] dx
- \int_\Omega \mathcal{G}_\varepsilon(x,y)\big[B_\varepsilon\nabla v + (c_\varepsilon + \lambda)v\big] dx ,
\end{eqnarray*}
which implies
\begin{eqnarray*}
|w_\varepsilon(y)| &\leq & C\int_\Omega |\nabla \mathcal{G}_\varepsilon(x,y)|[d_x]^{\sigma-1} dx
+ C\int_\Omega |\mathcal{G}_\varepsilon(x,y)|[d_x]^{\sigma-1} dx
+ C \int_\Omega \big(|\nabla \mathcal{G}_\varepsilon(x,y)| + |\mathcal{G}_\varepsilon(x,y)|\big) dx \nonumber\\
&=: & I_1 + I_2 + I_3,
\end{eqnarray*}
where we use the estimate $\eqref{f:4.32}$.

To estimate $I_1$, set $r = d_y/2$. It follows from  $\eqref{pri:4.10}$ and $\eqref{pri:4.17}$ that
\begin{eqnarray*}
\int_{B(y,r)}|\nabla_x \mathcal{G}_\varepsilon(x,y)|[d_x]^{\sigma-1} dx
&\leq & C\sum_{j=0}^\infty (2^{-j}r)^d\Big(\dashint_{B(y,2^{-j}r)\setminus B(y,2^{-j-1}r)}|\nabla_x \mathcal{G}_\varepsilon(x,y)|^2 dx\Big)^{\frac{1}{2}}r^{\sigma-1} \\
&\leq & C\sum_{j=0}^\infty (2^{-j}r)^{d-1}\Big(\dashint_{B(y,2^{-j+1}r)\setminus B(y,2^{-j-2}r)}|\mathcal{G}_\varepsilon(x,y)|^2 dx\Big)^{\frac{1}{2}}r^{\sigma-1}
\leq  Cr^\sigma,
\end{eqnarray*}
Next, we address ourselves to the integral on $\Omega\setminus B(y,r)$. Let $Q$ be a cube in $\mathbb{R}^d$ with the property that
$3Q\subset\Omega\setminus\{y\}$,
and $l(Q)$, $\text{dist}(Q,\partial\Omega)$ are comparable, where $l(Q)$ denotes the side length of $Q$ (see \cite[pp.167]{S}).
Thus, for fixed $z\in Q$ there exist $c_1,c_2>0$ such that
$c_1|z-y|\leq |x-y|\leq c_2|z-y|$ for any $x\in Q$, and we have
\begin{eqnarray*}
\int_Q |\nabla_x \mathcal{G}_\varepsilon(x,y)|[d_x]^{\sigma-1} dx
&\leq& C[l(Q)]^{\sigma-2}|Q|\Big(\dashint_{2Q}|\mathcal{G}_\varepsilon(x,y)|^2dx\Big)^{\frac{1}{2}} \\
&\leq& C[l(Q)]^{\sigma+\sigma_1-2}|Q|\frac{r^{\sigma_2}}{|z-y|^{d-2+\sigma_1+\sigma_2}}
\leq Cr^{\sigma_2} \int_Q \frac{[d_x]^{\sigma+\sigma_1-2}}{|x-y|^{d-2+\sigma_1+\sigma_2}}dx,
\end{eqnarray*}
where we use the estimate $\eqref{pri:4.10}$ in the first inequality, the estimate $\eqref{pri:4.17}$ in the second one,
and the Chebyshev's inequality in the last one. (Note that $\sigma_1$ and $\sigma_2$ will be given later.)
By decomposing $\Omega\setminus B(y,r)$ as a non-overlapping union of cubes $Q$ (see \cite[pp.167-170]{S}), we then obtain
\begin{eqnarray*}
\int_{\Omega\setminus B(y,r)}|\nabla_x \mathcal{G}_\varepsilon(x,y)|[d_x]^{\sigma-1} dx
\leq Cr^{\sigma_2}\int_{\Sigma_0\cup\Sigma_1\cup\Sigma_2} \frac{[d_x]^{\sigma+\sigma_1-2}}{(r+|x-y|)^{d-2+\sigma_1+\sigma_2}}dx =: I_{11} + I_{12} + I_{13}.
\end{eqnarray*}
Note that we add the additional distance $r$ in the denominator in the second inequality
and therefore the corresponding domain of integral becomes the union of $\Sigma_0 = B(y,r)$, $\Sigma_1 = \cup_{j=0}^\infty \Omega_j$ and $\Sigma_2 = \cup_{j=0}^\infty \cup_{m=0}^\infty\Omega_{j,m}$,
where
\begin{eqnarray*}
\Omega_j &=& \Omega \cap \{2^jr\leq |x-y|\leq 2^{j+1}r\}\cap \{2^jr\leq d_x \leq 2^{j+1}r+2r\}, \\
\Omega_{j,m} & = & \Omega \cap \{2^jr\leq |x-y|\leq 2^{j+1}r\}\cap \{2^{-m-1}(2^jr)\leq d_x \leq 2^{-m}(2^{j}r)\}.
\end{eqnarray*}
Then a routine computation gives rise to $I_{11}\leq Cr^\sigma$,
\begin{eqnarray*}
I_{12} &\leq&  Cr^{\sigma_2}\cdot \sum_{j=0}^\infty \frac{(2^jr)^{\sigma+\sigma_1-2}}{(2^jr)^{d-2+\sigma_1+\sigma_2}}\cdot(2^jr)^d
\leq Cr^\sigma \sum_{j=0}^\infty(2^j)^{\sigma-\sigma_2}\leq Cr^\sigma, \\
I_{13} &\leq& Cr^{\sigma_2}\cdot\sum_{j=0}^\infty \sum_{m=0}^\infty\frac{[2^{-m-1}(2^jr)]^{\sigma+\sigma_1-1}}{(2^jr)^{d-2+\sigma_1+\sigma_2}}\cdot (2^jr)^{d-1}
\leq Cr^\sigma\sum_{m=0}^\infty(2^{-m})^{\sigma+\sigma_1-1}\sum_{j=0}^\infty (2^j)^{\sigma-\sigma_2}\leq Cr^\sigma,
\end{eqnarray*}
provided we choose $\sigma_1,\sigma_2\in(0,1)$ such that $\sigma_1 + \sigma > 1$ and $\sigma_2 < \sigma$. Combining the above estimates,
we obtain $I_1 \leq C[d_y]^\sigma$. In view of $\eqref{pri:4.17}$, we obtain
\begin{equation*}
 I_2 \leq C[d_y]^\sigma\int_\Omega \frac{1}{|x-y|^{d-1}} dx \leq C[d_y]^\sigma,
\end{equation*}
and $I_3\leq C(I_1+I_2)$. Hence, for any $y\in\Omega$ we have
\begin{equation}\label{f:4.12}
 |w_\varepsilon(y)| \leq I_1 + I_2 + I_3 \leq C[d_y]^\sigma.
\end{equation}


Consider three cases: $(1)~|x-y|\leq d_x/4$;
$(2)~|x-y|\leq d_y/4$; $(3)~|x-y|>\max\{d_x/4,d_y/4\}$. In the first case, let $r=d_x$.
In view of $\eqref{f:4.32}$ and $\eqref{f:4.12}$ we first have
\begin{equation}\label{f:4.13}
 \sup_{z\in B(x,r/2)}|\nabla v(z)| \leq Cr^{\sigma-1}
 \qquad \text{and}\qquad
 \sup_{z\in B(x,r/2)}|w_\varepsilon(z)|\leq Cr^{\sigma}.
\end{equation}
It follows from $\eqref{pri:4.9}$, $\eqref{f:4.12}$ and $\eqref{f:4.13}$ that
\begin{equation*}
\begin{aligned}
|w_\varepsilon(x) - w_\varepsilon(y)|
&\leq [w_\varepsilon]_{C^{0,\sigma}(B(x,d_x/4))}|x-y|^\sigma  \\
&\leq C|x-y|^\sigma\bigg\{r^{-\sigma}\Big(\dashint_{2B}|w_\varepsilon|^2 dz\Big)^{\frac{1}{2}}
+ r^{1-\sigma}\Big(\dashint_{2B}(|\nabla v|^p + |v|^p) dz\Big)^{\frac{1}{p}}
+ r^{2-\sigma}\Big(\dashint_{2B}(|\nabla v|^q + |v|^q) dz\Big)^{\frac{1}{q}}\bigg\} \\
& \leq C|x-y|^\sigma,
\end{aligned}
\end{equation*}
where we also use the fact $\|v\|_{L^\infty(\Omega)}\leq C$ (the maximum principle) in the last inequality.
It is clear to see that we can handle the second case in the same manner.

In the case of (3), we derive
\begin{eqnarray*}
|w_\varepsilon(x) - w_\varepsilon(y)|
\leq |w_\varepsilon(x)| + |w_\varepsilon(y)|
\leq C|x-y|^\sigma
\end{eqnarray*}
from $\eqref{f:4.12}$. Thus we have proved $\|w_{\varepsilon}\|_{C^{0,\sigma}(\Omega)}\leq C\|g\|_{C^{0,\sigma}(\partial\Omega)}$.
This together with $\|v\|_{C^{0,\sigma}(\Omega)}\leq C\|g\|_{C^{0,\sigma}(\Omega)}$ gives
\begin{equation*}
\|u_{\varepsilon,1}\|_{C^{0,\sigma}(\Omega)}
\leq \|w_{\varepsilon}\|_{C^{0,\sigma}(\Omega)} + \|v\|_{C^{0,\sigma}(\Omega)}
\leq C\|g\|_{C^{0,\sigma}(\partial\Omega)}.
\end{equation*}

In addition, assume that $u_{\varepsilon,2}$ satisfies
$\mathcal{L}_\varepsilon(u_{\varepsilon,2}) = \text{div}{f} + F$ in $\Omega$ and $u_{\varepsilon,2} = 0$ on $\partial\Omega$.
It follows from Corollary $\ref{cor:4.1}$ that $\|u_{\varepsilon,2}\|_{C^{0,\sigma}(\Omega)}\leq C\{\|f\|_{L^p(\Omega)}+ \|F\|_{L^q(\Omega)}\}$.
Let $u_\varepsilon = u_{\varepsilon,1} + u_{\varepsilon,2}$, we finally obtain $\eqref{pri:1.5}$ and complete the proof.
\qed

\section{Lipschitz estimates \& Nontangential maximal function estimates}
\begin{lemma}\label{lemma:5.4}
Suppose $A\in\Lambda(\mu,\tau,\kappa)$.
Let $p>d$ and $\nu\in(0,\eta]$. Assume that $f = (f_i^\alpha)\in C^{0,\nu}(\Omega;\mathbb{R}^{md})$, $F\in L^p(\Omega;\mathbb{R}^m)$
and $g\in C^{1,\nu}(\partial\Omega;\mathbb{R}^m)$.
Then the unique solution $u_\varepsilon$ to
$L_\varepsilon(u_\varepsilon) = \emph{div}(f) + F $ in $\Omega$ and $u_\varepsilon =g$ on $\partial\Omega$
satisfies the uniform estimate
\begin{equation}\label{pri:5.12}
 \|\nabla u_\varepsilon\|_{L^\infty(\Omega)} \leq C \left\{\|f\|_{C^{0,\nu}(\Omega)} + \|F\|_{L^p(\Omega)}
 + \|g\|_{C^{1,\nu}(\partial\Omega)}\right\},
\end{equation}
where $C$ depends only on $\mu,\tau,\kappa, m,d,\eta,p,\nu$ and $\Omega$.
\end{lemma}

\begin{pf}
  See \cite[Remark 16]{CC}. In fact, $\eqref{pri:5.12}$ is a updated version of \cite[Theorem 2]{MAFHL}.
\qed
\end{pf}

\begin{lemma}
Suppose $A\in \Lambda(\mu,\tau,\kappa)$. Let $\Gamma_\varepsilon(x,y)$ denote the fundamental solution of $L_\varepsilon$,
then we have
\begin{equation}\label{f:5.4}
\max\{|\nabla_x\Gamma_\varepsilon(x,z)|,~|\nabla_z\Gamma_\varepsilon(x,z)|\}\leq C|x-z|^{1-d},
\qquad |\nabla_x\nabla_z\Gamma_\varepsilon(x,z)|\leq C|x-z|^{-d},
\end{equation}
where $C$ depends only on $\mu,\tau,\kappa,m,d$.
\end{lemma}
\begin{pf}
See \cite[pp.6]{SZW6}.
\qed
\end{pf}



\begin{lemma}
Suppose $A\in\Lambda(\mu,\tau,\kappa)$. Let $F\in L^p(\Omega;\mathbb{R}^{m})$ with $p>d$, $f\in C^{0,\sigma}(\Omega;\mathbb{R}^{md})$ with
$\sigma \in (0,1)$,
and $u_\varepsilon\in H_{loc}^1(\Omega;\mathbb{R}^m)$ be the weak solution of $L_\varepsilon(u_\varepsilon) = \emph{div}(f) + F $ in $\Omega$.
Then for any $B\subset2B\subset\Omega$, we have $|\nabla u_\varepsilon|\in L^\infty(B)$ and the uniform estimate
\begin{equation}\label{pri:5.8}
 \|\nabla u_\varepsilon\|_{L^\infty(B)} \leq \frac{C}{r}\Big(\dashint_{2B}|u_\varepsilon|^2\Big)^{\frac{1}{2}}
 + C\left\{\|f\|_{L^\infty(2B)} + r^\sigma[f]_{C^{0,\sigma}(2B)}
 + r\Big(\dashint_{2B}|F|^p\Big)^{\frac{1}{p}} \right\},
\end{equation}
where $C$ depends only on $\mu,\tau,\kappa, p,d,m$ and $\sigma$.
\end{lemma}

\begin{pf}
By rescaling we may assume $r=1$, and $\varphi\in C_0^\infty(2B)$ is a cut-off function
such that $\varphi = 1$ on $5/4B$, $\varphi = 0$ outside $3/2B$, and $|\nabla\varphi|\leq C$. Then we have
\begin{equation*}
 -\text{div}[A_\varepsilon^{\alpha\beta}\nabla(\varphi u_\varepsilon^\beta)]
 = \text{div}(f^\alpha\varphi)  + F^\alpha\varphi - f^\alpha\nabla\varphi
 - A_\varepsilon^{\alpha\beta}\nabla u_\varepsilon^\beta\nabla\varphi
 - \text{div}(A_\varepsilon^{\alpha\beta}\nabla\varphi u_\varepsilon^\beta).
\end{equation*}
It follows from the fundamental solution that for any $x\in B$,
\begin{equation}\label{f:4.33}
u_\varepsilon(x)
 =  -\int_{2B}\nabla_y\Gamma_{\varepsilon}(x,y)f(y)\varphi(y) dy
+ \int_{2B} \Gamma_{\varepsilon}(x,y)\big[F\varphi - f\nabla\varphi - A_\varepsilon \nabla u_\varepsilon\nabla\varphi\big] dy
+ \int_{2B} \nabla_y\Gamma_{\varepsilon}(x,y)A_\varepsilon\nabla\varphi u_\varepsilon dy,
\end{equation}
where we use integration by part in first and third term in right hand side of the above equality.
Then differentiating both sides of $\eqref{f:4.33}$ with respect to $x$ gives
\begin{eqnarray*}
 \nabla u_\varepsilon(x) &=& -\int_{2B}\nabla_x\nabla_y\Gamma_{\varepsilon}(x,y)\big[f(y)\varphi(y)-f(x)\varphi(x)\big] dy
 - f(x)\varphi(x)\int_{\partial(2B)}\nabla_x\Gamma_{\varepsilon}(x,y) n(y) dS(y) \\
 &+& \int_{2B} \nabla_x\Gamma_{\varepsilon}(x,y)\big[F\varphi - f\nabla\varphi - A_\varepsilon\nabla u_\varepsilon\nabla\varphi\big] dy
 + \int_{2B}\nabla_x\nabla_y\Gamma_{\varepsilon}(x,y)A_\varepsilon\nabla\varphi u_\varepsilon  dy,
\end{eqnarray*}
where $dS$ denotes the surface measure of $\partial\Omega$,
and $n$ is the outward unit normal to $\partial(2B)$.
We refer the reader to \cite[pp.55]{DGNT} for the skill used above.

Hence, in view of $\eqref{f:5.4}$, we obtain
\begin{eqnarray*}
|\nabla u_\varepsilon(x)|
&\leq & C\Bigg\{\int_{2B}\frac{|f(y)\varphi(y)-f(x)\varphi(x)|}{|x-y|^{d}} dy
+ \sup_{x\in B}|f(x)\varphi(x)|\int_{|x-y| = 1} \frac{dS(y)}{|x-y|^{d-1}} \\
& + & \int_{2B} \frac{|F(y)|}{|x-y|^{d-1}}dy + \int_{(3/2B)\setminus(5/4B)} \frac{|f(y)|+|\nabla u_\varepsilon(y)|}{|x-y|^{d-1}}dy
+ \int_{(3/2B)\setminus(5/4B)}\frac{|u_\varepsilon(y)|}{|x-y|^d}dy\Bigg\},
\end{eqnarray*}
where we use the observation of $\nabla\varphi = 0$ on $5/4B$ and $\varphi = 0$ outside $3/2B$ in last two terms. This leads to
\begin{equation*}
|\nabla u_\varepsilon(x)| \leq C\bigg\{[f]_{C^{0,\sigma}(2B)}+\|f\|_{L^\infty(2B)}+\|F\|_{L^p(2B)}
+\Big(\int_{3/2B}|\nabla u_\varepsilon|^2 dy\Big)^{\frac{1}{2}}
+\Big(\int_{2B}|u_\varepsilon|^2 dy\Big)^{\frac{1}{2}}\bigg\}
\end{equation*}
for any $x\in B$. Then it follows from the Cacciopolli's inequality that
\begin{equation*}
 \|\nabla u_\varepsilon\|_{L^\infty(B)} \leq C\{\|u_\varepsilon\|_{L^2(2B)} + \|f\|_{C^{0,\sigma}(2B)} + \|F\|_{L^p(2B)}\},
\end{equation*}
where $p>d$, and $C$ depends on $\mu,\tau,\kappa, p, d,m$ and $\sigma$.
By using the rescaling technique, $\eqref{pri:5.8}$ can be easily derived.
\qed
\end{pf}

\begin{thm}\label{thm:5.1}
$(\text{Interior Lipschitz estimates})$\textbf{.}
Suppose that $A\in\Lambda(\mu,\tau,\kappa)$, $V$ satisfies $\eqref{a:2}$, $\eqref{a:4}$, and $B,c$ satisfy $\eqref{a:3}$.
Let $p>d$ and $\sigma\in (0,1)$. Assume that $u_\varepsilon\in H^1_{loc}(\Omega;\mathbb{R}^m)$ is a weak solution to
$\mathcal{L}_\varepsilon(u_\varepsilon) = \emph{div}(f) + F$ in $\Omega$,
where $f\in C^{0,\sigma}(\Omega;\mathbb{R}^{md})$ and $F\in L^p(\Omega;\mathbb{R}^m)$.
Then for any $B\subset2B\subset\Omega$ with $0<r\leq 1$, we have the uniform estimate
\begin{equation}\label{pri:5.5}
 \|\nabla u_\varepsilon\|_{L^\infty(B)}
 \leq \frac{C}{r}\Big(\dashint_{2B} |u_\varepsilon|^2\Big)^{\frac{1}{2}}
 + C\left\{ \|f\|_{L^\infty(2B)}+ r^\sigma [f]_{C^{0,\sigma}(2B)} + r\Big(\dashint_{2B} |F|^p \Big)^{\frac{1}{p}}\right\},
\end{equation}
where $C$ depends only on $\mu, \tau, \kappa, \lambda, p,d,m$ and $\sigma$.
\end{thm}

\begin{pf}
We only need to prove $\eqref{pri:5.5}$ in the case of $\varepsilon < \varepsilon_0$, where $\varepsilon_0$ will be given later.
Since the estimate $\eqref{pri:5.5}$ immediately follows from the classical results when $\varepsilon \geq \varepsilon_0$.
Consider the transformation $T(x,\varepsilon)= [T^{\beta\gamma}(x,\varepsilon)]$ as follows
\begin{equation}\label{transform:1}
 u_\varepsilon^\beta = T^{\beta\gamma}(x,\varepsilon) v_\varepsilon^\gamma = \big[\delta^{\beta\gamma}+\varepsilon\chi_{0}^{\beta\gamma}(x/\varepsilon)\big]v_\varepsilon^\gamma.
\end{equation}
In view of $\eqref{f:2.9}$, it is not hard to see $T(x,\varepsilon)$ is a diagonally dominant matrix
whenever $\varepsilon<\varepsilon_1 = \varepsilon_1(\mu,\tau,\kappa,m,d)$. Hence we have the existence of
$T^{-1}(x,\varepsilon)$,
\begin{equation}\label{f:5.24}
  1/2 \leq \big\|T(\cdot,\varepsilon)\big\|_{L^\infty(\Omega)} \leq 3/2
  \quad \text{and} \quad
  2/3 \leq \big\|T^{-1}(\cdot,\varepsilon)\big\|_{L^\infty(\Omega)} \leq 2 \qquad \text{for}~\varepsilon\in(0,\varepsilon_1),
\end{equation}
where
$$\|T(\cdot,\varepsilon)\|_{L^\infty(\Omega)} = \sup_{x\in\Omega}\|T(x,\varepsilon)\|_\infty
\quad \text{and} \quad
\|T(x,\varepsilon)\|_\infty=\max_{1\leq \alpha,\beta\leq m}|T^{\alpha\beta}(x,\varepsilon)|.$$
Moreover, in view of
\begin{eqnarray*}
 T^{-1}(x,\varepsilon) - T^{-1}(y,\varepsilon) = T^{-1}(x,\varepsilon)\big[T(y,\varepsilon)-T(x,\varepsilon)\big]T^{-1}(y,\varepsilon),
\end{eqnarray*}
we have
\begin{eqnarray*}
\|T^{-1}(x,\varepsilon) - T^{-1}(y,\varepsilon)\|_\infty
\leq \|T^{-1}(x,\varepsilon)\|_\infty\|T^{-1}(y,\varepsilon)\|_\infty\|T(y,\varepsilon)-T(x,\varepsilon)\|_\infty
\leq C|y-x|^\sigma
\end{eqnarray*}
for any $\sigma\in(0,1]$ and $x,y\in B\subset\Omega$, where we use
$[T(\cdot,\varepsilon)]_{C^{0,\sigma}(\Omega)}\leq C\varepsilon^{1-\sigma}[\chi_0]_{C^{0,\sigma}(Y)}\leq C(\mu,\tau,\kappa,m,d,\sigma)$
which follows from $\eqref{f:2.11}$ and $\eqref{f:2.9}$.
Thus we obtain
\begin{equation}\label{f:5.5}
\| T^{-1}(\cdot,\varepsilon)\|_{C^{0,\sigma}(\Omega)} = \|T^{-1}(\cdot,\varepsilon)\|_{L^\infty(\Omega)}
+ \big[ T^{-1}(\cdot,\varepsilon)\big]_{C^{0,\sigma}(\Omega)}  \leq \max\big\{2, C(\mu,\tau,\kappa,\sigma,m,d)\big\}
\end{equation}
for $\varepsilon\in (0,\varepsilon_1)$.

Consider the new system
\begin{equation*}
-\text{div}\big(\tilde{A}_\varepsilon\nabla v_\varepsilon\big) = \text{div}(\tilde{f}) + \tilde{F} \quad \text{in} ~\Omega,
\end{equation*}
where
$\tilde{A}^{\alpha\gamma}_\varepsilon = A_\varepsilon^{\alpha\beta}\big[\delta^{\beta\gamma}+\varepsilon\chi_{0,\varepsilon}^{\beta\gamma}\big]$,
\begin{equation*}
\tilde{f}^\alpha = f^\alpha + \varepsilon V_\varepsilon^{\alpha\beta}\chi_{0,\varepsilon}^{\beta\gamma}v^\gamma_\varepsilon
\quad \text{and} \quad
 \tilde{F}^\alpha = F^\alpha + A_\varepsilon^{\alpha\beta}\nabla\chi_{0,\varepsilon}^{\beta\gamma}\nabla v_\varepsilon^\gamma
 + V_\varepsilon^{\alpha\beta}\nabla v_\varepsilon^\gamma
 - B_\varepsilon^{\alpha\beta}\nabla u_\varepsilon^\beta
 -\big(c_\varepsilon^{\alpha\beta} + \lambda\delta^{\alpha\beta}\big)u_\varepsilon^\beta.
\end{equation*}
Obviously, there exists $\varepsilon_2 = \varepsilon_2(\mu,\tau,\kappa,d,m)$ such that $\tilde{A}\in \Lambda(\frac{\mu}{2},\tau,\kappa+1)$
whenever $\varepsilon \leq \varepsilon_2$.

Let $\varepsilon_0 = \min\{\varepsilon_1,\varepsilon_2,1\}$ and $\varepsilon\in(0,\varepsilon_0]$. For any $r\in(0,1]$,
it follows from $\eqref{pri:5.8}$ that
\begin{equation}\label{f:5.21}
\|\nabla v_\varepsilon\|_{L^\infty(B)} \leq C\Big\{r^{-1}\|v_\varepsilon\|_{L^\infty(2B)}
+ \|\tilde{f}\|_{L^\infty(2B)} + r^{\sigma^\prime} [\tilde{f}]_{C^{0,\sigma^\prime}(2B)}
+ r\Big(\dashint_{2B}|\tilde{F}|^p\Big)^{\frac{1}{p}}\Big\},
\end{equation}
where $\sigma^\prime = \min\{\sigma,\tau\}$.
For convenience, we denote
\begin{equation*}
 \mathcal{R}(nB) = \frac{1}{r}\Big(\dashint_{nB}|u_\varepsilon|^2\Big)^{\frac{1}{2}} + \|f\|_{L^\infty(nB)}
 + r\Big(\dashint_{nB}|F|^p\Big)^{\frac{1}{p}}.
\end{equation*}
Hence, in view of $\eqref{pri:4.9}$, $\eqref{pri:4.8}$, $\eqref{transform:1}$ and $\eqref{f:5.5}$, we obtain
\begin{equation}\label{f:5.20}
\begin{aligned}
 &\|v_\varepsilon\|_{L^\infty(2B)}  \leq \|T^{-1}(\cdot,\varepsilon)\|_{L^\infty(2B)}\|u_\varepsilon\|_{L^\infty(2B)} \leq Cr\mathcal{R}(4B),\\
 &[v_\varepsilon]_{C^{0,\sigma^\prime}(2B)} \leq  [T^{-1}(\cdot,\varepsilon)]_{C^{0,\sigma^\prime}(2B)}\|u_\varepsilon\|_{L^\infty(2B)}
 + \|T^{-1}(\cdot,\varepsilon)\|_{L^\infty(2B)}[u_\varepsilon]_{C^{0,\sigma^\prime}(2B)}
\leq C\{r+r^{1-\sigma^\prime}\}\mathcal{R}(4B), \\
& \Big(\dashint_{2B}|\nabla v_\varepsilon|^p dx\Big)^{\frac{1}{p}}
\leq \|\nabla T^{-1}(\cdot,\varepsilon)\|_{L^\infty(2B)}\|u_\varepsilon\|_{L^\infty(2B)}+
\|T^{-1}(\cdot,\varepsilon)\|_{L^\infty(2B)}\Big(\dashint_{2B}|\nabla u_\varepsilon|^p dx\Big)^{\frac{1}{p}}
\leq C\{r+1\}\mathcal{R}(4B),
\end{aligned}
\end{equation}
where we use Theorem $\ref{thm:4.1}$ to estimate the term of $\big(\dashint_{2B}|\nabla u_\varepsilon|^p dx\big)^{1/p}$.

Note that
\begin{equation}\label{f:4.34}
\begin{aligned}
\big[\tilde{f}\big]_{C^{0,\sigma^\prime}(2B)}
&\leq [f]_{C^{0,\sigma^\prime}(2B)}
+  \varepsilon [V_\varepsilon\chi_{0,\varepsilon}v_\varepsilon]_{C^{0,\sigma^\prime}(2B)}\\
& \leq  [f]_{C^{0,\sigma^\prime}(2B)}
+ \varepsilon\|V\|_{L^\infty(\mathbb{R}^d)}\|\chi_0\|_{L^\infty(Y)}[v_\varepsilon]_{C^{0,\sigma^\prime}(2B)} \\
&+ \varepsilon^{1-\sigma^\prime}\Big\{[V]_{C^{0,\sigma^\prime}(\mathbb{R}^d)}\|\chi_0\|_{L^\infty(Y)}\|v_\varepsilon\|_{L^\infty(2B)}
+ [\chi_0]_{C^{0,\sigma^\prime}(Y)}\|V\|_{L^\infty(\mathbb{R}^d)}\|v_\varepsilon\|_{L^\infty(2B)}\Big\} \\
&\leq  [f]_{C^{0,\sigma^\prime}(2B)}
+ C\Big\{\|v_\varepsilon\|_{L^\infty(2B)} + [v_\varepsilon]_{C^{0,\sigma^\prime}(2B)}\Big\}
\end{aligned}
\end{equation}
where we use the condition $\eqref{a:4}$ and the estimate $\eqref{f:2.9}$ in the last inequality.
Moreover, it follows from $\eqref{f:5.20}$ and $\eqref{f:4.34}$ that
\begin{eqnarray}\label{f:5.22}
 \|\tilde{f}\|_{L^\infty(2B)} + r^{\sigma^\prime} [\tilde{f}]_{C^{0,\sigma^\prime}(2B)}
&\leq & \|f\|_{L^\infty(2B)} + r^{\sigma^\prime} [f]_{C^{0,\sigma^\prime}(2B)}
+ C\|v_\varepsilon\|_{L^\infty(2B)} + Cr^{\sigma^\prime}\|v_\varepsilon\|_{C^{0,\sigma^\prime}(2B)} \nonumber\\
&\leq & \|f\|_{L^\infty(2B)} + r^{\sigma^\prime} [f]_{C^{0,\sigma^\prime}(2B)} + C\{r+r^{1+{\sigma^\prime}}\}\mathcal{R}(4B)
\end{eqnarray}
and
\begin{eqnarray}\label{f:5.23}
r\Big(\dashint_{2B}|\tilde{F}|^p\Big)^{\frac{1}{p}}
&\leq& r\Big(\dashint_{2B}|F|^p\Big)^{\frac{1}{p}} + Cr\Big(\dashint_{2B}|\nabla v_\varepsilon|^p\Big)^{\frac{1}{p}}
+ Cr\Big(\dashint_{2B}|\nabla u_\varepsilon|^p\Big)^{\frac{1}{p}} + Cr\|u_\varepsilon\|_{L^\infty(2B)} \nonumber\\
&\leq& r\Big(\dashint_{2B}|F|^p\Big)^{\frac{1}{p}} + C\{r+r^2\}\mathcal{R}(4B)
\end{eqnarray}
Combining $\eqref{f:5.21},\eqref{f:5.20},\eqref{f:5.22}$ and $\eqref{f:5.23}$, we have
\begin{equation*}
\|\nabla v_\varepsilon\|_{L^\infty(B)} \leq C\bigg\{\|f\|_{L^\infty(2B)} + r^{\sigma^\prime} [f]_{C^{0,\sigma^\prime}(2B)}
+ r\Big(\dashint_{2B}|F|^p\Big)^{\frac{1}{p}}\bigg\} + C\big\{1+r+r^{1+\sigma^\prime}+r^2\big\}\mathcal{R}(4B),
\end{equation*}
where $C$ depends on $\mu,\tau,\kappa,\lambda,\sigma,p,m,d$. This, together with $\eqref{f:2.11}$, $\eqref{f:5.24}$ and $\eqref{f:5.20}$,
gives
\begin{equation*}
\|\nabla u_\varepsilon\|_{L^\infty(B)} \leq \|\nabla T(\cdot,\varepsilon)\|_{L^\infty(B)}\|v_\varepsilon\|_{L^\infty(B)}
+ \|T(\cdot,\varepsilon)\|_{L^\infty(B)}\|\nabla v_\varepsilon\|_{L^\infty(B)}
\leq C\big\{\mathcal{R}(4B) + r^{\sigma}[f]_{C^{0,\sigma}(4B)}\big\}.
\end{equation*}
Note that since $\sigma^\prime \leq \sigma$ we have
$r^{\sigma^\prime}[f]_{C^{0,\sigma^\prime}(4B)}
\leq Cr^\sigma [f]_{C^{0,\sigma}(4B)}$.

For any $B_r(x_0)\subset B_{2r}(x_0)\subset \Omega$, there exist $\{B_{\frac{r}{4}}(x_i)\}_{i=1}^{N}$ and $x_i\in B_r(x_0)$
such that $B_r(x_0)\subset\cup_{i=1}^N B_{\frac{r}{4}}(x_i)$. It is clear to see $B_{r}(x_i)\subset B_{2r}(x_0)$ for any $x_i\in B_r(x_0)$.
Hence we have
\begin{eqnarray*}
 \|\nabla u_\varepsilon\|_{L^\infty(B_r(x_0))} & \leq & \max_{1\leq i\leq N}\big\{\|\nabla u_\varepsilon\|_{L^\infty(B_{\frac{r}{4}}(x_i))}\big\}
 \leq C\big\{\mathcal{R}(B_r(x_i)) + r^\sigma [f]_{C^{0,\sigma}(B_r(x_i))}\big\}\\
 &\leq & C\Big\{r^{-1}\Big(\dashint_{B_{2r}(x_0)}|u_\varepsilon|^2\Big)^{\frac{1}{2}} + \|f\|_{L^\infty(B_{2r}(x_0))}
 + r^\sigma[f]_{C^{0,\sigma}(B_{2r}(x_0))}
 + r\Big(\dashint_{B_{2r}(x_0)}|F|^p\Big)^{\frac{1}{p}} \Big\},
\end{eqnarray*}
and we complete the proof.
\qed
\end{pf}

To prove the global Lipschitz estimates, we study some properties of the Dirichlet correctors $\Phi_{\varepsilon,k}$, $0\leq k\leq d$,
which actually play a similar role as $\chi_{k}$ in the interior Lipschitz estimates.

\begin{lemma}\label{lemma:5.1}
Suppose $A\in\Lambda(\mu,\tau,\kappa)$. Let $g\in C^{0,1}(\partial\Omega;\mathbb{R}^m)$,
and $u_\varepsilon\in H^1(\Omega;\mathbb{R}^m)$ be the solution of
$L_{\varepsilon}(u_\varepsilon) = 0$ in $\Omega$ and $u_\varepsilon = g$ on $\partial \Omega$. Then for any $Q\in \partial\Omega$ and
$\varepsilon\leq r < \emph{diam}(\Omega)$,
 \begin{equation}\label{pri:5.1}
 \left(\dashint_{B(Q,r)\cap\Omega}|\nabla u_\varepsilon|^2 dx\right)^{\frac{1}{2}}
 \leq C\left\{\|\nabla g\|_{L^\infty(\partial\Omega)}+\varepsilon^{-1}\|g\|_{L^\infty(\partial\Omega)}\right\},
 \end{equation}
 where $C$ depends only on $\mu, \tau, \kappa, d, m$ and $\Omega$~.
\end{lemma}

\begin{lemma}\label{lemma:5.2}
Let $A\in\Lambda(\mu,\tau,\kappa)$. Then we have
\begin{equation}\label{pri:5.3}
 \|\Phi_{\varepsilon,k}^\beta-P_k^\beta\|_{L^\infty(\Omega)}\leq  C\varepsilon
 \qquad \text{and}\qquad \|\nabla \Phi_{\varepsilon,k}\|_{L^\infty(\Omega)} \leq C
\end{equation}
for $1\leq k\leq d$, where $C$ depends only on $\mu, \tau, \kappa, d,m,\eta$ and $\Omega$.
\end{lemma}

\begin{remark}
\emph{ Lemmas $\ref{lemma:5.1}$ and $\ref{lemma:5.2}$ were proved in \cite{SZW0}, as well as in \cite{MAFHL}. Here we omit the proof. }
\end{remark}

\begin{lemma}\label{lemma:5.3}
Assume that $A\in\Lambda(\mu,\tau,\kappa)$, and $V$ satisfies $\eqref{a:2}$, $\eqref{a:4}$. Then for any $\sigma\in(0,1]$, we have
\begin{equation}\label{pri:5.2}
\|\Phi_{\varepsilon,0}-I\|_{L^\infty(\Omega)}\leq  C\varepsilon,
\qquad \|\Phi_{\varepsilon,0}-I\|_{C^{0,\sigma}(\Omega)} \leq C\varepsilon^{1-\sigma}
\qquad\text{and}\qquad\|\nabla \Phi_{\varepsilon,0}\|_{L^\infty(\Omega)} \leq C,
\end{equation}
where $C$ depends only on $\mu, \tau, \kappa, d,m,\eta$ and $\Omega$~.
\end{lemma}

\begin{pf}
Let $r_0 = \text{diam}(\Omega)$. If $\varepsilon \geq cr_0$, $\eqref{pri:5.2}$ follows from the classical boundary Lipschitz estimates for elliptic system in divergence form with
the H\"{o}lder continuous coefficients. If $0<\varepsilon < cr_0$, consider
\begin{equation*}
 u_\varepsilon(x) = \Phi_{\varepsilon,0}(x) - I - \varepsilon\chi_0(x/\varepsilon).
\end{equation*}
Then $L_{\varepsilon}(u_\varepsilon) = L_{\varepsilon}(\Phi_{\varepsilon,0}) - L_{\varepsilon}[\varepsilon\chi_0(x/\varepsilon)] = 0$ in $\Omega$, and
$u_\varepsilon = -\varepsilon\chi_0(x/\varepsilon)$ on $\partial \Omega$.
Hence, it follows from the Agmon-Miranda maximum principle (see \cite[Theorem 3]{MAFHL}
or \cite[Remark 3.4.4]{SZW0}) that
\begin{equation*}
 \sup_{x\in\Omega}|u_\varepsilon(x)| \leq C\sup_{x\in\partial\Omega}|u_\varepsilon(x)|
 \leq C\varepsilon\|\chi_0\|_{L^\infty(\mathbb{R}^d)}\leq  C\varepsilon,
\end{equation*}
which implies $\|\Phi_{\varepsilon,0}-I\|_{L^\infty(\Omega)} \leq  C\varepsilon$.
Additionally, for any $\sigma\in(0,1)$,
in view of Theorem $\ref{thm:1.4}$, we have $\|\Phi_{\varepsilon,0} - I\|_{C^{0,\sigma}(\Omega)}\leq C\varepsilon^{1-\sigma}$.
Note that $L_\varepsilon$ is the special
case of $\mathcal{L}_\varepsilon$,
and $C$ depends only on $\mu,\tau,\kappa, d,m,\sigma$ and $\Omega$.

Moreover, it follows from Lemma $\ref{lemma:5.1}$ that
\begin{equation*}
\Big(\dashint_{D(Q,r)} |\nabla u_\varepsilon|^2 dx \Big)^{1/2} \leq C,
\end{equation*}
which implies
\begin{equation}\label{f:4.1}
 \Big(\dashint_{D(Q,r)} |\nabla \Phi_{\varepsilon,0}|^2 dx\Big)^{1/2} \leq C
\end{equation}
for any $Q\in\partial\Omega$, and $\varepsilon \leq r < r_0$.
By the interior Lipschtiz estimate, we have
\begin{equation*}
 \sup_{\{d_x\geq\varepsilon\}\cap\Omega}|\nabla u_\varepsilon(x)|
 \leq \frac{C}{\varepsilon}\Big(\dashint_{\Omega}|u_\varepsilon|^2 dy\Big)^{1/2}
 \leq C,
\end{equation*}
which gives
\begin{equation}\label{f:5.16}
\sup_{\{d_x\geq\varepsilon\}\cap\Omega}|\nabla \Phi_{\varepsilon,0}(x)| \leq C.
\end{equation}

In the case of $\{d_x<\varepsilon\}\cap\Omega$, we apply the blow-up argument. Let
\begin{equation*}
 v(x) = \frac{1}{\varepsilon}\Phi_{\varepsilon,0}(\varepsilon x) - \frac{1}{\varepsilon}I,
\end{equation*}
then we have
\begin{equation*}
\left\{\begin{aligned}
 L_1(v) & = \text{div}(V) &\qquad &\text{in}~~~\Omega_{\varepsilon} , \\
     v  & = 0         &\qquad &\text{on}~~\partial\Omega_{\varepsilon},
\end{aligned}\right.
\end{equation*}
where $\Omega_\varepsilon = \{x\in \mathbb{R}^d :\varepsilon x \in \Omega\}$.
Note that although the character of the boundary varies, $\Omega_{\varepsilon}$ is still a bounded $C^{1,\eta}$ domain.
The boundary functions of $\Omega_\varepsilon$ are denoted by $\psi_{i,\varepsilon}(x)=\psi_i(\varepsilon x)$, $i=1,\cdots,n_0$
(recall Remark $\ref{rm:2.6}$), and
we fortunately have
$\|\psi_{i,\varepsilon}\|_{C^{1,\eta}(\mathbb{R}^{d-1})} \leq \varepsilon^{1+\eta}\|\psi_i\|_{C^{1,\eta}(\mathbb{R}^{d-1})}
\leq \varepsilon^{1+\eta}M_0$.
Hence, it follows from the (boundary) Lipschitz estimate $\eqref{f:2.15}$ that
\begin{equation*}
\sup_{B(0,1)\cap\Omega_\varepsilon}|\nabla v| \leq C\left\{\Big(\dashint_{B(0,2)\cap\Omega_\varepsilon}|\nabla v|^2dy\Big)^{\frac{1}{2}}
+ \|V\|_{C^{0,\tau}(\mathbb{R}^d)}\right\}.
\end{equation*}
This implies
\begin{eqnarray}\label{f:5.17}
 \sup_{x\in D(0,\varepsilon)}|\nabla \Phi_{\varepsilon,0}(x)|
\leq C\left\{\Big(\dashint_{D(0,2\varepsilon)}|\nabla \Phi_{\varepsilon,0}|^2dy\Big)^{\frac{1}{2}} + \|V\|_{C^{0,\tau}(\mathbb{R}^d)} \right\}
\leq C.
\end{eqnarray}
Note that we choose $r=2\varepsilon$ in $\eqref{f:4.1}$ to give the last inequality.
Combining $\eqref{f:5.16}$ and $\eqref{f:5.17}$,
we have $\|\nabla\Phi_{\varepsilon,0}\|_{L^\infty(\Omega)}\leq C$,
and this also implies the second estimate of $\eqref{pri:5.2}$ for $\sigma =1$.
We thus complete the proof.
\qed
\end{pf}

\begin{lemma}\label{cor:5.1}
Assume the same conditions as in Lemma $\ref{lemma:5.3}$. Then we have
\begin{equation}\label{pri:5.9}
\big\|\nabla \Phi_{\varepsilon,0}\big\|_{C^{0,\tau}(\Omega)} \leq C\max\{\varepsilon^{-\tau},1\}.
\end{equation}
Furthermore, $\Phi_{\varepsilon,0}^{-1}$ exists and satisfies
the following estimates:
\begin{equation}\label{pri:5.11}
  2/3 \leq \big\|\Phi_{\varepsilon,0}^{-1}\big\|_{L^\infty(\Omega)} \leq 2, \qquad
  \big\|\nabla (\Phi_{\varepsilon,0}^{-1})\big\|_{L^\infty(\Omega)} \leq C,
\end{equation}
whenever $\varepsilon\leq \varepsilon_0$,
where $\varepsilon_0 = \varepsilon_0(\mu,\tau,\kappa,d,m,\Omega)$ is sufficiently small,
and $C$ depends only on $\mu, \tau, \kappa, d,m,\eta$ and $\Omega$.
\end{lemma}

\begin{pf}
Let $\widetilde{\Phi}_{\varepsilon,0} = \Phi_{\varepsilon,0} - I$,
then $L_\varepsilon(\widetilde{\Phi}_{\varepsilon,0}) = \text{div}(V_\varepsilon)$ in $\Omega$, and $\widetilde{\Phi}_{\varepsilon,0} = 0$ on
$\partial\Omega$.
We first prove $\eqref{pri:5.9}$ in the case of $\varepsilon<1$. Set $U(\varepsilon) = \Omega\cap B(P,\varepsilon)$
for any $P\in\overline{\Omega}$.
By translation we may assume $P=0$. In view of the Schauder estimate $\eqref{f:2.16}$ and Lemma $\ref{lemma:5.3}$, we have
\begin{eqnarray*}
 \big[\nabla\widetilde{\Phi}_{\varepsilon,0}\big]_{C^{0,\tau}(U(\varepsilon))}
 &\leq& C\varepsilon^{-\tau}\left\{\Big(\dashint_{U(2\varepsilon)}|\nabla\widetilde{\Phi}_{\varepsilon,0}
 |^2\Big)^{\frac{1}{2}} + \|V_\varepsilon\|_{L^\infty(U(2\varepsilon))}
 + \varepsilon^{\tau}[V_\varepsilon]_{C^{0,\tau}(U(2\varepsilon))}\right\} \\
 &\leq&
 C\varepsilon^{-\tau}\Big\{\|\nabla\Phi_{\varepsilon,0}\|_{L^\infty(\Omega)} + \|V\|_{C^{0,\tau}(\mathbb{R}^d)}\Big\}
 \leq C\varepsilon^{-\tau},
\end{eqnarray*}
where $C$ depends only on $\mu,\tau,\kappa,m,d,\eta$ and $M_0$. Note that
$[\nabla\Phi_{\varepsilon,0}]_{C^{0,\tau}(U(\varepsilon))}
= [\nabla\widetilde{\Phi}_{\varepsilon,0}]_{C^{0,\tau}(U(\varepsilon))}$.
Thus by a covering argument (see \cite[pp.98]{DGNT}), we obtain $\|\nabla \Phi_{\varepsilon,0}\|_{C^{0,\tau}(\Omega)} \leq C\varepsilon^{-\tau}$.
The case of $\varepsilon >1$ is trivial, since we can derive $\eqref{pri:5.9}$ by using
the Schauder estimates $\eqref{pri:2.5}$ directly.

Next we prove $\eqref{pri:5.11}$.  It follows from $\eqref{pri:5.2}$ that
$\|\widetilde{\Phi}_{0,\varepsilon}\|_{L^\infty(\Omega)}\leq C\varepsilon$.
Since $\Phi_{\varepsilon, 0} = I + \widetilde{\Phi}_{0,\varepsilon}$,
we know that there exists $\Phi_{\varepsilon,0}^{-1}\in L^\infty(\Omega)$ such that
\begin{equation*}
 1/2 \leq \|\Phi_{\varepsilon,0}\|_{L^\infty(\Omega)} \leq 3/2
 \quad\text{and} \quad 2/3 \leq \|\Phi_{\varepsilon,0}^{-1}\|_{L^\infty(\Omega)} \leq 2
\end{equation*}
whenever $\varepsilon \leq \varepsilon_0(\mu,\tau,\kappa,d,m,\Omega)$, and $\varepsilon_0$ is sufficiently small.

Due to
$$\Phi_{\varepsilon,0}^{-1}(x) - \Phi_{\varepsilon,0}^{-1}(y)
= \Phi_{\varepsilon,0}^{-1}(x)\big[\Phi_{\varepsilon,0}(y) - \Phi_{\varepsilon,0}(x)\big]\Phi_{\varepsilon,0}^{-1}(y),$$
we have
\begin{eqnarray*}
\big| \Phi_{\varepsilon,0}^{-1}(x) - \Phi_{\varepsilon,0}^{-1}(y) \big|
 \leq \|\Phi_{\varepsilon,0}^{-1}\|_{L^\infty(\Omega)}\|\Phi_{\varepsilon,0}^{-1}\|_{L^\infty(\Omega)}
 \|\nabla \Phi_{\varepsilon,0}\|_{L^\infty(\Omega)}|x-y|
 \leq C|x-y|
\end{eqnarray*}
for $x,y\in\Omega$, and this implies $\|\nabla \Phi_{\varepsilon,0}^{-1}\|_{L^\infty(\Omega)} \leq C$.
The proof is complete.
\qed
\end{pf}

\begin{lemma}\label{lemma:5.7}
$(\text{A nonuniform estimate})$\textbf{.}
Suppose that $A\in\Lambda(\mu,\tau,\kappa)$, $V$ satisfies $\eqref{a:2}$, $\eqref{a:4}$, and $B,c$
satisfy $\eqref{a:3}$.
Let $p>d$ and $\sigma\in(0,\tau]$. Assume $f\in C^{0,\sigma}(\Omega;\mathbb{R}^{md})$ and $F\in L^p(\Omega;\mathbb{R}^m)$,
then the weak solution to $\mathcal{L}_\varepsilon(u_\varepsilon) = \emph{div}(f) + F$ in $\Omega$ and $u_\varepsilon = 0$ on $\partial\Omega$
satisfies the estimate
\begin{equation} \label{pri:5.10}
 \|\nabla u_\varepsilon\|_{C^{0,\frac{\sigma}{4}}(\Omega)}
 \leq C\max\big\{\varepsilon^{\frac{\sigma}{2}-1},1\big\}\big\{ \|f\|_{C^{0,\sigma}(\Omega)} + \|F\|_{L^p(\Omega)}\big\},
\end{equation}
where $C$ depends only on $\mu,\tau,\kappa,\lambda,p,d,m,\eta$ and $\Omega$.
\end{lemma}

\begin{pf}
If $\varepsilon \geq 1$, $\eqref{pri:5.10}$ follows directly from the Schauder estimate $\eqref{pri:2.7}$
and the Lipschitz estimate $\eqref{pri:2.6}$.

In the case of $0<\varepsilon<1$,
the main idea is based upon the following interpolation inequality
\begin{eqnarray}\label{f:5.12}
\big[\nabla u_\varepsilon\big]_{C^{0,\frac{\sigma}{4}}(\Omega)}
\leq 2\|\nabla u_\varepsilon\|_{L^\infty(\Omega)}^{\frac{3}{4}}\big[\nabla u_\varepsilon\big]_{C^{0,\sigma}(\Omega)}^{\frac{1}{4}}.
\end{eqnarray}

Set $U(\varepsilon) = \Omega\cap B(P,\varepsilon)$ for any $P\in\overline{\Omega}$, and by translation we may assume $P=0$.
We first study $\|\nabla u_\varepsilon\|_{L^\infty(\Omega)}$ through the uniform H\"older estimates. To do so,
let $v_\varepsilon = u_\varepsilon - [I+\varepsilon\chi_0(x/\varepsilon)]u_\varepsilon(0)$. Hence we have
\begin{equation*}
 \mathcal{L}_\varepsilon(v_\varepsilon) = \text{div}(f) + F
 +\big[ \text{div}(\varepsilon V_\varepsilon\chi_{0,\varepsilon})
 -B_\varepsilon\nabla_y\chi_{0,\varepsilon} - (c_\varepsilon+\lambda I)(I+\varepsilon\chi_{0,\varepsilon})\big]u_\varepsilon(0)
 \qquad \text{in}  ~~\Omega,
\end{equation*}
where $y=x/\varepsilon$.
If $0\in\partial\Omega$, we have $v_\varepsilon = 0$ on $\partial\Omega$.
From the Lipschitz estimate $\eqref{f:2.17}$ on $\varepsilon$ scale, we obtain
\begin{eqnarray*}
\|\nabla v_\varepsilon\|_{L^\infty(U(\varepsilon))}
&\leq& C\left\{\frac{1}{\varepsilon}\Big(\dashint_{U(2\varepsilon)}|v_\varepsilon|^2 \Big)^{\frac{1}{2}}
    + \|f\|_{L^\infty(U(2\varepsilon))} + \varepsilon^\sigma[f]_{C^{0,\sigma}(U(2\varepsilon))}
    +\varepsilon\Big(\dashint_{U(2\varepsilon)}|F|^p\Big)^{\frac{1}{p}} + |u_\varepsilon(0)| \right\}     \\
&\leq  & \frac{C}{\varepsilon}\Big(\dashint_{U(2\varepsilon)}|u_\varepsilon - u_\varepsilon(0)|^2\Big)^{\frac{1}{2}}
+ C\|u_\varepsilon\|_{L^\infty(\Omega)}
+ C\big\|f\big\|_{C^{0,\sigma}(\Omega)} + C\varepsilon^\sigma\big\|F\big\|_{L^{p}(U(2))}  \\
& \leq & \frac{C}{~\varepsilon^{1-\sigma}}\big\|u_\varepsilon\big\|_{C^{0,\sigma}(\Omega)}
+ C\Big\{\big\|f\big\|_{C^{0,\sigma}(\Omega)} + \big\|F\big\|_{L^p(\Omega)}\Big\},
\end{eqnarray*}
where we use $\eqref{pri:4.4}$ in the third inequality. This implies
\begin{equation}\label{f:5.14}
\|\nabla u_\varepsilon\|_{L^\infty(\Omega)}
\leq \|\nabla v_\varepsilon\|_{L^\infty(\Omega)} + \|\nabla \chi_0\|_{L^\infty(Y)}\|u_\varepsilon\|_{L^\infty(\Omega)}
\leq C\varepsilon^{\sigma-1}\Big\{\|f\|_{C^{0,\sigma}(\Omega)} + \|F\|_{L^p(\Omega)} \Big\}.
\end{equation}

Next, it directly follows  from the Schauder estimate $\eqref{f:2.18}$ that
\begin{eqnarray*}
\big[\nabla u_\varepsilon\big]_{C^{0,\sigma}(U(\varepsilon))}
&\leq & C\left\{\frac{1}{\varepsilon^{1+\sigma}}\Big(\dashint_{U(2\varepsilon)}|u_\varepsilon|^2\Big)^{\frac{1}{2}}
+ [f]_{C^{0,\sigma}(U(2\varepsilon))}
+\varepsilon^{-\sigma} \|f\|_{L^\infty(U(2\varepsilon))}
+ \varepsilon^{1-\sigma}\Big(\dashint_{U(2\varepsilon)}|F|^p\Big)^{\frac{1}{p}} \right\} \\
&\leq & \frac{C}{\varepsilon^{1+\sigma}}\Big\{\|u_\varepsilon\|_{L^\infty(\Omega)}
+  \|f\|_{C^{0,\sigma}(\Omega)}\Big\}
+  C\Big(\int_{\Omega}|F|^p\Big)^{\frac{1}{p}}  \\
&\leq & \frac{C}{\varepsilon^{1+\sigma}}
\Big\{ \|f\|_{C^{0,\sigma}(\Omega)} + \|F\|_{L^p(\Omega)} \Big\},
\end{eqnarray*}
where $C$ depends on $\mu,\tau,\kappa,\lambda, p,\sigma,m,d,\eta$ and $\Omega$. By a covering argument (see \cite[pp.98]{DGNT}), we have
\begin{equation}\label{f:5.15}
 [\nabla u_\varepsilon]_{C^{0,\sigma}(\Omega)} \leq C\varepsilon^{-1-\sigma}\Big\{ \|f\|_{C^{0,\sigma}(\Omega)} + \|F\|_{L^p(\Omega)} \Big\}.
\end{equation}

Finally, the estimate $\eqref{pri:5.10}$ follows from
$\eqref{f:5.12}$, $\eqref{f:5.14}$ and $\eqref{f:5.15}$.
We complete the proof.
\qed
\end{pf}

\begin{flushleft}
\textbf{Proof of Theorem \ref{thm:1.2}}\textbf{.}\quad
In the case of $g = 0$,
we only need to consider the following transformation
\begin{equation}\label{transform:2}
u_\varepsilon^\beta(x) = \Phi^{\beta\gamma}_{\varepsilon,0}(x)v^\gamma_\varepsilon(x)
\end{equation}
for $\varepsilon <\varepsilon_*$, where $\varepsilon_* = \min\{\varepsilon_0,\varepsilon_1,\varepsilon_2\}$,
and $\varepsilon_0$ is given in Lemma $\ref{cor:5.1}$ and $\varepsilon_1,\varepsilon_2$ can be chosen later.
Since it is clear to see that the estimate $\eqref{pri:1.2}$ immediately follows from the Lipschitz estimate $\eqref{pri:2.6}$ for
$\varepsilon \geq\varepsilon_*$.
Then the Dirichlet problem $\eqref{pde:1.1}$ can be transformed into
\begin{equation}\label{pde:1.5}
\left\{\begin{aligned}
 L_\varepsilon(v_\varepsilon)
&= \text{div}(\tilde{f}) + \tilde{F} &\quad& \text{in}~~\Omega, \\
 v_\varepsilon & = 0 &\quad& \text{on} ~~\partial\Omega,
\end{aligned}
\right.
\end{equation}
where
\begin{eqnarray*}
\tilde{f}^\alpha &=& f^\alpha + A_\varepsilon^{\alpha\beta}(\Phi_{\varepsilon,0}^{\beta\gamma} - \delta^{\beta\gamma})\nabla v_\varepsilon^\gamma
+ V_\varepsilon^{\alpha\beta}(\Phi_{\varepsilon,0}^{\beta\gamma}-\delta^{\beta\gamma})v_\varepsilon^\gamma,  \\
\tilde{F}^\alpha &=& F^\alpha + A_\varepsilon^{\alpha\beta}\nabla\Phi_{\varepsilon,0}^{\beta\gamma}\nabla v_\varepsilon^\gamma
+ V_\varepsilon^{\alpha\gamma}\nabla v_\varepsilon^\gamma
 - B_\varepsilon^{\alpha\beta}\nabla u_\varepsilon^\beta - (c_\varepsilon^{\alpha\beta}+\lambda\delta^{\alpha\beta})u_\varepsilon^\beta .
\end{eqnarray*}
\end{flushleft}
It follows from Theorem $\ref{thm:1.1}$ and Corollary $\ref{cor:4.1}$ that
\begin{equation*}
\max{\big\{ \|\nabla u_\varepsilon\|_{L^p(\Omega)},~\|u_\varepsilon\|_{C^{0,\sigma^\prime}(\Omega)}\big\}}
\leq C\left\{\|f\|_{L^\infty(\Omega)} + \|F\|_{L^p(\Omega)}\right\}.
\end{equation*}
This, together with Lemma $\ref{cor:5.1}$, gives
\begin{equation}\label{f:5.1}
\max{\big\{\|\nabla v_\varepsilon\|_{L^p(\Omega)}, \|v_\varepsilon\|_{C^{0,\sigma^\prime}(\Omega)}\big\}}
\leq C\left\{\|f\|_{L^\infty(\Omega)} + \|F\|_{L^p(\Omega)}\right\},
\end{equation}
where $\sigma^\prime = 1-d/p$, and $C$ depends on $\mu,\tau,\kappa,\lambda,\sigma,p,d,m$ and $\Omega$.
Here we use $v_\varepsilon = \Phi_{\varepsilon,0}^{-1}u_\varepsilon$, and
$\nabla v_\varepsilon = \nabla(\Phi_{\varepsilon,0}^{-1})u_\varepsilon + \Phi_{\varepsilon,0}^{-1}\nabla u_\varepsilon$.
However we need to rewrite
$\nabla v_\varepsilon = -\nabla\Phi_{\varepsilon,0}v_\varepsilon +(I-\Phi_{\varepsilon,0})\nabla v_\varepsilon + \nabla u_\varepsilon$
to handle the H\"older norm of $\nabla v_\varepsilon$.
Set $\nu = \min\{\tau,\sigma,\sigma^\prime\}/4$ and $\nu^\prime = \max\{\tau,1-2\nu\}$. Note that $0<\nu^\prime<1$, and we obtain
\begin{eqnarray*}
[\nabla v_\varepsilon]_{C^{0,\nu}(\Omega)}
&\leq &[\nabla \Phi_{\varepsilon,0}]_{C^{0,\nu}(\Omega)} \| v_\varepsilon\|_{L^\infty(\Omega)}
+\|\nabla\Phi_{\varepsilon,0}\|_{L^\infty(\Omega)}[v_\varepsilon]_{C^{0,\nu}(\Omega)}
  +  [I-\Phi_{\varepsilon,0}]_{C^{0,\nu}(\Omega)} \|\nabla v_\varepsilon\|_{L^{\infty}(\Omega)} \\
&+& \|I-\Phi_{\varepsilon,0}\|_{L^\infty(\Omega)} [\nabla v_\varepsilon]_{C^{0,\nu}(\Omega)} + [\nabla u_\varepsilon]_{C^{0,\nu}(\Omega)} \\
&\leq & C(\varepsilon^{-\tau}+\varepsilon^{2\nu-1})\big\{\|f\|_{C^{0,\sigma}(\Omega)} + \|F\|_{L^p(\Omega)}\big\}
+ C\varepsilon^{1-\nu}\|\nabla v_\varepsilon\|_{L^\infty(\Omega)} + C\varepsilon[\nabla v_\varepsilon]_{C^{0,\nu}(\Omega)}
\end{eqnarray*}
where we apply $\eqref{pri:5.2}$, $\eqref{pri:5.9}$, $\eqref{pri:5.10}$ and $\eqref{f:5.1}$ to the second inequality. This implies
\begin{eqnarray}\label{f:5.25}
[\nabla v_\varepsilon]_{C^{0,\nu}(\Omega)} \leq C\varepsilon^{-\nu^\prime}\big\{\|f\|_{C^{0,\sigma}(\Omega)} + \|F\|_{L^p(\Omega)}\big\}
+ C\|\nabla v_\varepsilon\|_{L^\infty(\Omega)}
\end{eqnarray}
whenever $\varepsilon < \varepsilon_1$, where
$\varepsilon_1 = \min\{1/(2C),1\}$.
Hence, we have
\begin{eqnarray*}
\|\tilde{f}\|_{C^{0,\nu}(\Omega)} 
&\leq & \|f\|_{C^{0,\sigma}(\Omega)} + C\varepsilon\|\nabla v_\varepsilon\|_{L^\infty(\Omega)}
+ \big[A_\varepsilon(\Phi_{\varepsilon,0}-I)\nabla v_\varepsilon\big]_{C^{0,\nu}(\Omega)}
+ \|V_\varepsilon(\Phi_{\varepsilon,0}-I)v_\varepsilon\|_{C^{0,\nu}(\Omega)} \\
 &\leq & C\varepsilon^{1-\nu^\prime}\|\nabla v_\varepsilon\|_{L^\infty(\Omega)} + C\big\{\|f\|_{C^{0,\sigma}(\Omega)}+\|F\|_{L^p(\Omega)}\big\},
\end{eqnarray*}
where we use $\eqref{a:4}$, $\eqref{pri:5.2}$, $\eqref{pri:5.9}$, $\eqref{f:5.1}$ and $\eqref{f:5.25}$ in the second inequality. In view of $\eqref{f:5.1}$, we also have
\begin{eqnarray*}
\|\tilde{F}\|_{L^p(\Omega)}
&\leq & \|F\|_{L^p(\Omega)}
+ C\big\{\|\nabla v_\varepsilon\|_{L^p(\Omega)} + \|\nabla u_\varepsilon\|_{L^p(\Omega)} + \|u_\varepsilon\|_{L^p(\Omega)} \big\}
\leq C\big\{\|f\|_{L^\infty(\Omega)} + \|F\|_{L^p(\Omega)}\big\}.
\end{eqnarray*}
We now apply Lemma $\ref{lemma:5.4}$ to $\eqref{pde:1.5}$ and obtain
\begin{eqnarray*}
\|\nabla v_\varepsilon\|_{L^{\infty}(\Omega)}
&\leq&  C \big\{\|\tilde{f}\|_{C^{0,\nu}(\Omega)} + \|\tilde{F}\|_{L^p(\Omega)}\big\}\\
&\leq& C\varepsilon^{1-\nu^\prime}\|\nabla v_\varepsilon\|_{L^\infty(\Omega)} + C\big\{ \|f\|_{C^{0,\sigma}(\Omega)} + \|F\|_{L^p(\Omega)}\big\},
\end{eqnarray*}
which gives $\|\nabla v_\varepsilon\|_{L^{\infty}(\Omega)} \leq C\big\{ \|f\|_{C^{0,\sigma}(\Omega)} + \|F\|_{L^p(\Omega)}\big\}$,
whenever $\varepsilon < \varepsilon_2 =\min\{1/(2C)^{\frac{1}{1-\nu^\prime}},1\}$.
So we have
\begin{eqnarray}\label{f:5.10}
\|\nabla u_\varepsilon\|_{L^{\infty}(\Omega)}
\leq \|\nabla (\Phi_{\varepsilon,0}v_\varepsilon)\|_{L^{\infty}(\Omega)}
\leq  C \{\|f\|_{C^{0,\sigma}} + \|F\|_{L^{p}(\Omega)}\},
\end{eqnarray}
where $C$ depends only on $\mu,\tau,\kappa,\lambda,p,\sigma,d,m,M_0,\eta$ and $|\Omega|$.

In the case of $g\not= 0$, consider the homogeneous system $\mathcal{L}_\varepsilon(u_\varepsilon) = 0$ in $\Omega$
and $u_\varepsilon = g$ on $\Omega$,
where $g\in C^{1,\sigma}(\partial\Omega;\mathbb{R}^{m})$ with $\sigma\in(0,\eta]$. Let $h_\varepsilon$ be the extension function of $g$, satisfying
\begin{equation*}
  -\text{div}(A_\varepsilon\nabla h_\varepsilon) = 0 \quad \text{in}~\Omega \qquad\text{and}\qquad h_\varepsilon = g \quad \text{on}~\partial\Omega.
\end{equation*}
It follows from Lemma $\ref{lemma:5.4}$ that $\|\nabla h_\varepsilon\|_{L^\infty(\Omega)}\leq C\|g\|_{C^{1,\sigma}(\partial\Omega)}$,
where $C$ depends only on $\mu,\tau,\kappa,\sigma,d,m,\eta$ and $\Omega$.
Let $\varrho = \min\{\tau,\sigma\}$.
By the argument applied to Lemma $\ref{lemma:5.7}$,
we obtain $[\nabla h_\varepsilon]_{C^{0,\varrho}(\Omega)}
\leq C\varepsilon^{-1-\varrho}\|g\|_{C^{1,\sigma}(\partial\Omega)}$.
Indeed, due to $\eqref{f:2.16}$ we have
\begin{equation*}
 \big[\nabla h_\varepsilon\big]_{C^{0,\varrho}(B(P,\varepsilon))}
 \leq C\varepsilon^{-\varrho}\Big(\dashint_{B(P,2\varepsilon)}|\nabla h_\varepsilon|^2\Big)^{\frac{1}{2}}
 \leq C\varepsilon^{-\varrho}\|\nabla h_\varepsilon\|_{L^\infty(\Omega)}
 \leq C\varepsilon^{-\varrho}\|g\|_{C^{1,\sigma}(\partial\Omega)}
\end{equation*}
for any $B(P,2\varepsilon)\subset\Omega$, while for the boundary estimates,
it follows from the (boundary) Schauder estimate $\eqref{f:2.19}$ that
\begin{eqnarray*}
 \big[\nabla h_\varepsilon\big]_{C^{0,\varrho}(D(\varepsilon))}
 &\leq & \frac{C}{\varepsilon^\varrho}\Big\{\|\nabla h_\varepsilon\|_{L^\infty(D(2\varepsilon))}
 + \|\nabla g\|_{L^\infty(\Delta(2\varepsilon))}
 + \varepsilon^\varrho[\nabla g]_{C^{0,\varrho}(\Delta(2\varepsilon))}
 + \varepsilon^{-1}\|g\|_{L^\infty(\Delta(2\varepsilon))}              \Big\} \\
 &\leq & \frac{C}{\varepsilon^{1+\varrho}}\Big\{\|\nabla h_\varepsilon\|_{L^\infty(\Omega)} + \|g\|_{C^{1,\varrho}(\partial\Omega)}\Big\}
  \leq C\varepsilon^{-1-\varrho}\|g\|_{C^{1,\sigma}(\partial\Omega)}
\end{eqnarray*}
for any $P\in\partial\Omega$, where $C$ depends on $\mu,\tau,\kappa,\varrho,d,m,M_0,\eta$, and $|\Omega|$.
Thus we have
\begin{equation}\label{f:5.6}
[\nabla h_\varepsilon]_{C^{0,\frac{\varrho}{2}}(\Omega)}
\leq 2\|\nabla h_\varepsilon\|_{L^\infty(\Omega)}^{\frac{1}{2}}[\nabla h_\varepsilon]_{C^{0,\varrho}(\Omega)}^{\frac{1}{2}}
\leq C\varepsilon^{-\frac{1+\varrho}{2}}\|g\|_{C^{1,\sigma}(\partial\Omega)}.
\end{equation}

Set $w^\beta_\varepsilon(x) = u^\beta_\varepsilon(x) - \Phi_{\varepsilon,0}^{\beta\gamma}(x)h_\varepsilon^\gamma$(x), we obtain
\begin{equation*}
 \left\{ \begin{aligned}
 \mathcal{L}_\varepsilon(w_\varepsilon) &= \text{div}(\tilde{f}) + \tilde{F}, &\quad &\text{in}~~\Omega, \\
  w_\varepsilon &= 0, &\quad &\text{on}~\partial\Omega,
 \end{aligned}
 \right.
\end{equation*}
where
\begin{eqnarray*}
\tilde{f}^\alpha & = & A_\varepsilon^{\alpha\beta}(\Phi_{\varepsilon,0}^{\beta\gamma}-\delta^{\beta\gamma})\nabla h_\varepsilon^\gamma
+ V_\varepsilon^{\alpha\beta}(\Phi_{\varepsilon,0}^{\beta\gamma}-\delta^{\beta\gamma})h_\varepsilon^\gamma , \\
\tilde{F}^\alpha & = & A_\varepsilon^{\alpha\beta}\nabla\Phi_{\varepsilon,0}^{\beta\gamma}\nabla h_\varepsilon^\gamma
+ V_\varepsilon^{\alpha\gamma}\nabla h_\varepsilon^\gamma - B_\varepsilon^{\alpha\beta}\nabla(\Phi_{\varepsilon,0}^{\beta\gamma}h_\varepsilon^\gamma)
- (c_\varepsilon^{\alpha\beta}+\lambda\delta^{\alpha\beta})\Phi_{\varepsilon,0}^{\beta\gamma}h_\varepsilon^\gamma.
\end{eqnarray*}
Now, let $\nu = \varrho/2$. In view of $\eqref{pri:5.2}$ and $\eqref{f:5.6}$, we have
\begin{equation*}
 \|\tilde{f}\|_{C^{0,\nu}(\Omega)}\leq C \|g\|_{C^{1,\sigma}(\partial\Omega)} \quad \text{and} \quad
 \|\tilde{F}\|_{L^p(\Omega)} \leq C\|g\|_{C^{1,\sigma}(\partial\Omega)}.
\end{equation*}
Note that $\|h_\varepsilon\|_{L^\infty(\Omega)}\leq C\|g\|_{L^\infty(\partial\Omega)}$, where $C$ depends only on $\mu,\tau,\kappa,\sigma,d,m$ and $\Omega$.
Hence, recalling $\eqref{f:5.10}$, we have
\begin{eqnarray}\label{f:5.7}
\|\nabla u_\varepsilon\|_{L^\infty(\Omega)}
 &\leq &  \|\nabla w_\varepsilon\|_{L^\infty(\Omega)} + \|\nabla (\Phi_{\varepsilon,0}h_\varepsilon)\|_{L^\infty(\Omega)} \nonumber \\
 &\leq &  C\{\|\tilde{f}\|_{C^{0,\nu}(\Omega)} +  \|\tilde{F}\|_{L^p(\Omega)} + \|g\|_{C^{1,\sigma}(\partial\Omega)}\}
 \leq C\|g\|_{C^{1,\sigma}(\partial\Omega)}.
\end{eqnarray}

Finally, $\eqref{pri:1.2}$ follows from $\eqref{f:5.10}$ and $\eqref{f:5.7}$ by writing $u_\varepsilon = u_{\varepsilon,1} + u_{\varepsilon,2}$,
where $u_{\varepsilon,1},u_{\varepsilon,2}$ respectively satisfy the homogeneous and non-homogeneous systems
(see $\eqref{pde:3.1}$). The proof is complete.
\qed

\begin{lemma}\label{lemma:5.5}
Suppose that the coefficients of $\mathcal{L}_\varepsilon$ satisfy the same conditions as in Theorem $\ref{thm:1.5}$.
Then $\mathcal{G}_\varepsilon(x,y)$ has the following estimates:
\begin{eqnarray}\label{pri:4.19}
|\nabla_x\nabla_y \mathcal{G}_\varepsilon(x,y)| \leq \frac{C}{|x-y|^d},
\end{eqnarray}
\begin{eqnarray}\label{pri:4.18}
|\nabla_x \mathcal{G}_\varepsilon(x,y)|\leq \frac{C}{|x-y|^{d-1}}\min\Big\{1,\frac{d_y}{|x-y|}\Big\} \quad\text{and}\quad
|\nabla_y \mathcal{G}_\varepsilon(x,y)|\leq \frac{C}{|x-y|^{d-1}}\min\Big\{1,\frac{d_x}{|x-y|}\Big\}
\end{eqnarray}
for any $x,y\in\Omega$ and $x\not=y$, where $C$ depends only on $\mu,\tau,\kappa,\lambda,d,m,\eta$ and $\Omega$.
\end{lemma}

\begin{pf}
For any $x,y\in \Omega$, let $r= |x-y|$. Due to
\begin{equation}\label{f:5.19}
\mathcal{L}_\varepsilon^*[\mathcal{G}_\varepsilon(\cdot,y)] = 0 \quad \text{in} ~\Omega\setminus B(y,\rho)
\end{equation}
for any $\rho>0$,
it follows from $\eqref{pri:1.2}$ and $\eqref{f:4.11}$ that
\begin{equation}\label{f:5.26}
|\nabla_x \mathcal{G}_\varepsilon(x,y)|
\leq \|\nabla \mathcal{G}_\varepsilon(\cdot,y) \|_{L^\infty(\Omega_{\frac{r}{6}}(x))}
\leq \frac{C}{r}\Big(\dashint_{\Omega_{\frac{r}{3}}(x)}|\mathcal{G}_\varepsilon(z,y)|^2 dz\Big)^{\frac{1}{2}}
\leq Cr^{1-d},
\end{equation}
where $x$ can be on $\partial\Omega$.
By applying the localization technique (as shown in Remark $\ref{rm:2.7}$)
to $\eqref{pri:1.2}$,
we have
\begin{equation*}
\|\nabla u_\varepsilon\|_{L^\infty(\Omega_{\frac{r}{6}}(x))}
\leq \frac{C}{r}\Big(\dashint_{\Omega_{\frac{r}{3}}(x)}|u_\varepsilon|^2 dy\Big)^{1/2}
\end{equation*}
for $u_\varepsilon$ satisfying $\mathcal{L}_\varepsilon(u_\varepsilon) = 0$
in $\Omega_{\frac{r}{3}}(x)$
and $u_\varepsilon = 0$ on $\partial(\Omega_{\frac{r}{3}}(x))\cap\partial\Omega$.
(We remark that we just consider the estimate at boundary,
and the interior one directly follows from $\eqref{pri:5.8}$.)
So, we can derive the second inequality of $\eqref{f:5.26}$.

For the adjoint Green's matrix $^*\mathcal{G}_\varepsilon(\cdot,x)$, we have
$|\nabla_y \mathcal{G}_\varepsilon(x,y)| = |\nabla_y ^*\mathcal{G}_\varepsilon(y,x)| \leq Cr^{1-d}$ by the same argument.
Moreover, since $\nabla_y \mathcal{G}_\varepsilon(\cdot,y)$ still satisfies $\eqref{f:5.19}$ for any $\rho>0$,
and $\nabla_y \mathcal{G}_\varepsilon(\cdot,y) = 0$ on $\partial\Omega$, we obtain
\begin{equation*}
 |\nabla_x \nabla_y\mathcal{G}_\varepsilon(x,y)|
 \leq \frac{C}{r}\Big(\dashint_{\Omega_{\frac{r}{3}}(x)}|z-y|^{2(1-d)} dz\Big)^{\frac{1}{2}} \leq Cr^{-d},
\end{equation*}
where $r/2<|z-y|<2r$. Observe that $\nabla_y \mathcal{G}_\varepsilon(\cdot,y) = 0$ and
$\nabla_x [^*\mathcal{G}_\varepsilon(\cdot,x)] = 0$ on $\partial\Omega$, we have
\begin{equation*}
 |\nabla_y\mathcal{G}_\varepsilon(x,y)| = |\nabla_y\mathcal{G}_\varepsilon(x,y) - \nabla_y\mathcal{G}_\varepsilon(\bar{x},y)|
 \leq |\nabla_x \nabla_y\mathcal{G}_\varepsilon(x,y)| |x-\bar{x}| \leq \frac{Cd_x}{|x-y|^d},
\end{equation*}
where $\bar{x}\in\partial\Omega$ such that $d_x=|x-\bar{x}|$. Similarly, we have $|\nabla_x\mathcal{G}_\varepsilon(x,y)|\leq Cd_y|x-y|^{-d}$,
and the proof is complete.
\qed
\end{pf}

\textbf{Proof of Theorem \ref{thm:1.5}}\textbf{.}\quad Define the conormal derivative $\frac{\partial}{\partial\nu}$,
$\left(\frac{\partial}{\partial\nu}\right)^*$ corresponding to $\mathcal{L}_\varepsilon$ and $\mathcal{L}_\varepsilon^*$ as follows:
\begin{eqnarray*}
\frac{\partial}{\partial\nu}  = -n_i A_{ij}(x/\varepsilon)\frac{\partial}{\partial x_j} - n_i V_i(x/\varepsilon), \qquad\quad
\left(\frac{\partial}{\partial\nu}\right)^* =  -n_j A_{ij}(x/\varepsilon)\frac{\partial}{\partial x_i} - n_j B_j(x/\varepsilon),
\end{eqnarray*}
where $n=(n_1,\cdots,n_d)$ denotes the outward unit normal vector to $\partial\Omega$.
Thus, define the Poisson kernel
$\mathcal{P}_\varepsilon(\cdot,y) = \big[\mathcal{P}_\varepsilon^{\gamma\beta}(\cdot,y)\big]$ associated with
$\mathcal{L}_\varepsilon$ as
\begin{equation*}
 \mathcal{P}_\varepsilon^{\gamma\beta}(x,y) = \left(\frac{\partial}{\partial \nu}\right)^*\big[\mathcal{G}_\varepsilon^{\alpha\gamma}(x,y)\big]
 = -n_j(x)a_{ij}^{\alpha\beta}(x/\varepsilon)\frac{\partial}{\partial x_i}\big\{\mathcal{G}_\varepsilon^{\alpha\gamma}(x,y)\big\}
  -n_j(x)B_{j}^{\alpha\beta}(x/\varepsilon)\mathcal{G}_\varepsilon^{\alpha\gamma}(x,y)
\end{equation*}
for $y\in\Omega$ and $x\in\partial\Omega$.  It follows from $\eqref{pri:4.18}$ that
\begin{equation}\label{f:5.18}
|\mathcal{P}_\varepsilon(x,y)| \leq Cd_y|x-y|^{-d} ,
\end{equation}
where $C$ depends only on $\mu,\tau,\kappa,\lambda,d,m,\eta$ and $\Omega$.
Thus for any $g\in L^p(\partial\Omega;\mathbb{R}^m)$ with $p\in (1,\infty]$, the solution to
$\mathcal{L}_\varepsilon(u_\varepsilon) = 0$ in $\Omega$ and $u_\varepsilon = g$ on $\partial\Omega$
can be written by
\begin{equation}\label{f:5.11}
  u_\varepsilon(y) = \int_{\partial\Omega} \mathcal{P}_\varepsilon(x,y)g(x)dS(x)
\end{equation}
for any $y\in\Omega$, and it follows from $\eqref{f:5.18}$ that
\begin{equation}
|u_\varepsilon(y)| \leq Cd_y\int_{\partial\Omega} \frac{|g(x)|}{|x-y|^d}dS(x)
\end{equation}

Recall that the nontangential maximal function of $u_\varepsilon$ is defined by
\begin{equation}\label{def:3}
(u_\varepsilon)^*(Q) = \sup\big\{ |u_\varepsilon(x)|:x\in\Omega ~\text{and} ~|x-Q|\leq N_0\text{dist}(x,\partial\Omega)\big\}
\qquad \text{for}~ Q\in\partial\Omega,
\end{equation}
where $N_0=N_0(\Omega)>1$ is sufficiently large.
Hence, if $|y-x_0|\leq N_0d_y$ for some $x_0\in\partial\Omega$, then we have
\begin{equation*}
\begin{aligned}
|u_\varepsilon(y)|
&\leq C\dashint_{\partial\Omega\cap B(x_0,r)}|g(x)|dS(x)
+ Cd_y\sum_{j=0}^\infty\int_{\Sigma_j}\frac{|g(x)|dS(x)}{|x-y|^d} \\
& \leq C\bigg\{\dashint_{\partial\Omega\cap B(x_0,r)}|g(x)|dS(x)
+ \sum_{j=0}^{\infty}2^{-j}\dashint_{\partial\Omega\cap B(x_0,2^{j+1}r)}|g(x)|dS(x)\bigg\} \\
&\leq C
\sup_{0 < r < \text{diam}(\Omega)}\Big\{\dashint_{\partial\Omega\cap B(x_0,r)}|g(x)|dS(x)\Big\}
\end{aligned}
\end{equation*}
where $r=d_y/2$ and $\Sigma_j = \partial\Omega\cap\big\{B(x_0,2^{j+1}r)\setminus B(x_0,2^jr)\big\}$.
Note that
\begin{equation*}
\mathcal{M}_{\partial\Omega}(|g|)(x_0) = \sup_{0 < r < \text{diam}(\Omega)}\Big\{\dashint_{\partial\Omega\cap B(x_0,r)}|g(x)|dS(x)\Big\}
\end{equation*}
is the Hardy-Littlewood maximal function of $g$ on $\partial\Omega$. Thus it is not hard to see that
\begin{equation*}
(u_\varepsilon)^*(x_0) \leq C\mathcal{M}_{\partial\Omega}(|g|)(x_0).
\end{equation*}

Due to the $L^p$ bounded properties of the Hardy-Littlewood maximal operator:
$\|\mathcal{M}_{\partial\Omega}(|g|)\|_{L^p(\partial\Omega)}\leq C\|g\|_{L^p(\partial\Omega)}$
for $1<p\leq \infty$ (see \cite[pp.5]{S}),
the estimates $\eqref{pri:1.6}$ and $\eqref{pri:1.7}$
can be derived immediately.

Now, we turn back to verify $\eqref{f:5.11}$.
Let $R = d_y/2$, and $\varphi\in C_0^1(B(y,R))$ be a cut-off function such that $\varphi= 1$ in $B(y,R/4)$ and $\varphi=0$ outside $B(y,R/2)$.
Then, since $\mathcal{L}_\varepsilon(u_\varepsilon) = 0$, we have
$\mathcal{L}_\varepsilon(\varphi u_\varepsilon) = -\mathcal{L}_\varepsilon[(1-\varphi)u_\varepsilon]$ in $\Omega$ and
$\varphi u_\varepsilon = 0$ on $\partial\Omega$. Hence, in view of $\eqref{pri:4.16}$, we obtain
\begin{eqnarray*}
  \big(\varphi u_\varepsilon\big)(y)
 &=& -\int_\Omega\mathcal{G}_\varepsilon(\cdot,y)\mathcal{L}_\varepsilon[(1-\varphi)u_\varepsilon] \\
 &=& -\lim_{r\to 0} \Big\{\int_{\Omega\setminus B(y,r)}\mathcal{G}_\varepsilon(\cdot,y)\mathcal{L}_\varepsilon[(1-\varphi)u_\varepsilon]
 - \int_{\Omega\setminus B(y,r)}\mathcal{L}_\varepsilon^*\big[\mathcal{G}_\varepsilon(\cdot,y)\big][(1-\varphi)u_\varepsilon]  \Big\} \\
 &=& \int_{\partial\Omega}\mathcal{P}_\varepsilon(\cdot,y)[(1-\varphi)u_\varepsilon]
  +\lim_{r\to0} \int_{\partial B(y,r)} \mathcal{G}_\varepsilon(\cdot,y)\frac{\partial}{\partial\nu}\big[(1-\varphi)u_\varepsilon\big]
  - \lim_{r\to0}\int_{\partial B(y,r)} \mathcal{P}_\varepsilon(\cdot,y)[(1-\varphi)u_\varepsilon] \\
 &=& \int_{\partial\Omega}\mathcal{P}_\varepsilon(\cdot,y)[(1-\varphi)u_\varepsilon].
\end{eqnarray*}
Note that $\mathcal{L}_\varepsilon^*\big[\mathcal{G}_\varepsilon(\cdot,y)\big] = 0$ in $\Omega\setminus B(y,r)$ for any $r>0$,
and $(1-\varphi)u_\varepsilon \equiv 0$ in $B(y,R/4)$.
The proof is complete.
\qed

\begin{remark}\label{rm:5.1}
\emph{Note that the same type of results
for $L_\varepsilon$ with Dirichlet boundary conditions and Neumann boundary conditions were shown in
\cite[Theorem 3]{MAFHL5} and \cite[Theorem 1.3]{SZW4}, respectively.
Also, we refer the reader to \cite[Theorem 1.3]{SZ} for the same type of result in the almost periodic setting.
In the case of $m=1$, when we
derive the estimate $\eqref{pri:1.7}$ with $C=1$, there is no regularity condition on the coefficients of $\mathcal{L}_\varepsilon$, but
some additional conditions on $V$ are inevitably required even when $\lambda\geq\lambda_0$, and $\varepsilon =1$ (see \cite[pp.179]{DGNT}). }
\end{remark}

\section{Convergence rates}
\begin{lemma}\label{lemma:6.5}
Suppose that $A\in\Lambda(\mu,\tau,\kappa)$, and $V$ satisfies $\eqref{a:2}$ and $\eqref{a:4}$. Let
\begin{equation*}
\Psi_{\varepsilon,0}^{\alpha\beta}(x) = \Phi_{\varepsilon,0}^{\alpha\beta}(x) - \delta^{\alpha\beta} -\varepsilon\chi_{0}^{\alpha\beta}(x/\varepsilon),\qquad
\Psi_{\varepsilon,k}^{\beta}(x) = \Phi_{\varepsilon,k}^{\beta}(x) - P_k^\beta -\varepsilon\chi_k^\beta(x/\varepsilon).
\end{equation*}
Then we have
\begin{eqnarray}\label{f:6.1}
\max_{0\leq k\leq d}\|\Psi_{\varepsilon,k}\|_{L^\infty(\Omega)}\leq C\varepsilon,
\qquad
\max_{0\leq k\leq d \atop
  x\in\Omega}|\nabla \Psi_{\varepsilon,k}(x)| \leq C\min\big\{1,~\varepsilon d_x^{-1}\big\},
\end{eqnarray}
where $C$ depends only on $\mu,\tau,\kappa,m,d$ and $\Omega$.
\end{lemma}

\begin{pf}
By the definition of Dirichlet correctors $\Phi_{\varepsilon,k}$, we have
$L_\varepsilon(\Psi_{\varepsilon,k}) = 0$ in $\Omega$ and $\Psi_{\varepsilon,k} = -\varepsilon\chi_{k,\varepsilon}$ on $\partial\Omega$.
Thus it follows from  the interior Lipschitz estimate $\eqref{pri:5.8}$ and Agmon-Miranda maximum principle (see \cite{MAFHL}) that,
\begin{equation*}
 |\nabla \Psi_{\varepsilon,k}(x)| \leq \frac{C}{d_x}\Big(\dashint_{B(x,d_x)} |\Psi_{\varepsilon,k}|^2dy\Big)^{1/2}
 \leq C\varepsilon d_x^{-1}, \quad\qquad \forall x\in \Omega.
\end{equation*}
The rest parts of the lemma follow from Lemmas $\ref{lemma:5.2}$ and $\ref{lemma:5.3}$, and Remark $\ref{rm:2.4}$.
\qed
\end{pf}

\begin{lemma}\label{lemma:6.1}
Suppose that $u_\varepsilon \in H^1(\Omega)$, $u\in H^2(\Omega)$ and $\mathcal{L}_\varepsilon(u_\varepsilon) = \mathcal{L}_0(u)$ in $\Omega$. Let
\begin{equation*}
 w_\varepsilon^\beta = u_\varepsilon^\beta - \Phi_{\varepsilon,0}^{\beta\gamma} u^\gamma -
 \big[\Phi_{\varepsilon,k}^{\beta\gamma}-x_k\delta^{\beta\gamma}\big]\frac{\partial u^\gamma}{\partial x_k},
\end{equation*}
where $1\leq k\leq d$. Then
\begin{eqnarray}\label{pri:6.1}
[\mathcal{L}_\varepsilon(w_\varepsilon)]^\alpha
& = & \frac{\partial}{\partial x_i}\left\{\widetilde{E}_{kij,\varepsilon}^{\alpha\gamma}\frac{\partial^2 u^\gamma}{\partial x_k\partial x_j}
+ \widetilde{F}_{ki,\varepsilon}^{\alpha\gamma}\frac{\partial u^\gamma}{\partial x_k}
+ \widetilde{H}_{i,\varepsilon}^{\alpha\gamma}u^\gamma \right\}
+ a_{ij,\varepsilon}^{\alpha\beta}\left\{\frac{\partial \Psi_{\varepsilon,k}^{\beta\gamma}}{\partial x_j}\frac{\partial^2 u^\gamma}{\partial x_k\partial x_i}
+ \frac{\partial \Psi_{\varepsilon,0}^{\beta\gamma}}{\partial x_j}\frac{\partial u^\gamma}{\partial x_i}\right\}   \nonumber \\
& - & B_{i,\varepsilon}^{\alpha\beta}\left\{\big[\Phi_{\varepsilon,k}^{\beta\gamma}-x_k\delta^{\beta\gamma}\big]\frac{\partial^2 u^\gamma}{\partial x_i\partial x_k}
+ \big[\Phi_{\varepsilon,0}^{\beta\gamma}-\delta^{\beta\gamma}\big]\frac{\partial u^\gamma}{\partial x_i} \right\}
- B_{i,\varepsilon}^{\alpha\beta}\left\{\frac{\partial \Psi_{\varepsilon,k}^{\beta\gamma}}{\partial x_i}\frac{\partial u^\gamma}{\partial x_k}
+ \frac{\partial \Psi_{\varepsilon,0}^{\beta\gamma}}{\partial x_i}u^\gamma \right\}   \nonumber \\
& - & \big[c^{\alpha\beta}_\varepsilon + \lambda\delta^{\alpha\beta}\big]\left\{
\big[\Phi_{\varepsilon,k}^{\beta\gamma}-x_k\delta^{\beta\gamma}\big]\frac{\partial u^\gamma}{\partial x_k}
+\big[\Phi_{\varepsilon,0}^{\beta\gamma}-\delta^{\beta\gamma}\big]u^\gamma  \right\}
- \varepsilon\frac{\partial\vartheta_k^{\alpha\gamma}}{\partial y_i}\frac{\partial^2 u^\gamma}{\partial x_k\partial x_i}
- \varepsilon\frac{\partial\zeta^{\alpha\gamma}}{\partial y_i}\frac{\partial u^\gamma}{\partial x_i},
\end{eqnarray}
where $y=x/\varepsilon$, and
\begin{eqnarray*}
\widetilde{E}_{kij,\varepsilon}^{\alpha\gamma} & = & \varepsilon E_{kij,\varepsilon}^{\alpha\gamma}
+ a_{ij,\varepsilon}^{\alpha\beta}\big[\Phi_{\varepsilon,k}^{\beta\gamma}-x_k\delta^{\beta\gamma}\big], \\
\widetilde{F}_{ki,\varepsilon}^{\alpha\gamma} ~& = & \varepsilon F_{ki,\varepsilon}^{\alpha\gamma}
+ a_{ik,\varepsilon}^{\alpha\beta}\big[\Phi_{\varepsilon,0}^{\beta\gamma}-\delta^{\beta\gamma}\big]
+ V_{i,\varepsilon}^{\alpha\beta}\big[\Phi_{\varepsilon,k}^{\beta\gamma}-x_k\delta^{\beta\gamma}\big]
+ \varepsilon\frac{~\partial\vartheta_k^{\alpha\gamma}}{\partial y_i}, \\
\widetilde{H}_{i,\varepsilon}^{\alpha\gamma}~ & = &  V_{i,\varepsilon}^{\alpha\beta}\big[\Phi_{\varepsilon,0}^{\beta\gamma}-\delta^{\beta\gamma}\big]
+ \varepsilon\frac{~\partial\zeta^{\alpha\gamma}}{\partial y_i}.
\end{eqnarray*}
Note that $E_{kij,\varepsilon}^{\alpha\gamma}, F_{ki,\varepsilon}^{\alpha\gamma}, \vartheta_{k}^{\alpha\gamma}, \zeta^{\alpha\gamma}$
are defined in Remark $\ref{rm:2.4}$.
\end{lemma}

\begin{pf}
From $\mathcal{L}_\varepsilon(u_\varepsilon) = \mathcal{L}_0(u)$, it follows that
\begin{eqnarray}{\label{f:6.3}}
 [\mathcal{L}_\varepsilon(w_\varepsilon)]^\alpha & = &  [\mathcal{L}_0(u)]^\alpha - [\mathcal{L}_\varepsilon(\Phi_{\varepsilon,0}u)]^\alpha
 - \big[\mathcal{L}_\varepsilon\big((\Phi_{\varepsilon,k}-x_kI)\frac{\partial u}{\partial x_k}\big)\big]^\alpha \nonumber\\
 & = & -\frac{\partial}{\partial x_i}\left\{\widehat{a}_{ij}^{\alpha\beta}\frac{\partial u^\beta}{\partial x_j}\right\}
 + \big(\widehat{B}_i^{\alpha\beta} - \widehat{V}_i^{\alpha\beta}\big)\frac{\partial u^\beta}{\partial x_i}
 + \big[\widehat{c}^{\alpha\beta}+\lambda\delta^{\alpha\beta}\big]u^\beta
 - [I_1]^\alpha - [I_2]^\alpha,
\end{eqnarray}
where $I_1 = \mathcal{L}_\varepsilon(\Phi_{\varepsilon,0}u)$ and
$I_2 =  \mathcal{L}_\varepsilon\big[(\Phi_{\varepsilon,k}-x_kI)\frac{\partial u}{\partial x_k}\big]$.
By the definition of $\Phi_{\varepsilon,k}$, and $\chi_{k}$, $0\leq k\leq d$,  we obtain
\begin{eqnarray*}
[I_1]^\alpha
& = & -\frac{\partial}{\partial x_i}\Big\{a_{ij,\varepsilon}^{\alpha\beta}\frac{\partial u^\beta}{\partial x_j}\Big\}
-\frac{\partial}{\partial x_i}\Big\{a_{ij,\varepsilon}^{\alpha\beta}\big[\Phi_{\varepsilon,0}^{\beta\gamma}-\delta^{\beta\gamma}\big]\frac{\partial u^\gamma}{\partial x_j}\Big\}
- a_{ij,\varepsilon}^{\alpha\beta}\frac{\partial\Phi_{\varepsilon,0}^{\beta\gamma}}{\partial x_j}\frac{\partial u^\gamma}{\partial x_i} \\
& - & \frac{\partial}{\partial x_i}\big\{V_{i,\varepsilon}^{\alpha\beta}\big\}\big[\Phi_{\varepsilon,0}^{\beta\gamma}-\delta^{\beta\gamma}\big]u^\gamma
- V_{i,\varepsilon}^{\alpha\beta}\frac{\partial}{\partial x_i}\big\{\Phi_{\varepsilon,0}^{\beta\gamma}u^\gamma\big\}
+ B_{i,\varepsilon}^{\alpha\beta}\frac{\partial}{\partial x_i}\big\{\Phi_{\varepsilon,0}^{\beta\gamma}u^\gamma\big\}
+ \big[c_\varepsilon^{\alpha\beta}+\lambda\delta^{\alpha\beta}\big]\Phi_{\varepsilon,0}^{\beta\gamma}u^\gamma \\
&=& -\frac{\partial}{\partial x_i}\Big\{a_{ij,\varepsilon}^{\alpha\beta}\frac{\partial u^\beta}{\partial x_j}\Big\}
- \frac{\partial}{\partial x_i}\Big\{a_{ij,\varepsilon}^{\alpha\beta}\big(\Phi_{\varepsilon,0}^{\beta\gamma}-\delta^{\beta\gamma}\big)\frac{\partial u^\gamma}{\partial x_j}\Big\}
- a_{ij,\varepsilon}^{\alpha\beta}\frac{\partial\Psi_{\varepsilon,0}^{\beta\gamma}}{\partial x_j}\frac{\partial u^\gamma}{\partial x_i}
- a_{ij,\varepsilon}^{\alpha\beta}\frac{\partial\chi_0^{\beta\gamma}}{\partial y_j}\frac{\partial u^\gamma}{\partial x_i}   \\
&-&  \frac{\partial}{\partial x_i}\Big\{V_{i,\varepsilon}^{\alpha\beta}\big[\Phi_{\varepsilon,0}^{\beta\gamma}-\delta^{\beta\gamma}\big]u^\gamma\Big\}
- V_{i,\varepsilon}^{\alpha\beta}\frac{\partial u^\beta}{\partial x_i}
+ B_{i,\varepsilon}^{\alpha\beta}\Big\{\frac{\partial u^\beta}{\partial x_i}
+ \big[\Phi_{\varepsilon,0}^{\beta\gamma}-\delta^{\beta\gamma}\big]\frac{\partial u^\gamma}{\partial x_i}
+ \frac{\partial\Psi_{\varepsilon,0}^{\beta\gamma}}{\partial x_i}u^\gamma
+ \frac{\partial\chi_0^{\beta\gamma}}{\partial y_i}u^\gamma\Big\} \\
&+& \big[c_\varepsilon^{\alpha\beta}+\lambda\delta^{\alpha\beta}\big]\Phi_{\varepsilon,0}^{\beta\gamma}u^\gamma
\end{eqnarray*}
and
\begin{eqnarray*}
[I_2]^\alpha
& = & -\frac{\partial}{\partial x_i}\Big\{a_{ij,\varepsilon}^{\alpha\beta}\frac{\partial\chi_k^{\beta\gamma}}{\partial y_j}\frac{\partial u^\gamma}{\partial x_k}\Big\}
-\frac{\partial}{\partial x_i}\Big\{a_{ij,\varepsilon}^{\alpha\beta}\frac{\partial}{\partial x_j}\big[\Phi_{\varepsilon,k}^{\beta\gamma}-x_k\delta^{\beta\gamma}-\varepsilon\chi_{k,\varepsilon}^{\beta\gamma}\big]
\frac{\partial u^\gamma}{\partial x_k}\Big\}
- \frac{\partial}{\partial x_i}\Big\{V_{i,\varepsilon}^{\alpha\beta}\big[\Phi_{\varepsilon,k}^{\beta\gamma}-x_k\delta^{\beta\gamma}\big]\frac{\partial u^\gamma}{\partial x_k}\Big\} \\
&-& \frac{\partial}{\partial x_i}\Big\{a_{ij,\varepsilon}^{\alpha\beta}\big[\Phi_{\varepsilon,k}^{\beta\gamma}-x_k\delta^{\beta\gamma}\big]\frac{\partial^2 u^\gamma}{\partial x_j\partial x_k}\Big\}
+ B_{i,\varepsilon}^{\alpha\beta}\frac{\partial}{\partial x_i}\Big\{\big[\Phi_{\varepsilon,k}^{\beta\gamma}-x_k\delta^{\beta\gamma}\big]\frac{\partial u^\gamma}{\partial x_k}\Big\}
+ \big[c_\varepsilon^{\alpha\beta}+\lambda\delta^{\alpha\beta}\big]\big[\Phi_{\varepsilon,k}^{\alpha\beta}-x_k\delta^{\alpha\beta}\big]\frac{\partial u^\gamma}{\partial x_k} \\
& = & -\frac{\partial}{\partial x_i}\Big\{a_{ik,\varepsilon}^{\alpha\beta}\frac{\partial\chi_j^{\beta\gamma}}{\partial y_k}\frac{\partial u^\gamma}{\partial x_j}\Big\}
- a_{ij,\varepsilon}^{\alpha\beta}\frac{\partial\Psi_{\varepsilon,k}^{\beta\gamma}}{\partial x_j}\frac{\partial^2 u^\gamma}{\partial x_k\partial x_i}
- \frac{\partial}{\partial x_i}\Big\{a_{ij,\varepsilon}^{\alpha\beta}\big[\Phi_{\varepsilon,k}^{\beta\gamma}-x_k\delta^{\beta\gamma}\big]\frac{\partial^2 u^\gamma}{\partial x_j\partial x_k}\Big\} \\
& - & \frac{\partial}{\partial x_i}\Big\{V_{i,\varepsilon}^{\alpha\beta}\big[\Phi_{\varepsilon,k}^{\beta\gamma}-x_k\delta^{\beta\gamma}\big]\frac{\partial u^\gamma}{\partial x_k}\Big\}
+ B_{i,\varepsilon}^{\alpha\beta}\Big\{\frac{\partial\chi_k^{\beta\gamma}}{\partial y_i}\frac{\partial u^\gamma}{\partial x_k}
+ \frac{\partial\Psi_{\varepsilon,k}^{\beta\gamma}}{\partial x_i}\frac{\partial u^\gamma}{\partial x_k}
+ \big[\Phi_{\varepsilon,k}^{\beta\gamma} - x_k\delta^{\beta\gamma}\big]\frac{\partial^2 u^\gamma}{\partial x_i\partial x_k}\Big\} \\
& + & \big[c_\varepsilon^{\alpha\beta}+\lambda\delta^{\alpha\beta}\big]\big[\Phi_{\varepsilon,k}^{\alpha\beta}-x_k\delta^{\alpha\beta}\big]\frac{\partial u^\gamma}{\partial x_k},
\end{eqnarray*}
where $y= x/\varepsilon$. Put $I_1$ and $I_2$ into $\eqref{f:6.3}$, and then we have
\begin{eqnarray}\label{f:6.6}
[\mathcal{L}_\varepsilon(w_\varepsilon)]^\alpha
& = & -\frac{\partial}{\partial x_i}\Big\{b_{ij}^{\alpha\gamma}(y)\frac{\partial u^\gamma}{\partial x_j}\Big\}
+ \frac{\partial}{\partial x_i}\Big\{ a_{ij,\varepsilon}^{\alpha\beta}\big[\Phi_{\varepsilon,k}^{\beta\gamma}-x_k\delta^{\beta\gamma}\big]\frac{\partial^2 u^\gamma}{\partial x_k\partial x_j}
+ a_{ij,\varepsilon}^{\alpha\beta}\big[\Phi_{\varepsilon,0}^{\beta\gamma}-\delta^{\beta\gamma}\big]\frac{\partial u^\gamma}{\partial x_j} \Big\}  \nonumber\\
& + &  a_{ij,\varepsilon}^{\alpha\beta}\Big\{\frac{\partial \Psi_{\varepsilon,k}^{\beta\gamma}}{\partial x_j}\frac{\partial^2 u^\gamma}{\partial x_k\partial x_i}
+ \frac{\partial \Psi_{\varepsilon,0}^{\beta\gamma}}{\partial x_j}\frac{\partial u^\gamma}{\partial x_i}\Big\}
+ \frac{\partial}{\partial x_i}\Big\{V_{i,\varepsilon}^{\alpha\beta}\big[\Phi_{\varepsilon,k}^{\beta\gamma}-x_k\delta^{\beta\gamma}\big]\frac{\partial u^\gamma}{\partial x_k}
+ V_{i,\varepsilon}^{\alpha\beta}\big[\Phi_{\varepsilon,0}^{\beta\gamma}-\delta^{\beta\gamma}\big]u^\gamma\Big\}  \nonumber\\
& - & B_{i,\varepsilon}^{\alpha\beta}\left\{\big[\Phi_{\varepsilon,k}^{\beta\gamma}-x_k\delta^{\beta\gamma}\big]\frac{\partial^2 u^\gamma}{\partial x_i\partial x_k}
+ \big[\Phi_{\varepsilon,0}^{\beta\gamma}-\delta^{\beta\gamma}\big]\frac{\partial u^\gamma}{\partial x_i}
+ \frac{\partial \Psi_{\varepsilon,k}^{\beta\gamma}}{\partial x_i}\frac{\partial u^\gamma}{\partial x_k}
+ \frac{\partial \Psi_{\varepsilon,0}^{\beta\gamma}}{\partial x_i}u^\gamma \right\}  \nonumber\\
& -& \big[c^{\alpha\beta}_\varepsilon + \lambda\delta^{\alpha\beta}\big]\Big\{\big[\Phi_{\varepsilon,0}^{\beta\gamma}-\delta^{\beta\gamma}\big]u^\gamma + \big[\Phi_{\varepsilon,k}^{\beta\gamma}-x_k\delta^{\beta\gamma}\big]\frac{\partial u^\gamma}{\partial x_k}\Big\}
-U_{i}^{\alpha\gamma}(y)\frac{\partial u^\gamma}{\partial x_i}
 +  W_{k,\varepsilon}^{\alpha\gamma}\frac{\partial u^\gamma}{\partial x_k}
 +  Z_\varepsilon^{\alpha\gamma} u^\gamma,
\end{eqnarray}
where $b_{ij}^{\alpha\gamma}$,  $U_i^{\alpha\gamma}$, $W_{i,\varepsilon}^{\alpha\gamma}$, $Z_\varepsilon^{\alpha\gamma}$ are defined in
Remark $\ref{rm:2.4}$. Besides, the following identities hold.
\begin{eqnarray*}
-\frac{\partial}{\partial x_i}\left\{b_{ij}^{\alpha\gamma}(y)\frac{\partial u^\gamma}{\partial x_j}\right\}
&=& -\varepsilon\frac{\partial}{\partial x_i}\left\{\frac{\partial}{\partial x_k}\big[E_{kij,\varepsilon}^{\alpha\gamma}\big]\frac{\partial u^\gamma}{\partial x_j}\right\}\\
&=& -\varepsilon\frac{\partial}{\partial x_i}\left\{\frac{\partial}{\partial x_k}\big[E_{kij,\varepsilon}^{\alpha\gamma}\frac{\partial u^\gamma}{\partial x_j}\big]\right\}
+\varepsilon\frac{\partial}{\partial x_i}\left\{E_{kij,\varepsilon}^{\alpha\gamma}\frac{\partial^2 u^\gamma}{\partial x_k\partial x_j}\right\}
= \varepsilon\frac{\partial}{\partial x_i}\left\{E_{kij,\varepsilon}^{\alpha\gamma}\frac{\partial^2 u^\gamma}{\partial x_k\partial x_j}\right\} ;\\
U_{i}^{\alpha\gamma}(y)\frac{\partial u^\gamma}{\partial x_i}
&=& \varepsilon\frac{\partial}{\partial x_k}\big\{F_{ki,\varepsilon}^{\alpha\gamma}\big\}\frac{\partial u^\gamma}{\partial x_i}
= \varepsilon\frac{\partial}{\partial x_k}\left\{F_{ki,\varepsilon}^{\alpha\gamma}\frac{\partial u^\gamma}{\partial x_i}\right\}
-\varepsilon F_{ki,\varepsilon}^{\alpha\gamma}\frac{\partial^2 u^\gamma}{\partial x_k\partial x_i}
= \varepsilon\frac{\partial}{\partial x_k}\left\{F_{ki,\varepsilon}^{\alpha\gamma}\frac{\partial u^\gamma}{\partial x_i}\right\};\\
W_{k,\varepsilon}^{\alpha\gamma}\frac{\partial u^\gamma}{\partial x_k}
&=& \varepsilon\frac{\partial}{\partial x_i}\left\{\frac{\partial\vartheta_k^{\alpha\gamma}}{\partial y_i}\frac{\partial u^\gamma}{\partial x_k}\right\}
-\varepsilon\frac{~\partial\vartheta_k^{\alpha\gamma}}{\partial y_i}\frac{\partial^2 u^\gamma}{\partial x_k\partial x_i}; \\
Z_\varepsilon^{\alpha\gamma}u^\gamma
&=&\varepsilon\frac{\partial}{\partial x_i}\left\{\frac{\partial\zeta^{\alpha\gamma}}{\partial y_i} u^\gamma\right\}
- \varepsilon\frac{~\partial\zeta^{\alpha\gamma}}{\partial x_i}\frac{\partial u^\gamma}{\partial x_i}.
\end{eqnarray*}
These together with $\eqref{f:6.6}$ give
the formula $\eqref{pri:6.1}$, and we complete the proof.
\qed
\end{pf}


\begin{flushleft}
\textbf{Proof of Theorem \ref{thm:1.3}}\textbf{.}\quad
 Let $w_{\varepsilon,1},w_{\varepsilon,2}\in H^1_0(\Omega;\mathbb{R}^m)$
 satisfy $w_\varepsilon^\beta = w_{\varepsilon,1}^\beta + w_{\varepsilon,2}^\beta$, such that
 \begin{equation}\label{f:6.7}
 \mathcal{L}_\varepsilon(w_{\varepsilon,1}) = \mathcal{L}_\varepsilon(w_\varepsilon)-\Theta
 \quad \text{in}~~ \Omega,
 \qquad
 \mathcal{L}_\varepsilon(w_{\varepsilon,2}) = \Theta \quad \text{in}~~ \Omega,
 \end{equation}
 where $w_\varepsilon$ is given in Lemma $\ref{lemma:6.1}$, and $\Theta = (\Theta^\alpha)$ satisfies
\end{flushleft}
 \begin{equation*}
  \Theta^\alpha = a_{ij,\varepsilon}^{\alpha\beta}\left\{\frac{\partial \Psi_{\varepsilon,k}^{\beta\gamma}}{\partial x_j}\frac{\partial^2 u^\gamma}{\partial x_k\partial x_i}
+ \frac{\partial \Psi_{\varepsilon,0}^{\beta\gamma}}{\partial x_j}\frac{\partial u^\gamma}{\partial x_i}\right\}
- B_{i,\varepsilon}^{\alpha\beta}\left\{
\frac{\partial \Psi_{\varepsilon,k}^{\beta\gamma}}{\partial x_i}\frac{\partial u^\gamma}{\partial x_k}
+ \frac{\partial \Psi_{\varepsilon,0}^{\beta\gamma}}{\partial x_i}u^\gamma \right\}.
 \end{equation*}
 For the first equation of $\eqref{f:6.7}$, it immediately follows from Theorem $\ref{thm:1.1}$, Lemmas $\ref{lemma:5.2}$, $\ref{lemma:5.3}$, and Remark $\ref{rm:2.4}$ that
 \begin{equation}\label{f:22}
  \|\nabla w_{\varepsilon,1}\|_{L^p(\Omega)} \leq C\varepsilon\|u\|_{W^{2,p}(\Omega)}.
 \end{equation}
 For the second equation, in view of $\eqref{pri:2.1}$, we have
 \begin{eqnarray*}
  c_0\|\nabla w_{\varepsilon,2}\|_{L^2(\Omega)}^2 \leq \mathrm{B}_\varepsilon[w_{\varepsilon,2},w_{\varepsilon,2}]
  = \int_{\Omega} \Theta^\alpha w_{\varepsilon,2}^\alpha.
 \end{eqnarray*}
For the right hand side, it follows from $\eqref{f:6.1}$ and Cauchy's inequality that
\begin{eqnarray*}
 \int_{\Omega} \Theta^\alpha w_{\varepsilon,2}^\alpha \leq C\varepsilon\int_\Omega \big(|\nabla^2 u| + |\nabla u| + |u|\big) |w_{\varepsilon,2}|~\frac{dx}{d_x}
 \leq C\varepsilon\|u\|_{H^2(\Omega)} \|\nabla w_{\varepsilon,2}\|_{L^2(\Omega)},
\end{eqnarray*}
where we use Hardy's inequality in the last inequality.
Hence we have
$$\|w_{\varepsilon}\|_{H^1_0(\Omega)}
\leq C\big\{\|\nabla w_{\varepsilon,1}\|_{L^2(\Omega)}
+ \|\nabla w_{\varepsilon,2}\|_{L^2(\Omega)}\big\}\leq C\varepsilon \|u\|_{H^2(\Omega)},$$
and this gives the estimate $\eqref{pri:1.3}$.

We now turn to show the estimate $\eqref{pri:1.7}$. First, by recalling the estimate $\eqref{f:4.2}$, we have
\begin{equation}\label{pri:6.2}
|\mathcal{G}_\varepsilon(x,y)| \leq \frac{Cd_y}{|x-y|^{d-1}}.
\end{equation}
Additionally, in view of $\eqref{f:6.1}$, we have
\begin{equation}\label{pri:6.3}
|\Theta(y)|\leq \frac{C\varepsilon}{d_y}\big\{|\nabla^2 u(y)|+|\nabla u(y)|+ |u(y)|\big\}
\end{equation}
for any $y\in\Omega$.

Since
\begin{eqnarray*}
w_{\varepsilon,2}^\beta(x) = \int_\Omega \mathcal{G}_\varepsilon(x,y)^{\alpha\beta}\Theta^\alpha(y) dy
\qquad \text{for}~ x\in\Omega,
\end{eqnarray*}
it follows from $\eqref{pri:6.2}$ and $\eqref{pri:6.3}$ that
\begin{equation}
|w_{\varepsilon,2}(x)|\leq C\varepsilon\int_\Omega \frac{1}{|x-y|^{d-1}}
\big\{|\nabla^2 u(y)| + |\nabla u(y)| + |u(y)|\big\} dy.
\end{equation}

Thus by the Hardy-Littlewood-Sobolev theorem of fractional integration (see \cite[pp.119]{S}), we obtain
\begin{equation*}
\|w_{\varepsilon,2}\|_{L^q(\Omega)} \leq C\varepsilon\|u\|_{W^{2,p}(\Omega)}
\end{equation*}
with $1/q=1/p-1/d$ when $1<p<d$. For $p>d$, we can straightforward use H\"older's inequality to arrive at
\begin{eqnarray*}
\|w_{\varepsilon,2}\|_{L^\infty(\Omega)} \leq C\varepsilon\|u\|_{W^{2,p}(\Omega)}.
\end{eqnarray*}

Besides, it follows from $\eqref{f:22}$ and the Sobolev inequality that
$\|w_{\varepsilon,1}\|_{L^q(\Omega)}\leq C\|\nabla w_{\varepsilon,1}\|_{L^p(\Omega)}\leq C\varepsilon\|u\|_{W^{2,p}(\Omega)}$. Thus we obtain
\begin{equation*}
 \|w_\varepsilon\|_{L^q(\Omega)} \leq \big\{\|w_{\varepsilon,1}\|_{L^q(\Omega)} + \|w_{\varepsilon,2}\|_{L^q(\Omega)}\big\}
 \leq C\varepsilon\|u\|_{W^{2,p}(\Omega)},
\end{equation*}
which implies the estimate $\eqref{pri:1.4}$, and the proof is complete.
\qed

\begin{remark}\label{rm:6.1}
\emph{ In view of Lemma $\ref{lemma:5.5}$ and Lemma $\ref{lemma:6.1}$,
we can actually derive the following estimates by the arguments developed in \cite{SZW1},
 \begin{equation*}
 \Big|\mathcal{G}_\varepsilon(x,y)- \mathcal{G}_0(x,y)\Big|\leq \frac{C\varepsilon}{|x-y|^{d-1}},
 \qquad \quad \forall~x,y\in \Omega ~~\text{and}~x\not=y.
 \end{equation*}
 Then we have
 \begin{equation*}
  \|u_\varepsilon - u\|_{L^\infty(\Omega)} \leq C\varepsilon\big[\ln(R_0\varepsilon^{-1}+2)\big]^{1-\frac{1}{d}}\|F\|_{L^d(\Omega)},
 \end{equation*}
 where $R_0$ denotes the diameter of $\Omega$. Moreover, let $\Omega$ be a bounded $C^{2,\eta}$ domain, we have
 \begin{equation*}
 \Big|\frac{\partial}{\partial x_i}\big\{\mathcal{G}_\varepsilon^{\alpha\beta}(x,y)\big\}
 -\frac{\partial}{\partial x_i}\big\{\Phi_{\varepsilon,0}^{\alpha\gamma}(x)\big\}\mathcal{G}_0^{\gamma\beta}(x,y)
 -\frac{\partial}{\partial x_i}\big\{\Phi_{\varepsilon,k}^{\alpha\gamma}(x)\big\}\frac{\partial}{\partial x_k}\big\{\mathcal{G}_0^{\alpha\gamma}(x,y)\big\}\Big|
 \leq \frac{C\varepsilon\ln(\varepsilon^{-1}|x-y|+2)}{|x-y|^d}
 \end{equation*}
 for any $x,y\in \Omega$ and $x\not=y$. Then it follows that for any $1<p<\infty$,
 \begin{equation*}
  \big\|u_\varepsilon - \Phi_{\varepsilon,0}u - (\Phi_{\varepsilon,k}^\beta - P_k^\beta)\frac{\partial u^\beta}{\partial x_k}\big\|_{W^{1,p}_0(\Omega)}
\leq C\varepsilon \big[\ln(R_0\varepsilon^{-1}+2)\big]^{2|\frac{1}{p}-\frac{1}{2}|}\|F\|_{L^p(\Omega)},
 \end{equation*}
 where $C$ depends only on $\mu,\tau,\kappa,\lambda,m,d,p$ and $\Omega$. The details are left to readers
 (or see \cite{SZW1}). }
\end{remark}

\begin{center}
\textbf{Acknowledgement }
\end{center}

The author is grateful to Professor Zhongwei Shen for warm hospitality and continued guidance in homogenization theory.
This paper can not be accomplished without his supervision.
The author also wishes to express his sincere appreciation to the Department of Mathematics in University of Kentucky for supplying him with the very good study facilities
when he visited there during 2013-2015. Besides, the author is indebted to the referees for many useful comments and suggestions.
This work was supported by the Chinese Scholar Council (File No. 201306180043),
and in part by the National Natural Science Foundation of China (Grant No. 11471147).